\documentclass[10pt,a4paper]{article}
\usepackage[utf8]{inputenc}
\usepackage[T1]{fontenc}
\usepackage[english]{babel}
\usepackage[bottom=3cm, top=3cm, left=2cm, right=2cm]{geometry}
\usepackage{amsmath}
\usepackage{amsfonts}
\usepackage{amssymb}
\usepackage{amsthm}
\usepackage{graphicx}
\usepackage{color}
\usepackage{mathrsfs}
\usepackage{lmodern}
\usepackage[most]{tcolorbox}
\usepackage{framed}
\usepackage{fancybox}
\usepackage{tikz}
\usepackage{bbm}
\usepackage{dsfont}
\usepackage{wrapfig}
\usepackage{hyperref}
\usepackage[pagewise]{lineno}
\newtcolorbox{boxrd}[2][]{enhanced,colback=white,width={\textwidth},
attach boxed title to top left={yshift={-0.5\baselineskip},xshift=1cm}, 
title={#2},
boxrule=0.5pt,
coltitle=black,
boxed title style={
  borderline={-0.5mm}{black}
  colframe=white,
  colback=gray!50,
  colupper={black},
},
}

\newtcolorbox{boxsp}[2][]{%
  enhanced,colback=white,colframe=black,coltitle=black,
  boxrule=0.4pt,
  fonttitle=\itshape,
  attach boxed title to top left={yshift=-0.5\baselineskip-0.3pt,xshift=2mm},
  boxed title style={tile,size=minimal,left=0.5mm,right=0.5mm,
    colback=white,before upper=\strut},
  title=#2,#1
}

\newcommand{\N}{\mathbb{N}}
\newcommand{\T}{\top}

\newcommand{\R}{\mathbb{R}}
\newcommand{\C}{\mathbb{C}}

\usepackage{fancyhdr}
\usepackage{changepage}

\title{\textbf{Asymptotic stability of solitary waves for the 1D near-cubic Schrödinger equation in the presence of an internal mode}}
\author{Guillaume Rialland}
\date{\footnotesize{Université de Paris-Saclay, UVSQ, CNRS, Laboratoire de Mathématiques de Versailles, 78000 Versailles \\ \texttt{guillaume.rialland@uvsq.fr}}}

\begin{document}
\maketitle

\begin{adjustwidth}{80pt}{80pt}
\small{\textsc{Abstract.} We consider perturbations of the one-dimensional cubic Schrödinger equation, of the form $i \, \partial_t \psi + \partial_x^2 \psi + |\psi|^2 \psi + g( |\psi|^2 ) \psi = 0$. Under hypotheses on the function $g$ that can be easily verified in some cases (such as $g(s) = s^\sigma$ with $\sigma >1$), we show that the linearized problem around a small solitary wave presents a unique internal mode. Moreover, under an additional hypothesis (the Fermi golden rule) that can also be verified in the case of powers $g(s) = s^\sigma$, we prove the asymptotic stability of the solitary waves with small frequencies.}
\end{adjustwidth}

\textcolor{white}{a} \\ \\ \textcolor{white}{a}

\section{Introduction}
\noindent We consider the non-linear Schrödinger equation
\begin{equation}
    i \, \partial_t \psi + \partial_x^2 \psi + |\psi|^2 \psi + g( |\psi|^2 ) \psi = 0, \ \, \, \, \, \, \, \, \, (t \, , x) \in \R \times \R,
    \label{NLS}
\end{equation}
which is a perturbation of the cubic NLS equation $i \, \partial_t \psi + \partial_x^2 \psi + | \psi |^2 \psi = 0$. Here, $g \, : \, \R_+ \to \R$ is a function so that the term $g ( |\psi|^2 ) \psi$ is small compared to $| \psi|^2 \psi$ for $| \psi |$ small. We refer to \cite{Pe} or \cite{Ki} for the physical interest of such equations; it is a classical and important matter to perturb the Schrödinger equation near the cubic non-linearity, and here we study the semi-linear perturbations of that equation. \\
\\ The corresponding Cauchy problem is globally well-posed in the energy space $H^1 ( \R )$ (see for example \cite{Ca2}). We recall that, for any solution $\psi \in H^1 ( \R )$, as long as it exists, the mass, momentum and energy are conserved:
\[ \int_{\R} | \psi |^2 , \, \, \, \, \text{Im} \int_{\R} \psi \, \overline{\partial_y \psi} , \, \, \, \, \int_{\R} \left ( \frac{1}{2} | \partial_x \psi |^2 - \frac{| \psi |^4}{4} - \frac{G ( | \psi |^2 )}{2} \right ), \]
where $G(s) := \int_0^s g$. We also recall the Galilean transform, translation and phase invariances of this equation: if $\psi (t \, , x)$ is a solution then, for any $\beta,\sigma,\gamma \in \R$, $\widetilde{\psi} (t \, , x) = e^{i( \beta x - \beta^2 t + \gamma )} \psi (t \, , x-2 \beta t - \sigma )$ is also a solution to the same equation. \\
\\ Solitary waves are solutions of \eqref{NLS} which take the form $\psi (t \, , x) = e^{i \omega t} \phi_\omega (x)$ where 
\begin{equation}
    \phi_\omega '' = \omega \phi_\omega - \phi_\omega^3 - \phi_\omega g(\phi_\omega^2). 
    \label{eqphi}
\end{equation}
Below we introduce the first elementary hypothesis:
\[ \begin{array}{ccl} (H_1) & : & g \in \mathscr{C}^5 ((0 \, , + \infty )) \cap \mathscr{C}^1 ( [0 \, , + \infty )) \, , \, \, g^{(k)}(s) \, \underset{s \to 0^+}{=} \, o \left ( s^{1-k} \right ) \, \, \text{for all $k \in \{ 0 \, , 1 \, , 2 \, , 3 \, , 4 \}$,} \\ & & \text{$g^{(5)} (s) \, \underset{s \to 0}{=} \, \mathcal{O} \left ( s^{-4} \right )$ and $g \not\equiv 0$ near $0$.} \end{array} \]
Assuming hypothesis $(H_1)$ holds and provided $\omega >0$ is small enough, the equation \eqref{eqphi} has a unique solution $\phi_\omega \in H^1 ( \R )$ that is positive and even. The invariances previously described generate a family of traveling waves given by $\psi (t \, , x) = e^{i ( \beta x - \beta^2 t + \omega t + \gamma )} \phi_\omega (x-2 \beta t - \sigma )$. To begin with, we recall the following standard orbital stability result (see \cite{Ca}, \cite{Gr}, \cite{Il}, \cite{We2}).

\begin{leftbar}
\noindent \textbf{Proposition 1.} For $\omega_0$ small enough and any $\epsilon > 0$, there exists $\delta > 0$ so that, for any $\psi_0 \in H^1 ( \R )$ satisfying $|| \psi_0 - \phi_{\omega_0} ||_{H^1 ( \R )} \leqslant \delta$, if $\psi$ is the solution of \eqref{NLS} with initial data $\psi (0) = \psi_0$, then
\[ \sup_{t \in \R} \inf_{( \gamma , \sigma ) \in \R^2} || \psi (t \, , \cdot + \sigma ) - e^{i \gamma} \phi_{\omega_0} ||_{H^1 ( \R )} \leqslant \epsilon.  \]
\end{leftbar}

\noindent The present paper establishes a result of asymptotic stability of small solitary waves for the equation \eqref{NLS}, under hypotheses that will be presented further. A vast literature deals with the asymptotic stability of solitary waves for nonlinear Schrödinger equations, in different cases (various nonlinearities, with or without potential, in different dimensions), see for example \cite{Co}, \cite{Cu1}, \cite{Cu2}, \cite{Ma1}, \cite{Ma0} and the review \cite{Ma4}. Recently, the asymptotic stability of small solitary waves for the equation $i \, \partial_t \psi + \partial_x^2 \psi + |\psi|^{p-1} \psi = 0$ has been established for $p \simeq 3$ in \cite{Cu3}. \\
\\ Depending on the function $g$, \eqref{NLS} may (or may not) involve \textit{internal modes}, that is to say, non-trivial solutions $(\mathcal{V}_1 \, , \mathcal{V}_2 \, , \lambda ) \in H^2 ( \R )^2 \times [0 \, , + \infty )$ to the system
\begin{equation}
    \left \{ \begin{array}{ccl} \mathcal{L}_+ \mathcal{V}_1 &=& \lambda \mathcal{V}_2 \\ \mathcal{L}_- \mathcal{V}_2 &=& \lambda \mathcal{V}_2 \end{array} \right.
    \label{sysL2}
\end{equation}
where $\mathcal{L}_+ = - \partial_x^2 + \omega - 3 \phi_\omega^2 - g( \phi_\omega^2 ) - 2 \phi_\omega^2 g'( \phi_\omega^2 )$ and $\mathcal{L}_- = - \partial_x^2 + \omega - \phi_\omega^2 - g( \phi_\omega^2 )$ are the operators that appear when we linearize \eqref{NLS} around $e^{i \omega t} \phi_\omega$ (decomposing the solution into its real and imaginary parts). The existence of internal modes generates time-periodic solutions to the linearized equation around the solitary wave, which constitute potential obstacles to the asymptotic stability of solitons. Indeed, if \eqref{sysL2} has a solution $(\mathcal{V}_1 \, , \mathcal{V}_2 \, , \lambda )$ then 
\[ u_1 = \sin ( \lambda t ) \mathcal{V}_1 \, \, \, \, \, \, \, \text{and} \, \, \, \, \, \, \, u_2 = \cos ( \lambda t ) \mathcal{V}_2 \]
solve the system $\left \{ \begin{array}{l} \partial_t u_1 = \mathcal{L}_- u_2 \\ \partial_t u_2 = - \mathcal{L}_+ u_1 \end{array} \right.$ which is precisely the linearized system around the soliton $\phi_\omega$. As examples, $g(s) = -s^2$ is a case without internal mode (see \cite{Pe} and \cite{Ma1}) while $g(s) = s^2$ is a case with an internal mode (see \cite{Pe} and \cite{Ma0}). In the case $g=0$ (integrable case), there is a resonance (see \cite{Ch}), which justifies why we ask for $g \not\equiv 0$ in hypothesis $(H_1)$. Thus the sign of the perturbation determines whether there exists an internal mode or not; see \cite{Ch}, \cite{CG} and \cite{Pe} for related discussions. \\
\\ A general analysis of the case without internal mode has been conducted in \cite{Ri} (see Theorems 1 and 2): under a certain hypothesis on the function $g$, it is shown that there is no internal mode and that asymptotic stability holds. This hypothesis encompasses in particular the case $g(s) = -s^{\sigma}$ with $\sigma > 1$. In this paper we prove that, for $g(s) = s^{\sigma}$ with $\sigma > 1$, just like $g(s) = s^2$, there exists a unique internal mode and that, despite this internal mode, asymptotic stability holds. More generally, we introduce the following hypothesis, that we will comment later on:
\[ \begin{array}{rcl} (H_2) & : & \displaystyle{\lim_{\omega \to 0} \frac{1}{\varepsilon_\omega^2 \sqrt{\omega}} \int_{\R} \mathfrak{B} ( \phi_\omega^2 ) \, \text{d}x = + \infty,}
\\ & & \displaystyle{\text{where $\mathfrak{B} (s) := -3g(s) + sg'(s) + 4 \frac{G(s)}{s}$ and $\varepsilon_\omega := \max\limits_{0 \leqslant k \leqslant 4} \sup\limits_{0 \leqslant s \leqslant 3 \omega} | s^{k-1} g^{(k)} (s)|$.}} \end{array} \]
In the definition of $\varepsilon_\omega$, $3 \omega$ can be replaced by $2^+ \omega$ where $2^+$ is any constant strictly greater than $2$. Note that hypothesis $(H_1)$ implies that $\varepsilon_\omega \to 0$ as $\omega \to 0$. We shall prove that hypotheses $(H_1)$ and $(H_2)$ are enough to ensure the existence of a unique internal mode.

\begin{leftbar}
\noindent \textbf{Theorem 1.} Assume that hypotheses $(H_1)$ and $(H_2)$ hold. Then, for $\omega >0$ small enough, the system \eqref{sysL2} has a solution $(\mathcal{V}_1 \, , \mathcal{V}_2 \, , \omega \lambda ) \in H^2 ( \R )^2 \times [0 \, , + \infty )$ where $( \mathcal{V}_1 \, , \mathcal{V}_2 ) \neq (0 \, , 0)$ and $\lambda \to 1^-$ as $\omega \to 0$. Moreover, the only solutions $( \tilde{\mathcal{V}}_1 \, , \tilde{\mathcal{V}}_2 \, , \omega \tilde{\lambda}  ) \in H^2 ( \R )^2 \times [0 \, , + \infty )$ of the system \eqref{sysL2} are:
\begin{itemize}
    \item $(0 \, , 0 \, , \mu)$ for any $\mu \geqslant 0$,
    \item $(a \phi_\omega ' \, , b \phi_\omega \, , 0 )$ for any $a,b \in \R$,
    \item $(c \mathcal{V}_1 \, , c \mathcal{V}_2 \, , \omega \lambda )$ for any $c \in \R$.
\end{itemize}
\end{leftbar}

\noindent \textit{Remark 1.} Properties and estimates of this internal mode $( \mathcal{V}_1 \, , \mathcal{V}_2 )$ can be found in Proposition 2 in section 2 (for their rescaled counterparts $(V_1 \, , V_2)$, which will be introduced at the beginning of section 2). \\
\\ \textit{Remark 2.} We easily check that hypotheses $(H_1)$ and $(H_2)$ hold in the case $g(s) = s^\sigma$ with $\sigma > 1$. Indeed, we have $\mathfrak{B} (s) = \frac{( \sigma -1)^2}{\sigma +1} s^\sigma$, $\varepsilon_\omega = C_\sigma \omega^{\sigma -1}$ and $\phi_\omega (x) \geqslant c \sqrt{\omega} e^{- \sqrt{\omega} |x|}$, with $C_\sigma >0$ and $c>0$ constants that do not depend on $\omega$. Therefore
\[ \frac{1}{\varepsilon_\omega^2 \sqrt{\omega}} \int_{\R} \mathfrak{B} ( \phi_\omega^2 ) \, \text{d}x = \frac{C_\sigma ( \sigma -1)^2}{\sigma +1} \omega^{\frac{3}{2} - 2 \sigma} \int_{\R} \phi_\omega^{2 \sigma} \geqslant \tilde{C}_\sigma \omega^{1-\sigma} \, \underset{\omega \to 0^+}{\longrightarrow} \, + \infty, \]
which proves that $(H_2)$ holds in this case. Hypothesis $(H_2)$ still holds in the case $g(s) = a_1 s^{\sigma_1} + \cdots + a_N s^{\sigma_N}$ with $1 < \sigma_1 < \cdots < \sigma_N$, $a_1 > 0$ and $a_i \in \R$ for $i \geqslant 2$. \\
\\ \textit{Remark 3.} The hypothesis $(H_2)$ echoes to the hypothesis $(H_2)$ in \cite{Ri}. We sum up both cases with the notation of the present paper as follows:
\begin{itemize}
    \item if $\frac{1}{\varepsilon_\omega^2 \sqrt{\omega}} \int_{\R} \mathfrak{B} ( \phi_\omega^2 ) \, \text{d}x \, \underset{\omega \to 0}{\longrightarrow} \, - \infty$, then there is no internal mode and the asymptotic stability result holds for small $\omega$;
    \item if $\frac{1}{\varepsilon_\omega^2 \sqrt{\omega}} \int_{\R} \mathfrak{B} ( \phi_\omega^2 ) \, \text{d}x \, \underset{\omega \to 0}{\longrightarrow} \, + \infty$, then we are in the situation of the present paper and there exists a unique internal mode (see Theorem 1);
    \item if none of these limits hold, which includes in particular the integrable case, the situation is expected to be more complicated.
\end{itemize}
The fact that the same integral appears in both cases is natural. Indeed, the construction of the internal mode (or the proof of the absence of internal modes in \cite{Ri}) relies on a factorization introduced in \cite{Ma1}: $\mathcal{S}^2 \mathcal{L}_+ \mathcal{L}_- = \mathcal{M}_+ \mathcal{M}_- \mathcal{S}^2$, where $\mathcal{S} = \phi_\omega \cdot \partial_x \cdot \phi_\omega^{-1}$ and $\mathcal{M}_{\pm} = - \partial_x^2 + \omega + \mathfrak{a}_\omega^{\pm}$, with $\mathfrak{a}_\omega^+ = g ( \phi_\omega^2 ) - 2 \frac{G ( \phi_\omega^2 )}{\phi_\omega^2}$ and $\mathfrak{a}_\omega^- = 5 g ( \phi_\omega^2 ) -6 \frac{G ( \phi_\omega^2 )}{\phi_\omega^2} - 2 \phi_\omega^2 g'( \phi_\omega^2 )$. The fact that such a factorization still exists in our perturbative case is a remarkable key ingredient. The analysis of the internal mode (or the absence of internal modes) involves the integral $\int_{\R} ( \mathfrak{a}_\omega^+ + \mathfrak{a}_\omega^- )$, which is precisely the integral involved in the hypothesis $(H_2)$ since $\mathfrak{B} ( \phi_\omega^2 ) = - \frac{\mathfrak{a}_\omega^+ + \mathfrak{a}_\omega^-}{2}$. The arguments linking the existence of an internal mode to the sign of this integral come from \cite{Si} and \cite{Me}. Roughly, if this integral is positive in some sense, we have an hypothesis of repulsivity: there is no internal mode, we can directly use virial arguments on the transformed problem that involves $(\mathcal{M}_- \, , \mathcal{M}_+ )$ and establish the asymptotic stability that way. On the other hand, if this integral is negative in some sense (case of the present paper), we do not have repulsivity on the potentials of $(\mathcal{M}_- \, , \mathcal{M}_+)$: there is an internal mode and we use a second factorization in order to end up with a repulsive potential and use virial arguments to prove the asymptotic stability. This second factorization will be displayed in section 2 below (Lemma 2). \\
\\ The internal mode will be constructed and studied in the section 2 below, and in particular in Proposition 2. Its understanding is the first one of the two main ingredients of the proof for the asymptotic stability of the small solitons. The second main ingredient is the \textit{Fermi golden rule}, which aims at proving that the internal mode component of the solution is nonlinearly damped. The approach here is inspired by \cite{Ko1}, \cite{Ko2} and \cite{Ma0}. The idea is that, in the proof, it is crucial that a certain constant does not vanish. \\
\\ To introduce the Fermi golden rule hypothesis, we need some quantities which may appear cryptic for the moment, but they will be explained and related to the proof of the asymptotic stability in sections 5 and 6. In section 5 we will show that there exist $\mathfrak{g}_1$ and $\mathfrak{g}_2$ non-trivial bounded even solutions of the system
\[ \left \{ \begin{array}{ccl} \mathcal{L}_+ \mathfrak{g}_1 &=& 2 \omega \lambda \mathfrak{g}_2 \\ \mathcal{L}_- \mathfrak{g}_2 &=& 2 \omega \lambda \mathfrak{g}_1, \end{array} \right. \]
where $\lambda$ is the eigenvalue introduced in Theorem 1. Now we introduce
\[ \mathcal{G}_1 = \mathcal{V}_1^2 \phi_\omega ( 3 + 3 g'( \phi_\omega^2 ) + 2 \phi_\omega^2 g''( \phi_\omega^2 )) - \mathcal{V}_2^2 \phi_\omega ( 1 + g'( \phi_\omega^2)) \, \, \, \, \, \, \, \text{and} \, \, \, \, \, \, \, \mathcal{G}_2 = 2 \mathcal{V}_1 \mathcal{V}_2 \phi_\omega (1 + g'( \phi_\omega^2 )), \]
where $\langle \cdot , \cdot \rangle$ denotes the scalar product in $L^2 ( \R )$. The Fermi golden rule hypothesis we will need is the following:
\[ \begin{array}{ccc} (H_3) & : & \text{there exists a positive quantity $\underline{\textbf{FGR}} ( \omega_0 )$ depending only on $\omega_0$ such that,} \\
\\ & & | \omega - \omega_0 | \leqslant \frac{\omega_0}{2} \, \, \, \, \, \Longrightarrow \, \, \, \, \, \int_{\R} \left ( \mathcal{G}_1 \mathfrak{g}_1 + \mathcal{G}_2 \mathfrak{g}_2 \right ) \geqslant \underline{\textbf{FGR}} ( \omega_0) > 0. \end{array} \]
Note that $\underline{\textbf{FGR}} (0)=0$ (integrable case). Combining the hypotheses $(H_2)$ - the control of the internal mode - and $(H_3)$ - the Fermi golden rule -, we are able to prove an asymptotic stability result. 

\begin{leftbar}
\noindent \textbf{Theorem 2.} Assume that hypotheses $(H_1)$ and $(H_2)$ hold. Assume that the Fermi golden rule hypothesis $(H_3)$ also holds. Then, for $\omega_0 > 0$ small enough, there exists $\delta > 0$ with the following property: for any even function $\psi_0 \in H^1 ( \R )$ with $|| \psi_0 - \phi_{\omega_0} ||_{H^1 ( \R )} < \delta$, there exist $\omega_+ > 0$ and a $\mathscr{C}^1$ function $\gamma \, : \, [ 0 \, , + \infty ) \to \R$ with $\lim\limits_{+ \infty} \gamma ' = \omega_+$ such that, if $\psi$ denotes the solution of \eqref{NLS} with initial data $\psi(0) = \psi_0$, then, for any bounded interval $I \subset \R$,
\[ \lim_{t \to + \infty} || \psi (t) - e^{i \gamma (t)} \phi_{\omega_+} ||_{L^\infty (I)} = 0. \]
\end{leftbar}

\noindent \textit{Remark 4.} As it is pointed out in \cite{Ma0}, the symmetry assumption in Theorem 2 is technical, in the sense that it simplifies the proof, but no deep additional difficulty is expected in the non symmetric case. \\
\\ \textit{Remark 5.} It can be shown that, for any $\eta >0$, $\delta >0$ may be chosen small enough so that $\omega_+ \in [\omega_0 - \eta \, , \omega_0 + \eta ]$. \\
\\ As it will be shown in section 5, although hypothesis $(H_3)$ may appear difficult to check in general, it is rigorously justified for $\sigma \sim 1^+$ and numerically verified for $g(s) = s^\sigma$, where $\sigma > 1$. Henceforth, the asymptotic stability result holds for such cases, since all three hypotheses $(H_1)$, $(H_2)$ and $(H_3)$ are verified. \\
\\ The layout of this paper is globally adapted from \cite{Ma0} and \cite{Ri}. As said previously, the main arguments are the understanding of the internal mode and the Fermi golden rule. Once these two points studied, the rest of the proof is almost unchanged. In section 2, we construct the internal mode and give some of its properties. Using the identity $\mathcal{S}^2 \mathcal{L}_+ \mathcal{L}_- = \mathcal{M}_+ \mathcal{M}_- \mathcal{S}^2$ mentioned above, we first construct the internal mode for $(\mathcal{M}_+ \, , \mathcal{M}_-)$ then come back to $(\mathcal{L}_+ \, , \mathcal{L}_-)$. We then introduce the second factorization that will lead further to the second transformed problem, based on a new differential operator $K$. Finally we prove a sort of coercivity property on the operator $K$, and the uniqueness of the internal mode. In section 3, we introduce the rescaled modulation decomposition of the solution, a standard decomposition for stability arguments, and in particular we introduce the internal mode component of the solutions (that will be denoted $b$). In section 4, we prove a first virial argument directly on the solution, without transformation of the linearized operators. In section 5, we study the second main point of the proof: the Fermi golden rule. We will explain how hypothesis $(H_3)$ can be explicitly checked for $g(s) = s^\sigma$ with $\sigma >1$ using simple numerical computations. In section 6, we control the internal mode component of the solution: more precisely, we control $\int_0^s |b|^4$. This is the estimate that requires the Fermi golden rule. In section 7, we introduce the setting of the double transformed problem and technical results related to it; in section 8, we prove coercivity results that will enable us to go back from the transformed problem to the original problem. In section 9, we prove the second virial argument, on the transformed problem this time. Gathering all previous results, in section 10 we finally prove the Theorem 1, assuming all three hypotheses $(H_1)$, $(H_2)$ and $(H_3)$ hold. \\
\\ In all this paper, assume hypothesis $(H_1)$ holds. \\
\\ The letters $u$, $v$, $w$ and $z$ will denote complex-valued functions; we will index by $1$ their real part and by $2$ their imaginary part. The $L^2$ scalar product will be denoted by $\langle u \, , v \rangle = \text{Re} \left ( \int_{\R} u \overline{v} \, \text{d}x \right )$ and the $L^2$ norm by $|| \cdot ||$. The $H^1$ norm will be denoted by $|| \cdot ||_{H^1 ( \R )}$. The scalar product in $\R^N$ will be denoted by $\cdot$. Lastly, the letter $C$ will denote various positive constants whose expression change from one line to another. The concerned constants do not depend on the parameters $\omega_0$, $\epsilon$, $\theta$, $\vartheta$, $A$ and $B$ (that will be introduced in sections 3, 4 and 7, except in the proof of Proposition 5 and in section 10, when some of these parameters are already fixed. \\
\\ This paper is the result of many discussions with Yvan Martel. The motivation of this paper and its proof are inspired by his paper \cite{Ma0}. May he be warmly thanked for it here.

\section{Construction of the internal mode}

\subsection{Properties of the solitons}
\noindent We begin by recalling some properties of the solitons, and proving an additional estimate for future use. First, and until the end of this paper, let us rescale the solitons: $\phi_\omega (x) = \sqrt{\omega} \, Q_\omega ( \sqrt{\omega} \, x)$. We denote $y = x / \sqrt{\omega}$ the rescaled variable. Now $Q_\omega$ is solution of the equation
\begin{equation}
    Q_\omega '' = Q_\omega - Q_\omega^3 - \frac{g(\omega Q_\omega^2)}{\omega} Q_\omega
    \label{eqQ}
\end{equation}
Integrating this equation, we find the following useful relation:
\begin{equation}
    (Q_\omega ')^2 = Q_\omega^2 - \frac{1}{2} Q_\omega^4 - \frac{G(\omega Q_\omega^2)}{\omega^2}.
    \label{eqQ'2}
\end{equation}
From \cite{Ri} we recall the following estimates: for $\omega >0$ small enough, for any $k \geqslant 0$ there exist positive constants $c_k$ and $C_k$ such that $c_k e^{-|y|} \leqslant |Q_\omega^{(k)} (y)| \leqslant C_k e^{-|y|}$ for all $y \in \R$. We also recall that, for $\omega >0$ small enough,
\[ \left | \frac{g( \omega Q_\omega^2 )}{\omega} \right | + \left | \frac{G ( \omega Q_\omega^2 )}{\omega^2 Q_\omega^2} \right | + \left | Q_\omega^2 g'(\omega Q_\omega^2) \right | + | \omega Q_\omega^4 g''( \omega Q_\omega^2 ) | \leqslant C \varepsilon_\omega Q_\omega^2 \leqslant C \varepsilon_\omega e^{-2|y|}. \]
Recall that $\varepsilon_\omega$ is defined in hypothesis $(H_2)$. These quantities are involved in the linearized operators we will have to deal with. In our linearization, we decompose the solution $\psi$ into its real and imaginary parts. Linearizing that way the equation \eqref{NLS} and rescaling, we obtain the operators
\[ \begin{array}{rl} & \displaystyle{L_+ = - \partial_y^2 + 1 -3Q_\omega^2 - \frac{g( \omega Q_\omega^2)}{\omega} - 2 Q_\omega^2 g'(\omega Q_\omega^2)} \\ \\ \text{and} & \displaystyle{L_- = - \partial_y^2 + 1 - Q_\omega^2 - \frac{g(\omega Q_\omega^2)}{\omega}.} \end{array} \]
Spectral properties of the operators $L_+$ and $L_-$ can be found in \cite{We1}. Let $S = Q_\omega \cdot \partial_y \cdot \frac{1}{Q_\omega}$ and $S^* = - \frac{1}{Q_\omega} \cdot \partial_y \cdot Q_\omega$. We recall from \cite{Ri} (Lemma 6) the key relation $S^2 L_+ L_- = M_+ M_- S^2$ with $M_{\pm} = - \partial_y^2 + 1 + a_\omega^{\pm}$, where
\[ \begin{array}{rl} & \displaystyle{a_\omega^+ = \frac{g(\omega Q_\omega^2)}{\omega} -2 \frac{G(\omega Q_\omega^2)}{\omega^2 Q_\omega^2}} \\ \\ \text{and} & \displaystyle{a_\omega^- = 5 \frac{g(\omega Q_\omega^2)}{\omega^2 Q_\omega^2} - 6 \frac{G( \omega Q_\omega^2)}{\omega^2 Q_\omega^2} - 2 Q_\omega^2 g'(\omega Q_\omega^2).} \end{array} \]
The previous bound shows the following crucial estimate: for $\omega >0$ small enough and all $y \in \R$,
\begin{equation}
    | a_\omega^{\pm} (y) | \leqslant C \varepsilon_\omega Q_\omega^2 (y) \leqslant C \varepsilon_\omega e^{-2|y|}.
    \label{estia}
\end{equation}
We simply write $Q :=Q_0$, solution of $Q'' = Q - Q^3$. We can write $Q_\omega$ as an expansion of $Q$ in the following sense.

\begin{leftbar}
\noindent \textbf{Lemma 1.} Assume that hypothesis $(H_1)$ holds. Let $D_\omega := Q_\omega - Q$. For $\omega > 0$ small enough and any $k \in \{ 0 \, , ... \, , 6 \}$,
\[ \forall y \in \R , \, \, \, |D_\omega^{(k)} (y)| \leqslant C \varepsilon_\omega e^{-|y|}. \]
\end{leftbar}

\noindent \textit{Proof.} We compute $D_\omega '' = D_\omega - D_\omega ( Q_\omega^2 + Q Q_\omega + Q^2) - \frac{g( \omega Q_\omega^2 )}{\omega} Q_\omega$. Defining $L_+^0 = - \partial_y^2 + 1 - 3 Q^2$, the equation satisfied by $D_\omega$ can be rewritten as
\[ \begin{array}{rcl} L_+^0 D_\omega &=& \displaystyle{D_\omega ( Q_\omega^2 + QQ_\omega - 2 Q^2 ) + \frac{g(\omega Q_\omega^2)}{\omega} Q_\omega} \\ \\ &=& \displaystyle{D_\omega^2 ( Q_\omega + 2Q) + \frac{g( \omega Q_\omega^2 )}{\omega} Q_\omega.} \end{array} \]
We know (see \cite{Ch}) that $L_+^0$ has only one negative eigenvalue: $-3$, associated with the eigenfunction $Q^2$. It is also known that the kernel of $L_+^0$ is generated by $Q'$. We recall the following spectral inequality from \cite{Tao}: there exists positive constants $c_1,c_2,c_3$ such that
\[ \langle L_+^0 D_\omega \, , D_\omega \rangle \geqslant c_1 || D_\omega ||_{H^1}^2 - c_2 | \langle D_\omega \, , Q^2 \rangle |^2 - c_3 | \langle D_\omega \, , Q' \rangle |^2. \]
Here, $D_\omega$ is even and $Q'$ is odd, thus $\langle D_\omega \, , Q' \rangle = 0$. In order to estimate the other terms, we recall the following result from Lemma 2 in \cite{Ri}: for any $\delta > 0$, for $\omega > 0$ small enough we have $|D_\omega (y)| \leqslant \delta e^{-|y|}$. Let $\delta \in (0 \, , 1)$, to be fixed later. This implies that $||D_\omega||_\infty \leqslant \delta$ and $||D_\omega||^2 \leqslant \delta^2 \leqslant \delta$. First,
\[ \left | \langle (Q_\omega + 2 Q) D_\omega^2 \, , D_\omega \rangle \right | =  \left | \int_{\R} (Q_\omega + 2Q) D_\omega^3 \right | \leqslant C \delta ||D_\omega||^2. \]
Second, 
\[ \left | \left \langle \frac{g( \omega Q_\omega^2 )}{\omega} Q_\omega \, , D_\omega \right \rangle \right | \leqslant C \varepsilon_\omega \int_{\R} Q_\omega^2 D_\omega \leqslant C \varepsilon_\omega || Q_\omega^2 || \, || D_\omega || \leqslant C \varepsilon_\omega || D_\omega ||. \]
Hence, $| \langle L_+^0 D_\omega \, , D_\omega \rangle | \leqslant C ( \varepsilon_\omega ||D_\omega || + \delta ||D_\omega||^2 )$. Now, let us estimate the projection $\langle D_\omega \, , Q^2 \rangle$. Using the facts that $L_+^0 Q^2 = -3 Q^2$ and that $L_+^0$ is self-adjoint, we write that
\[ \begin{array}{rcl} | \langle D_\omega \, , Q^2 \rangle | &=& \displaystyle{\frac{1}{3} | \langle D_\omega \, , -3 Q^2 \rangle | \, \, = \, \, \frac{1}{3} | \langle D_\omega \, , L_+^0 Q^2 \rangle | \, \, = \, \, C | \langle L_+^0 D_\omega \, , Q^2 \rangle |} \\
\\ & \leqslant & \displaystyle{C \left ( \left | \left \langle (Q_\omega + 2Q) D_\omega^2 \, , Q^2 \right \rangle \right | + \left | \left \langle \frac{g(\omega Q_\omega^2)}{\omega} Q_\omega \, , Q^2 \right \rangle \right | \right )} \end{array} \]
where
\[ \begin{array}{rl} & \displaystyle{\left | \left \langle (Q_\omega + 2Q) D_\omega^2 \, , Q^2 \right \rangle \right | = \int_{\R} (Q_\omega + 2Q)Q^2 D_\omega^2 \leqslant C ||D_\omega||^2} \\
\\ \text{and} & \displaystyle{\left | \left \langle \frac{g(\omega Q_\omega^2)}{\omega} Q_\omega \, , Q^2 \right \rangle \right | = \left | \int_{\R} \frac{g(\omega Q_\omega^2)}{\omega} Q_\omega Q^2 \right | \leqslant C \varepsilon_\omega.} \end{array} \]
Thus, $| \langle D_\omega \, , Q^2 \rangle |^2 \leqslant C ( ||D_\omega||^2 + \varepsilon_\omega )^2 \leqslant C ( \varepsilon_\omega^2 + \varepsilon_\omega ||D_\omega||^2 + \delta ||D_\omega||^2 )$. Using the spectral inequality, we find that
\[ ||D_\omega||^2 \leqslant ||D_\omega||_{H^1}^2 \leqslant C \left ( |\langle L_+^0 D_\omega \, , D_\omega \rangle | + | \langle D_\omega \, , Q^2 \rangle | \right ) \leqslant C ( \varepsilon_\omega + \delta ) ||D_\omega||^2 + C \varepsilon_\omega ||D_\omega|| + C \varepsilon_\omega^2. \]
Recalling that all the letters $C$ refer to constant that do not depend on $\omega$, we fix $\delta \in (0 \, , 1)$ such that $C \delta < \frac{1}{4}$. Then we take $\omega > 0$ small enough such that $C \varepsilon_\omega < \frac{1}{4}$ and $|D_\omega (y)| \leqslant \delta e^{-|y|}$. We get
\[ \begin{array}{rl} & \displaystyle{\frac{1}{2} ||D_\omega||^2 - C \varepsilon_\omega ||D_\omega|| - C \varepsilon_\omega^2 \leqslant \left ( 1 - C \varepsilon_\omega - C \delta \right ) ||D_\omega||^2 - C \varepsilon_\omega ||D_\omega|| - C \varepsilon_\omega^2 \leqslant 0} \\
\\ \text{thus} & ||D_\omega||^2 - C \varepsilon_\omega ||D_\omega|| - C \varepsilon_\omega^2 \leqslant 0. \end{array} \]
The only positive root of the polynomial $X^2 - C \varepsilon_\omega X - C \varepsilon_\omega^2$ is $C \varepsilon_\omega$ (where $C$ is a different positive constant), thus we obtain $||D_\omega|| \leqslant C \varepsilon_\omega$. This leads to
\[ ||D_\omega||_{H^1}^2 \leqslant C ( \varepsilon_\omega + \delta ) ||D_\omega||^2 + C \varepsilon_\omega ||D_\omega|| + C \varepsilon_\omega^2 \leqslant C \varepsilon_\omega^2 \]
and then, using Sobolev's inequality, $||D_\omega||_\infty \leqslant C ||D_\omega||_{H^1} \leqslant C \varepsilon_\omega$. \\
\\ By the variation of the constants, we find the following expressions of $D_\omega$ and $D_\omega '$: for $y>0$,
\[ \begin{array}{rl} & D_\omega (y) = \frac{e^{-y}}{2} \left ( D_\omega (0) - \int_0^{+ \infty} Z_\omega(z) e^z \, \text{d}z \right ) - \frac{e^y}{2} \int_y^{+ \infty} Z_\omega (z) e^{-z} \, \text{d}z + \frac{e^{-y}}{2} \int_y^{+ \infty} Z_\omega (z) e^z \, \text{d}z \\
\\ \text{and} & D_\omega '(y) = \frac{e^{-y}}{2} \left ( \int_0^{+ \infty} Z_\omega(z) e^z \, \text{d}z  -D_\omega (0) \right ) - \frac{e^y}{2} \int_y^{+ \infty} Z_\omega (z) e^{-z} \, \text{d}z - \frac{e^{-y}}{2} \int_y^{+ \infty} Z_\omega (z) e^z \, \text{d}z, \end{array} \]
where $Z_\omega := -D_\omega ( Q_\omega^2 + Q Q_\omega + Q^2 ) - \frac{g( \omega Q_\omega^2 )}{\omega} Q_\omega$. Using the estimate $| D_\omega | \leqslant C \varepsilon_\omega$, we find that $|Z_\omega (z)| \leqslant C \varepsilon_\omega e^{-2z}$ for $z>0$, then it leads to $|D_\omega (y)| + |D_\omega (y)| \leqslant C \varepsilon_\omega e^{-y}$. Using the equation $D_\omega '' = D_\omega + Z_\omega$ and the fact that $D_\omega$ is even, we conclude that $| D_\omega^{(k)} (y) | \leqslant C \varepsilon_\omega e^{-|y|}$ for any $y \in \R$ and any $k \in \{ 0 \, , ... \, , 6 \}$. \hfill \qedsymbol

\subsection{The internal mode}
\noindent The spectral problem we study here, whose solutions are called \textit{internal modes}, is
\begin{equation}
    \left \{ \begin{array}{ccl} L_+ V_1 & = & \lambda V_2 \\ L_- V_2 &=& \lambda V_1 \end{array} \right.
    \label{sysL}
\end{equation}
with an eigenvalue $\lambda$ close to $1$. Let us define the following related system:
\begin{equation}
    \left \{ \begin{array}{ccl} M_+ W_1 & = & \lambda W_2 \\ M_- W_2 &=& \lambda W_1. \end{array} \right.
    \label{sysM}
\end{equation}
Observe that, if \eqref{sysM} is satisfied by $(W_1 \, , W_2)$, then \eqref{sysL} is satisfied by $V_1 = (S^*)^2 W_1$ and $V_2 = \lambda^{-1} L_+ V_1$. \\
\\ We introduce the following notation that will be used throughout the entire article:
\[ I_\omega = - \int_{\R} ( a_\omega^+ + a_\omega^- ) (y) \, \text{d}y \, \, \, \, \, \, \, \, \, \, \text{and} \, \, \, \, \, \, \, \, \, \, \varrho_\omega = \frac{\varepsilon_\omega^2}{I_\omega}. \]
Hypothesis $(H_2)$ exactly means that: 
\[ \frac{1}{\varrho_\omega} = \frac{I_\omega}{\varepsilon_\omega^2} \, \underset{\omega \to 0}{\longrightarrow} \, + \infty. \]
Thus, hypothesis $(H_2)$ implies in particular that $I_\omega > 0$ for $\omega > 0$ small enough. It also implies that $\varrho_\omega \to 0$ as $\omega \to 0$. These assumptions are crucial in what follows. We also notice from \eqref{estia} that $I_\omega \leqslant C \varepsilon_\omega$ and thus $\varepsilon_\omega \leqslant C \varrho_\omega$. \\
\\ We prove the existence of an internal mode.

\begin{leftbar}
\noindent \textbf{Proposition 2.} Assume that hypotheses $(H_1)$ and $(H_2)$ hold. There exists $\omega_1 > 0$, a $\mathscr{C}^1$ function $\alpha \, : \, (0 \, , \omega_1) \to [0 \, , + \infty )$ and smooth real-valued $\omega$-dependent functions $W_1 (y)$ and $W_2 (y)$ such that the following properties hold for all $\omega \in (0 \, , \omega_1 )$. 
\begin{itemize}
    \item \textit{(Expansion of the eigenvalue.)} As $\omega \to 0$, $\alpha ( \omega ) = \frac{I_\omega}{4} \left ( 1 + \mathcal{O} ( \varrho_\omega ) \right )$.
    \item \textit{(Solutions to the spectral problem.)} The pair of functions $(W_1 \, , W_2)$ solves \eqref{sysM} with $\lambda = 1 - \alpha^2$ and the pair of functions $(V_1 \, , V_2) = ((S^*)^2 W_1 \, , \lambda^{-1} L_+ V_1)$ solves \eqref{sysL}. 
    \item \textit{(Expansion of the eigenfunctions.)} For $j \in \{ 1 \, , 2 \}$, $W_j = 1 + S_j + \widetilde{W_j}$ and
    \[ V_1 = 1 - Q^2 + R_1 + \widetilde{V_1} , \, \, \, \, \, V_2 = 1 + R_2 + \widetilde{V_2}, \]
    where, for all $y \in \R$,
    \[ \begin{array}{llll} |R_j^{(k)} (y)| \leqslant C \varepsilon_\omega (1 + |y|) & \text{for all $k \in \{ 0 \, , ... \, , 4 \}$,} & |S_j^{(k)} (y)| \leqslant C \varepsilon_\omega (1 + |y|) & \text{for all $k \in \N$,} \\ |\widetilde{W}_j^{(k)} (y)| \leqslant C \varepsilon_\omega^2 ( 1 + y^2 ) & \text{for all $k \in \N$,} & | \widetilde{V}_j^{(k)} (y)| \leqslant C \varepsilon_\omega^2 ( 1+y^2) & \text{for all $k \in \{ 0 \, , ... \, , 2 \}$.} \end{array} \]
    \item \textit{(Decay properties of the eigenfunctions.)} For $j \in \{ 1 \, , 2 \}$ and all $y \in \R$,
    \[ \begin{array}{l} |W_j(y)| \leqslant C e^{- \alpha |y|} \, \, \, \, \, \text{and, for any $k \geqslant 1$,} \, \, |W_j^{(k)} (y)| \leqslant C \left ( \varepsilon_\omega I_\omega^{k-1} e^{- \alpha |y|} + \varepsilon_\omega e^{- \kappa |y|} \right ), \\ |V_j(y)| \leqslant C e^{- \alpha |y|} \, \, \, \, \, \text{and, for any $k \in \{ 1 \, , 2 \}$,} \, \, |V_j^{(k)} (y)| \leqslant C \left ( \varepsilon_\omega I_\omega^{k-1} e^{- \alpha |y|} + e^{- |y|} \right ), \end{array} \]
    where $\kappa = \sqrt{2 - \alpha^2}$. Moreover, for any $k \geqslant 0$ and all $y \in \R$,
    \[ | ( W_1 - W_2 )^{(k)} (y) | \leqslant C \varepsilon_\omega e^{- \kappa |y|}. \]
    \item \textit{(Asymptotics of the eigenfunctions.)} For $j \in \{ 1 \, , 2 \}$ and all $y \in \R$,
    \[ |W_j (y) - e^{- \alpha |y|} | +  |V_1(y) - (1-Q^2(y))e^{-\alpha |y|}| +  |V_2 (y) - e^{-\alpha |y|}| \leqslant C \varrho_\omega e^{- \alpha |y|}. \]
    Moreover, for all $y > 0$,
    \[ |W_j ' (y) + \alpha e^{- \alpha y} | \leqslant C \varrho_\omega I_\omega e^{- \alpha y} + C \varepsilon_\omega e^{- \kappa y}. \]
    Finally,
    \[ \left | \langle W_1 \, , W_2 \rangle - \frac{1}{\alpha} \right | +  \left | \langle V_1 \, , V_2 \rangle - \frac{1}{\alpha} \right | \leqslant C \, \frac{\varrho_\omega}{I_\omega}. \]
    \item \textit{(Derivatives with regards to $\omega$.)} First,
    \[ |\omega \, \alpha '( \omega) | \leqslant C \varepsilon_\omega \, \, \, \, \, \, \, \, \, \text{and} \, \, \, \, \, \, \, \, \, \omega \left | \alpha ' ( \omega ) - \frac{1}{4} \partial_\omega I_\omega \right | \leqslant C \varepsilon_\omega^2. \]
    Moreover, for $j \in \{ 1 \, , 2 \}$, any $k \geqslant 1$ and all $y \in \R$,
    \[ \begin{array}{rl} & | \partial_\omega W_j | \leqslant \frac{C \varepsilon_\omega \varrho_\omega}{\omega I_\omega} (1 + |y|) e^{- \alpha |y|} + \frac{C \varepsilon_\omega}{\omega} (1 + |y|) e^{- \kappa |y|} \leqslant \frac{C \varepsilon_\omega \varrho_\omega}{\omega I_\omega} (1 + |y|) e^{- \alpha |y|} \\
    \\ \text{and} & | \partial_y^k \partial_\omega W_j | \leqslant \frac{C \varepsilon_\omega^k}{\omega} (1 + |y|) e^{- \alpha |y|} + \frac{C \varepsilon_\omega}{\omega} (1 + |y|) e^{- \kappa |y|} \leqslant \frac{C \varepsilon_\omega}{\omega} (1 + |y|) e^{- \alpha |y|}. \end{array} \]
    Finally, for $j \in \{ 1 \, , 2 \}$ and all $y \in \R$,
    \[ | \partial_\omega V_j | + | \partial_y \partial_\omega V_j | \leqslant \frac{C}{\omega} \left ( \frac{\varepsilon_\omega \varrho_\omega}{I_\omega} +1 \right ) (1 + |y|) e^{- \alpha |y|}. \]
\end{itemize}

\end{leftbar}

\noindent \textit{Proof.} We define $\alpha > 0$ and $\kappa > 0$ such that $\lambda = 1 - \alpha^2$ and $\kappa^2 = 1 + \lambda = 2 - \alpha^2$. We introduce $X_1 = \frac{W_1 + W_2}{2}$ and $X_2 = \frac{W_1 - W_2}{2}$. The system \eqref{sysM} on $(W_1 \, , W_2)$ is equivalent to the following system on $(X_1 \, , X_2)$:
\begin{equation}
    \left \{ \begin{array}{ccl} - \partial_y^2 X_1 + \alpha^2 X_1 + b_\omega^+ X_1 + b_\omega^- X_2 &=& 0 \\ - \partial_y^2 X_2 + \kappa^2 X_2 + b_\omega^- X_1 + b_\omega^+ X_2 &=& 0 \end{array} \right.
    \label{sysX}
\end{equation}
where $b_\omega^{\pm} = \frac{a_\omega^+ \pm a_\omega^-}{2}$. We introduce the following matrix notation:
\[ X = \left ( \begin{array}{c} X_1 \\ X_2 \end{array} \right ) , \, \, \, \, H_\alpha = \left ( \begin{array}{cc} - \partial_y^2 + \alpha^2 & 0 \\ 0 & - \partial_y^2 + \kappa^2 \end{array} \right ) , \, \, \, \, P_\omega = \left ( \begin{array}{cc} b_\omega^+ & b_\omega^- \\ b_\omega^- & b_\omega^+ \end{array} \right ). \]
The system \eqref{sysX} is equivalent to the matrix equation
\begin{equation}
    (H_\alpha + P_\omega) X = 0.
    \label{sysH}
\end{equation}
We use Birman-Schwinger arguments similar to the ones developed in \cite{Me}. To this end, let us introduce
\[ \begin{array}{rl} & \displaystyle{ |P_\omega|^{\frac{1}{2}} = \frac{1}{2} \left ( \begin{array}{cc} \sqrt{|a_+|} + \sqrt{|a_-|} & \sqrt{|a_+|} - \sqrt{|a_-|} \\ \\ \sqrt{|a_+|} - \sqrt{|a_-|} & \sqrt{|a_+|} + \sqrt{|a_-|} \end{array} \right ) } \\
\\ \text{and} & \displaystyle{ P_\omega^{\frac{1}{2}} = \frac{1}{2} \left ( \begin{array}{cc} a_+^{\frac{1}{2}} + a_-^{\frac{1}{2}} & a_+^{\frac{1}{2}} - a_-^{\frac{1}{2}} \\ \\ a_+^{\frac{1}{2}} - a_-^{\frac{1}{2}} & a_+^{\frac{1}{2}} + a_-^{\frac{1}{2}} \end{array} \right ) } \end{array} \]
where $x^{\frac{1}{2}} := \text{sgn} (x) \sqrt{|x|}$ is a continuous function of $x$. The important relation satisfied by these two matrices is that $P_\omega^{\frac{1}{2}} | P_\omega |^{\frac{1}{2}} = | P_\omega |^{\frac{1}{2}} P_\omega^{\frac{1}{2}} = P_\omega$. Moreover, we recall from the estimates on $a_\omega^{\pm}$ that $|P_\omega| \leqslant C \varepsilon_\omega Q_\omega^2$ in the sense that all coefficients of the matrix $P_\omega$ satisfy this inequality. We define the operator $K_{\alpha , \omega} = P_\omega^{\frac{1}{2}} H_\alpha^{-1} |P_\omega|^{\frac{1}{2}}$ on $L^2 ( \R ) \times L^2 ( \R )$, with integral kernel
\[ K_{\alpha , \omega} (y \, , z) = \frac{1}{2} P_\omega^{\frac{1}{2}} (y) \left ( \begin{array}{cc} \frac{e^{- \alpha |y-z|}}{\alpha} & 0 \\ 0 & \frac{e^{- \kappa |y-z|}}{\kappa} \end{array} \right ) |P_\omega|^{\frac{1}{2}} (z). \]
We decompose $K_{\alpha , \omega} = L_{\alpha , \omega} + M_{\alpha , \omega}$ where
\[ \begin{array}{rl} & \displaystyle{L_{\alpha , \omega} (y \, , z) = \frac{1}{2 \alpha} P_\omega^{\frac{1}{2}} (y) \left ( \begin{array}{cc} 1 & 0 \\ 0 & 0 \end{array} \right ) |P_\omega|^{\frac{1}{2}} (z)} \\
\\ \text{and} & \displaystyle{M_{\alpha , \omega} (y \, , z) = P_\omega^{\frac{1}{2}} (y) N_\alpha (y \, , z) |P_\omega|^{\frac{1}{2}} (z) \, \, \, \, \text{with} \, \, \, \, N_\alpha (y \, , z) = \frac{1}{2} \left ( \begin{array}{cc} \frac{e^{-\alpha |y-z|} -1}{\alpha} & 0 \\ 0 & \frac{e^{- \kappa |y-z|}}{\kappa} \end{array} \right ).} \end{array} \]
We can extend these operators for $\alpha = 0$, defining
\[ M_{0,\omega} (y \, , z) = P_\omega^{\frac{1}{2}} (y) N_0 (y \, , z) |P_\omega|^{\frac{1}{2}} (z) \, \, \, \, \text{where} \, \, \, \, N_0 (y \, , z) = \frac{1}{2} \left ( \begin{array}{cc} -  |y-z| & 0 \\ 0 & \frac{1}{\sqrt{2}} e^{- \sqrt{2} |y-z|} \end{array} \right ). \]
That way, the map $\alpha \mapsto M_{\alpha , \omega}$ is well-defined and analytic (in the Hilbert-Schmidt norm) in a neighborhood of $\alpha = 0$. From the estimate $|a_\omega^{\pm}(y)| \leqslant C \varepsilon_\omega Q_\omega^2 (y)$ we deduce that $M_{\alpha,\omega}$ converges to $0$ as $\omega \to 0$. Hence, $(\alpha \, , \omega) \mapsto M_{\alpha , \omega}$ is continuous in $\omega$ and analytical in $\alpha$ in a neighborhood of $(0 \, , 0)$; and we cannot say more \textit{a priori} in terms of regularity in $\omega$. \\
\\ We observe that \eqref{sysH} is satisfied by $(\alpha \, , X)$ if, and only if, the function $\Psi = P_\omega^{\frac{1}{2}} X$ solves $\Psi = - \omega K_{\alpha , \omega} \Psi$ i.e. $\Psi + (1 + M_{\alpha , \omega})^{-1} L_{\alpha , \omega} \Psi = 0$. The existence and the analytic regularity of $(1 + M_{\alpha , \omega})^{-1}$ make sense since $||M_{\alpha , \omega}|| < 1$ for $\omega > 0$ small enough. Indeed, writing that $\left | \frac{e^{- \alpha |y-z|} - 1}{\alpha} \right | \leqslant |y-z|$ and $\left | \frac{e^{- \kappa |y-z|}}{\kappa} \right | \leqslant \frac{1}{\kappa} \leqslant C$, we see that $|M_{\alpha , \omega} (y \, , z)| \leqslant C \varepsilon_\omega Q_\omega (y) Q_\omega (z) (1 + |y-z|)$. This leads to
\[ ||M_{\alpha , \omega}|| \leqslant C \varepsilon_\omega \left ( \int_{\R} \int_{\R} Q_\omega^2 (y) Q_\omega^2 (z) (1 + |y-z|)^2 \, \text{d}y  \, \text{d}z \right )^{1/2} \leqslant C \varepsilon_\omega \left ( \int_{\R} \int_{\R} e^{-2|y|} e^{-2|z|} (1+|y-z|)^2 \, \text{d}y \, \text{d}z \right )^{1/2} \leqslant C \varepsilon_\omega \]
which proves that $|| M_{\alpha , \omega} || \leqslant C \varepsilon_\omega < 1$ for $\omega > 0$ small enough. \\
\\ Therefore, we aim at finding $\alpha > 0$ small such that $-1$ is an eigenvalue of the operator $(1 + M_{\alpha , \omega})^{-1} L_{\alpha , \omega}$. More generally, let us consider the eigenvalue problem
\begin{equation}
    (1 + M_{\alpha , \omega} )^{-1} L_{\alpha , \omega} \Psi = \mu \Psi.
    \label{eqmu}
\end{equation}
By definition, $L_{\alpha , \omega}$ is a rank one operator and, for any $\varphi \in L^2 ( \R )$, we have
\[ (L_{\alpha , \omega} \varphi)(y) = \frac{p_\omega ( \varphi )}{2 \alpha} P_\omega^{\frac{1}{2}} (y) e_1 \, \, \, \, \, \text{where} \, \, \, p_\omega ( \varphi ) = \int_{\R} e_1 \cdot \left ( |P_\omega|^{\frac{1}{2}} \varphi \right ) \, \, \text{and} \, \, e_1 = \left ( \begin{array}{c} 1 \\ 0 \end{array} \right ). \]
Here, the dot $\cdot$ denotes the usual scalar product in $\R^2$. We find that $( \mu \, , \Psi )$ satisfies \eqref{eqmu} if and only if
\begin{equation}
    p_\omega ( \Psi ) (1 + M_{\alpha , \omega})^{-1} (P_\omega^{\frac{1}{2}} e_1 ) = 2 \alpha \mu \Psi.
    \label{eqPsi}
\end{equation}
We define the function $r ( \alpha \, , \omega ) = p_\omega \left ( (1 + M_{\alpha , \omega})^{-1} ( P_\omega^{\frac{1}{2}} e_1 ) \right )$. Hence, $(\mu \, , \Psi)$ solves \eqref{eqmu} if, and only if, $r ( \alpha \, , \omega ) = 2 \alpha \mu$. Therefore, $-1$ is an eigenvalue of the operator $(1 + M_{\alpha , \omega})^{-1} L_{\alpha , \omega}$ if, and only if, $s ( \alpha \, , \omega ) = 0$, where 
\[ s ( \alpha \, , \omega ) = \alpha + \frac{r( \alpha \, , \omega )}{2}. \]
We easily see that $\frac{\partial s}{\partial \alpha} (0 \, , 0) = 1$. By the implicit function theorem, there exists a continuous function $\omega \mapsto \alpha ( \omega )$ defined in a neighborhood of $0$ such that $s(a \, , \omega ) = 0$ if, and only if, $a = \alpha ( \omega )$. Moreover, since $s$ is $\mathscr{C}^1$ with regards to $\omega$ on $(0 \, , \omega^*)$ for a certain $\omega^* > 0$, we know from the implicit function theorem that $\alpha$ is $\mathscr{C}^1$ in a neighborhood of $0$, excepted (possibly) at the point $0$ precisely. We now expand $r(\alpha ( \omega ) \, , \omega )$ as follows:
\[ r ( \alpha ( \omega ) \, , \omega ) = r(0 \, , \omega ) + \int_0^{\alpha ( \omega )} \frac{\partial r}{\partial \alpha} ( \tilde{\alpha} \, , \omega ) \, \text{d} \tilde{\alpha}. \]
Let us estimate $r(0 \, , \omega)$ first. We write
\[ r(0 \, , \omega) = e_1 \cdot \left [ \int_{\R} | P_\omega |^{\frac{1}{2}} (y) (1 + M_{0 , \omega})^{-1} (P_\omega^{\frac{1}{2}} e_1 ) (y) \text{d}y \right ] = e_1 \cdot \int_{\R} |P_\omega|^{\frac{1}{2}} (y) P_\omega^{\frac{1}{2}} (y) e_1 \, \text{d}y + \text{IR}_\omega = - \frac{I_\omega}{2} + \text{IR}_\omega \]
where $\displaystyle{\text{IR}_\omega = e_1 \cdot \left [ \int_{\R} | P_\omega |^{\frac{1}{2}} (y) ((1 + M_{0 , \omega})^{-1} - 1) (P_\omega^{\frac{1}{2}} e_1 ) (y) \text{d}y \right ]}$. We recall that $||M_{0 , \omega}|| \leqslant C \varepsilon_\omega < 1$. This shows, by Neumann expansion, that $||(1 + M_{0,\omega})^{-1} -1|| \leqslant 2 ||M_{0 , \omega}|| \leqslant C \varepsilon_\omega$. Then,
\[ \left | \int_{\R} | P_\omega |^{\frac{1}{2}} (y) ((1 + M_{0 , \omega})^{-1} - 1) (P_\omega^{\frac{1}{2}} e_1 ) (y) \text{d}y \right | \leqslant || \, |P_\omega|^{\frac{1}{2}} (y) || \, || (1 + M_{0 , \omega}^{-1} - 1) || \, ||P_\omega^{\frac{1}{2}} e_1|| \leqslant C \sqrt{\varepsilon_\omega} \cdot C \varepsilon_\omega \cdot C \sqrt{\varepsilon_\omega} \leqslant C \varepsilon_\omega^2. \]
Thus, we have proven that $r(0 \, , \omega) = - \frac{I_\omega}{2} + \mathcal{O} ( \varepsilon_\omega^2 )$. \\
\\ Now let us take care of the integral involving $\frac{\partial r}{\partial \alpha}$ in the Taylor expansion of $r ( \alpha ( \omega ) \, , \omega )$. We have
\[ \frac{\partial r}{\partial \alpha} ( \tilde{\alpha} \, , \omega ) = - e_1 \cdot \left [ \int_{\R} | P_\omega |^{\frac{1}{2}} (y) \left ( \left ( \frac{\partial M_{\tilde{\alpha},\omega}}{\partial \alpha} (1 + M_{\tilde{\alpha} , \omega})^{-2} \right ) (P_\omega^{\frac{1}{2}} e_1) \right ) (y) \, \text{d}y \right ]. \]
We notice that
\[ \frac{\partial N_{\tilde{\alpha} , \omega}}{\partial \alpha} (y \, , z) = \frac{1}{2} \left ( \begin{array}{cc} |y-z|^2 \theta ( \tilde{\alpha} |y-z| ) & 0 \\ 0 & - \frac{\tilde{\alpha} (1 + \kappa |y-z| )}{\kappa^3} e^{- \kappa |y-z|} \end{array} \right ) \, \, \, \, \, \, \text{where} \, \, \theta (w) = \frac{1-(1+w) e^{-w}}{w^2}. \]
We can show that $0 \leqslant \theta (w) \leqslant \frac{1}{2}$ for all $w \geqslant 0$, and thus $\left | \frac{\partial N_{\tilde{\alpha} , \omega}}{\partial \alpha} (y \, , z) \right | \leqslant C \left ( 1 + |y-z|^2 \right ) \leqslant C (1 + y^2 + z^2)$ in the sense that each coefficient of the matrix satisfies this inequality. Therefore, 
\[ \left | \frac{\partial M_{\tilde{\alpha} , \omega}}{\partial \alpha} (y \, , z) \right | = \left | P_\omega^{\frac{1}{2}} (y) \frac{\partial N_{\tilde{\alpha} , \omega}}{\partial \alpha} (y \, , z) |P_\omega|^{\frac{1}{2}} (z) \right | \leqslant C \varepsilon_\omega (1 + y^2 + z^2 ) Q_\omega (y) Q_\omega (z). \]
We recall that $||M_{\tilde{\alpha} , \omega} || \leqslant C \varepsilon_\omega < 1$ thus $||(1 + M_{\tilde{\alpha} , \omega})^{-2}|| \leqslant C$ and we have
\[ \begin{array}{rcl} \displaystyle{\left | \frac{\partial r}{\partial \alpha} ( \tilde{\alpha} \, , \omega ) \right |} & \leqslant & \displaystyle{\int_{\R} C \sqrt{\varepsilon_\omega} Q_\omega (y) \left | \int_{\R} \frac{\partial M_{\tilde{\alpha} , \omega}}{\partial \alpha} (y \, , z) \left ( (1 + M_{\tilde{\alpha} , \omega} )^{-2} (P_\omega^{\frac{1}{2}} e_1 ) \right ) (z) \, \text{d}z \right |} \\
\\ & \leqslant & \displaystyle{\int_{\R} C \sqrt{\varepsilon_\omega} Q_\omega (y) \left \| \frac{\partial M_{\tilde{\alpha} , \omega}}{\partial \alpha} (y \, , \cdot) \right \| \, \left \| (1 + M_{\tilde{\alpha} , \omega})^{-2} \right \| \, \left \| P_\omega^{\frac{1}{2}} e_2 \right \| \, \text{d}y} \\
\\ & \leqslant & \displaystyle{\int_{\R} C \sqrt{\varepsilon_\omega} Q_\omega (y) \cdot C \varepsilon_\omega (1+y^2) Q_\omega (y) \cdot C \cdot C \sqrt{\varepsilon_\omega} \, \text{d}y} \\
\\ & \leqslant & C \varepsilon_\omega^2. \end{array} \]
Getting back to the Taylor expansion, we get
\[ \left | \int_0^{\alpha ( \omega )} \frac{\partial r}{\partial \alpha} ( \tilde{\alpha} \, , \omega ) \, \text{d} \tilde{\alpha} \right | \leqslant C \varepsilon_\omega^2 \alpha ( \omega ) \leqslant C \varepsilon_\omega^2 \]
and then
\[ r(\alpha ( \omega ) \, , \omega ) = \frac{I_\omega}{2} + \mathcal{O} ( \varepsilon_\omega^2 ) \, \, \, \, \, \text{thus} \, \, \, \, \, \alpha ( \omega ) = - \frac{1}{2} r( \alpha ( \omega ) \, , \omega ) = \frac{I_\omega}{4} + \mathcal{O} ( \varepsilon_\omega^2 ) = \frac{I_\omega}{4} (1 + \mathcal{O} ( \varrho_\omega )) \]
where we recall that $\varrho_\omega \longrightarrow 0$ as $\omega \to 0$. \\
\\ \textit{(Expansion of the eigenfunctions.)} From now on, $\alpha$ denotes $\alpha ( \omega )$. We compute the expansion of the eigenfunction $X$ of \eqref{sysH} corresponding to the eigenfunction $\Psi = P_\omega^{\frac{1}{2}} X$ of \eqref{eqmu} chosen with the normalisation $p_\omega ( \Psi ) = -2 \alpha$. From \eqref{eqPsi} with $\mu = -1$ we obtain $\Psi = (1 + M_{\alpha , \omega})^{-1} ( P_\omega^{\frac{1}{2}} e_1 )$. This leads to
\[ X = e_1 - N_\alpha Y_\omega \, \, \, \, \, \text{with} \, \, Y_\omega = | P_\omega |^{\frac{1}{2}} \left ( 1 + A_{\alpha , \omega} \right ) (P_\omega^{\frac{1}{2}} e_1) \]
where $A_{\alpha , \omega} = (1 + M_{\alpha , \omega})^{-1} - 1 = \sum\limits_{j=1}^{+ \infty} (-1)^j M_{\alpha , \omega}^j$. Writing that $Y_\omega = P_\omega e_1 + | P_\omega |^{\frac{1}{2}} A_{\alpha , \omega} \left ( P_\omega^{\frac{1}{2}} e_1 \right )$ and using the expression of $N_0$, we see that
\[ N_0 Y_\omega = \left ( \begin{array}{c} -T_1 \\ -T_2 \end{array} \right ) + N_0 | P_\omega |^{\frac{1}{2}} A_{\alpha , \omega} \left ( P_\omega^{\frac{1}{2}} e_1 \right ) \]
where
\[ T_1 (y) = \frac{1}{4} \int_{\R} |y-z| \left ( a_\omega^+ (z) + a_\omega^- (z) \right ) \, \text{d}z \, \, \, \, \text{and} \, \, \, \, T_2 (y) = - \frac{\sqrt{2}}{8} \int_{\R} e^{- \sqrt{2} |y-z|} \left ( a_\omega^+ (z) - a_\omega^-(z) \right ) \, \text{d}z. \]
Hence, the expansion of $X$ can be written as
\[ X = e_1 + \left ( \begin{array}{c} X_1 \\ X_2 \end{array} \right ) + \widetilde{X} \, \, \, \, \, \text{where} \, \, \widetilde{X} = - N_0 | P_\omega |^{\frac{1}{2}} A_{\alpha , \omega} \left ( P_\omega^{\frac{1}{2}} e_1 \right ) + (N_0 - N_\alpha ) Y_\omega . \]
First, we have $|Y_\omega| \leqslant C |P_\omega|^{\frac{1}{2}} |P_\omega|^{\frac{1}{2}} e_1 \leqslant C ( |a_\omega^+| +  |a_\omega^-| ) \leqslant C \varepsilon_\omega Q_\omega^2$. Moreover, we have
\[ \left | e^{- \alpha |y-z|} -1+ \alpha |y-z| \right | \leqslant C \alpha^2 \left ( 1 + |y| + |z| \right )^2 \, \, \, \, \, \text{and} \, \, \, \, \, \left | \frac{e^{- \kappa |y-z|}}{\kappa} - \frac{e^{- \sqrt{2} |y-z|}}{\sqrt{2}} \right | \leqslant C \alpha^2 \]
thus $|N_0 - N_\alpha| (y \, , z) \leqslant C \alpha \left ( 1 + |y| + |z| \right )^2 \leqslant C I_\omega \left ( 1 + |y| + |z| \right )^2 \leqslant C \varepsilon_\omega (1 + |y| + |z|)^2$. Henceforth,
\[ | (N_\alpha - N_0 )Y_\omega (y)| \leqslant C \varepsilon_\omega^2 \int_{\R} (1 + |y| + |z| )^2 Q_\omega^2 (z) \, \text{d}z \leqslant C \varepsilon_\omega^2 \int_{\R} (1 + |y| + |z|)^2 e^{-2 |z|} \, \text{d}z \leqslant C \varepsilon_\omega^2 (1 +  |y|)^2 \leqslant C \varepsilon_\omega^2 (1+y^2). \]
Now, let us control the other term. We know that $||A_{\alpha , \omega}|| \leqslant C \varepsilon_\omega$ thus $||A_{\alpha , \omega} (P_\omega^{\frac{1}{2}} e_1)|| \leqslant C \varepsilon_\omega^{3/2}$. We also know that $|N_0 (y \, , z)| \leqslant 1 + |y| + |z|$. This leads, thanks to Cauchy-Schwarz inequality, to
\[ \left | N_0 |P_\omega|^{\frac{1}{2}} A_{\alpha , \omega} \left ( P_\omega^{\frac{1}{2}} e_1 \right ) (y) \right | \leqslant \left \| N_0 (y \, , \cdot) P_\omega^{\frac{1}{2}} \right \| \, \left \| A_{\alpha , \omega} (P_\omega^{\frac{1}{2}} e_1) \right \| \leqslant C \sqrt{\varepsilon_\omega} (1 + |y|) \cdot C \varepsilon_\omega^{3/2} \leqslant C \varepsilon_\omega^2 (1+y^2). \]
Therefore,
\[ X = \left ( \begin{array}{c} 1+T_1 \\ T_2 \end{array} \right ) + \widetilde{X} \, \, \, \, \, \text{where} \, \, | \widetilde{X} (y) | \leqslant C \varepsilon_\omega^2 (1 + y^2). \]
Defining $S_1 = T_1 + T_2$ and $S_2 = T_1 - T_2$, and recalling that $W_1 = X_1+X_2$ and $W_2 = X_1-X_2$, we obtain
\[ W = \left ( \begin{array}{c} 1 + S_1 \\ 1 + S_2 \end{array} \right ) + \widetilde{W} \, \, \, \, \, \text{where} \, \, | \widetilde{W} (y) | \leqslant C \varepsilon_\omega^2 (1 + y^2). \]
One can notice that the functions $T_j$ and $S_j$ verify the following differential equations: $T_1 '' = \frac{a_\omega^+ + a_\omega^-}{2}$, $T_2 '' - 2T_2 = \frac{a_\omega^+ - a_\omega^-}{2}$, $S_1'' = a_\omega^+ + 2T_2$ and $S_2'' = a_\omega^- - 2T_2$. Using $|a_\omega^{\pm}| \leqslant C \varepsilon_\omega Q_\omega^2$, we can also see that 
\[  |T_1(y)| + |T_2(y)| + |S_1(y)| + |S_2 (y)| \leqslant C \varepsilon_\omega (1 + |y|). \]
Using the differential equations satisfied by $T_j$ and $S_j$, we see that the previous bounds on $T_j$ and $S_j$ still hold for $T_j^{(k)}$ and $S_j^{(k)}$, for any $k \geqslant 0$. \\
\\ About the derivatives of $\widetilde{X}$, there is no additional difficulty. For example, 
\[ \frac{\partial (N_0 - N_\alpha)}{\partial y} (y \, , z) = \frac{1}{2} \left ( \begin{array}{cc} \text{sgn} (y-z) \left ( e^{- \alpha |y-z|} -1 \right ) & 0 \\ 0 & \displaystyle{\text{sgn} (y-z) \left ( e^{- \kappa |y-z|} - e^{- \sqrt{2}|y-z|} \right )} \end{array} \right ) \]
where $| e^{- \alpha |y-z|} -1 | \leqslant \alpha |y-z| \leqslant \alpha (1 + |y| + |z|)$ and $|e^{- \kappa |y-z|} - e^{- \sqrt{2} |y-z|}| \leqslant C |y-z| \, | \kappa - \sqrt{2}| \leqslant C \alpha^2 (1 + |y| + |z|)$. We could discuss the other derivatives with similar arguments and show that, for all $k \geqslant 1$,
\[ \left | \frac{\partial^k (N_0 - N_\alpha)}{\partial \alpha^k} \right | (y \, , z) \leqslant C \alpha (1 + |y| + |z|) \leqslant C \varepsilon_\omega (1 + |y| + |z|) \]
which leads to
\[ \left | \frac{\partial^k}{\partial y^k} \left ( (N_0 - N_\alpha) Y_\omega (y) \right ) \right | \leqslant \int_{\R} C \varepsilon_\omega (1 + |y| + |z|) \cdot \varepsilon_\omega Q_\omega^2 (z) \, \text{d}z \leqslant C \varepsilon_\omega^2 (1+ |y|). \]
About the other term in $\widetilde{X}$, we easily see that $\left | \frac{\partial^k N_0}{\partial y^k} (y \, , z) \right | \leqslant C$ for any $k \geqslant 1$, using the explicit expression of $N_0$. Using this estimate, we find that
\[ \left | \frac{\partial^k}{\partial y^k} \left [ N_0 |P_\omega|^{\frac{1}{2}} A_{\alpha , \omega} \left ( P_\omega^{\frac{1}{2}} e_1 \right ) \right ] \right | = \left | \int_{\R} \frac{\partial^k N_0}{\partial y^k} (y \, , z) |P_\omega|^{\frac{1}{2}} (z) A_{\alpha , \omega} \left ( P_\omega^{\frac{1}{2}} e_1 \right ) (z) \, \text{d}z \right | \leqslant C \varepsilon_\omega^2 (1 + |y| ). \]
Finally, for all $k \geqslant 1$, $| \widetilde{X}^{(k)} (y)| \leqslant C \varepsilon_\omega^2 (1 + |y|)$ and thus $| \widetilde{W}^{(k)} (y)| \leqslant C \varepsilon_\omega^2 (1+|y|)$. \\
\\ Now, let us deduce the expansions of $V_1$ and $V_2$. In what follows, the notation $\widetilde{\mathcal{O}}_p (1)$ refers to any function $w$ such that $|w^{(k)} (y)| \leqslant C_k (1+y^2)$ for all $k \in \{ 0 \, , ... \, , p \}$ (here the constants $C_k$ can depend on $k$ but do not depend on $\omega$). We recall that $V_1 = (S^*)^2 W_1 = \frac{Q_\omega ''}{Q_\omega} W_1 + 2 \frac{Q_\omega '}{Q_\omega} W_1 ' + W_1''$. Using the equations \eqref{eqQ} and \eqref{eqQ'2}, we compute
\[ \begin{array}{rl} & \displaystyle{ V_1 = 1-Q^2 + \widetilde{R}_1 + \varepsilon_\omega^2 \widetilde{\mathcal{O}}_4 (1)} \\ \\ \text{where} & \displaystyle{\widetilde{R}_1 = -2 Q D_\omega + (1-Q^2)S_1 + \frac{2Q'}{Q} S_1' - 2 \frac{G(\omega Q_\omega^2)}{\omega^2 Q_\omega^2} + 2T_2.} \end{array} \]
Here we recall that $D_\omega = Q_\omega - Q$. An elementary Taylor expansion shows that $\left | \frac{G ( \omega Q_\omega^2)}{\omega^2 Q_\omega^2} - \frac{G ( \omega Q^2)}{\omega^2 Q^2} \right | = \varepsilon_\omega^2 Q \, \widetilde{\mathcal{O}}_4 ( 1)$. Thus we can write $V_1 = 1-Q^2 + R_1 + \varepsilon_\omega^2 \widetilde{\mathcal{O}}_4 (1)$, where
\[ R_1 = -2 Q D_\omega + (1-Q^2)S_1 + \frac{2Q'}{Q} S_1' - 2 \frac{G(\omega Q^2)}{\omega^2 Q^2} + 2T_2. \]
We can easily observe that $|R_1^{(k)} (y)| \leqslant C \varepsilon_\omega (1 + |y|)$ for any $k \in \{ 0 \, , ... \, , 5 \}$. \\
\\ In order to establish the expansion of $V_2 = \lambda^{-1} L_+ V_1$, we recall that $\alpha = \mathcal{O} ( I_\omega ) = \mathcal{O} ( \varepsilon_\omega )$, thus $\lambda^{-1} = (1 - \alpha^2)^{-1} = 1 + \mathcal{O} ( \varepsilon_\omega^2 )$. After computations and using the equations \eqref{eqQ} and \eqref{eqQ'2}, we obtain
\[ \begin{array}{rl} & V_2 = 1 + \widetilde{R}_2 + \varepsilon_\omega^2 \widetilde{\mathcal{O}}_2 (1) \\
\\ \text{where} & \displaystyle{\widetilde{R}_2 = -R_1'' + R_1 - 3Q^2 R_1 - 6Q(1-Q^2)D_\omega - (1-Q^2) \left ( \frac{g( \omega Q^2)}{\omega} + 2 Q^2 g'(\omega Q^2) \right ).} \end{array}  \]
Here too, elementary Taylor expansions have allowed us to write $g(\omega Q^2)$ instead of $g(\omega Q_\omega^2)$ and $Q^2 g'(\omega Q^2)$ instead of $Q_\omega^2 g'(\omega Q_\omega^2)$. The cost is absorbed in the term $\varepsilon_\omega^2 \widetilde{\mathcal{O}}_2 (1)$. Now we compute $R_1''$ using the expression above. We find, after lengthy computations, $\widetilde{R}_2 = R_2 + \varepsilon_\omega^2 \widetilde{\mathcal{O}}_2 (1)$, where
\[ R_2 = -4(1-Q^2) QD_\omega + 4Q' D_\omega ' + T_1 -3 T_2 + \frac{2Q'}{Q} T_1' - \frac{2Q'}{Q} T_2' + 2 \frac{g ( \omega Q^2)}{\omega} -4 \frac{G ( \omega Q^2)}{\omega^2 Q^2} + 2 \frac{G( \omega Q^2)}{\omega^2} . \]
At last, we have the following expansion of $V_2$: $V_2 = 1 + R_2 + \varepsilon_\omega^2 \widetilde{\mathcal{O}}_2 (1)$. Moreover, $|R_2^{(k)} (y)| \leqslant C \varepsilon_\omega (1 + |y|)$ for any $k \in \{ 0 \, , ... \, , 4 \}$. \\
\\ \textit{(Decay properties.)} It is clear that $|X_j| \leqslant 1 + |T_j| + C \varepsilon_\omega^2 (1+y^2) \leqslant C (1 + y^2)$, thus $|P_\omega X| \leqslant C \varepsilon_\omega Q_\omega^2 (y) (1 + y^2) \leqslant C \varepsilon_\omega e^{-|y|}$. Recalling \eqref{sysH} and the expression of $H_\alpha^{-1} (y \, , z)$, we see that
\[ \begin{array}{rl} & \displaystyle{X = - H_\alpha^{-1} (P_\omega X) = - \frac{1}{2} \int_{\R} \left ( \begin{array}{cc} \frac{e^{- \alpha |y-z|}}{\alpha} & 0 \\ 0 & \frac{e^{- \kappa |y-z|}}{\kappa} \end{array} \right ) P_\omega (z) X(z) \, \text{d}z} \\
\\ \text{therefore} & \displaystyle{|X_1(y)| \leqslant \frac{C \varepsilon_\omega}{\alpha} \int_{\R} e^{- \alpha |y-z|} e^{-|z|} \, \text{d}z \leqslant \frac{C \varepsilon_\omega}{\alpha} \, e^{-\alpha |y|} \leqslant \frac{C \varepsilon_\omega}{I_\omega} e^{-\alpha |y|}.} \end{array} \]
This does not prove the estimate $|X_1(y)| \leqslant C e^{-\alpha |y|}$ that we claim to be true in Proposition 2, since $\varepsilon_\omega / I_\omega$ has no reason to be bounded: this particular estimate will have to wait a little more. However, the argument above proves similarly that
\[ |X_2 (y)| \leqslant \frac{C \varepsilon_\omega}{\kappa} \, e^{- \kappa |y|} \leqslant C \varepsilon_\omega e^{- \kappa |y|} \leqslant C \varepsilon_\omega e^{-|y|}  \]
since we can assume that $\omega$ is small enough so that $\kappa > 1$. The argument is the same in order to estimate $X_j^{(k)}$ for $k \geqslant 1$. Indeed,
\[ X^{(k)} (y) = - \frac{1}{2} \int_{\R} \frac{\partial^k}{\partial y^k}\left ( \begin{array}{cc} \frac{e^{- \alpha |y-z|}}{\alpha} & 0 \\ 0 & \frac{e^{- \kappa |y-z|}}{\kappa} \end{array} \right ) P_\omega (z) X(z) \, \text{d}z. \]
Up to constants, differentiating the matrix $H_\alpha^{-1}$ only makes an $\alpha^k$ and a $\kappa^k$ appear (respectively in the control of $X_1^{(k)}$ and in the control of $X_2^{(k)}$). That way, we obtain, for any $k \geqslant 1$,
\[ \begin{array}{rl} & |X_1^{(k)} (y)| \leqslant C_k \varepsilon_\omega \alpha^{k-1} e^{- \alpha |y|} \leqslant C_k \varepsilon_\omega I_\omega^{k-1} e^{- \alpha |y|} \\ \\ \text{and} & |X_2^{(k)} (y) | \leqslant C_k \kappa^k \varepsilon_\omega e^{- \kappa |y|} \leqslant C_k \varepsilon_\omega e^{-|y|} \end{array} \]
and the bound on $X_2^{(k)}$ remains true when $k=0$. This proves the decay properties of $W_j^{(k)}$ for any $k \geqslant 1$, and it also proves the decay properties of $(W_1 - W_2)^{(k)}$ for any $k \geqslant 0$, since $W_1 - W_2 = 2X_2$. \\
\\ The similar bounds on the functions $V_j$ do not present additional difficulties, they simply stem from the expressions $V_1 = (S^*)^2 W_1$ and $V_2 = \lambda^{-1} L_+ V_1$. We need the estimate $|W_j(y)| \leqslant C e^{- \alpha |y|}$, which is not proven yet but which will be proven in the following paragraph. \\
\\ \textit{(Asymptotics of the eigenfunctions.)} The equality \eqref{sysH} can be written as
\[ \left \{ \begin{array}{ccl} X_1 (y) &=& \displaystyle{- \frac{1}{4 \alpha} \int_{\R} e^{- \alpha |y-z|} \left [ a_\omega^+ (z) ( X_1 (z) + X_2(z)) + a_\omega^- (z) (X_1 (z) - X_2 (z)) \right ] \, \text{d}z} \\
\\ X_2 (y) &=& \displaystyle{- \frac{1}{4 \kappa} \int_{\R} e^{- \kappa |y-z|} \left [ a_\omega^+ (z) ( X_1 (z) + X_2(z)) - a_\omega^- (z) (X_1 (z) - X_2 (z)) \right ] \, \text{d}z.} \end{array} \right. \]
Let us write
\[ \begin{array}{rcl} |X_1 (y) - e^{- \alpha |y|}| & \leqslant & \displaystyle{\left | - \frac{1}{4 \alpha} \int_{\R} e^{- \alpha |y-z|} \left [ a_\omega^+ (z) (X_1 (z) + X_2(z) -1) + a_\omega^- (z) (X_1(z) - X_2(z) - 1) \right ] \, \text{d}z \right | } \\
\\ & & \displaystyle{\, \, \, \, \, \, \, + \, \, \left | - \frac{1}{4 \alpha} \int_{\R} e^{- \alpha |y-z|} \left ( a_\omega^+ (z) + a_\omega^- (z) \right ) \, \text{d}z + \frac{1}{4 \alpha} \, e^{- \alpha |y|} \int_{\R} \left ( a_\omega^+ (z) + a_\omega^- (z) \right ) \, \text{d}z \right |} \\
\\ & & \displaystyle{\, \, \, \, \, \, \, + \, \, \left | - \frac{1}{4 \alpha} \, e^{- \alpha |y|} \int_{\R} \left ( a_\omega^+ (z) + a_\omega^- (z) \right ) \, \text{d}z - e^{- \alpha |y|} \right |.} \end{array} \]
Let us estimate these three terms separately. For the first one, we recall that $|X_1(z)-1| = |T_1(z) + \widetilde{X}_1 (z)| \leqslant C \varepsilon_\omega (1+z^2)$ and $|X_2(z)| = |T_2(z) + \widetilde{X}_2(z)| \leqslant C \varepsilon_\omega (1+z^2)$. Thus,
\[ \begin{array}{rl} & \displaystyle{\left | - \frac{1}{4 \alpha} \int_{\R} e^{- \alpha |y-z|} \left [ a_\omega^+ (z) (X_1 (z) + X_2(z) -1) + a_\omega^- (z) (X_1(z) - X_2(z) - 1) \right ] \, \text{d}z \right |} \\
\\ \leqslant & \displaystyle{\frac{C}{\alpha} \int_{\R} e^{- \alpha |y-z|} \varepsilon_\omega Q_\omega^2 (z) \varepsilon_\omega (1+z^2) \, \text{d}z \, \, \leqslant \, \, \frac{C \varepsilon_\omega^2}{\alpha} \int_{\R} e^{- \alpha |y-z| -2 |z|} (1+z^2) \, \text{d}z} \\
\\ \leqslant & \displaystyle{\frac{C \varepsilon_\omega^2}{\alpha} \, e^{- \alpha |y|} \, \, \leqslant \, \, C \varrho_\omega e^{- \alpha |y|}.} \end{array} \]
For the second term, we use the inequalities $|e^{-w}-1| \leqslant |w| e^{|w|}$ and $||y-z|-|y|| \leqslant |z|$ that hold for all $w,y,z \in \R$, as well as the monotonicity of $w \mapsto we^w$ on $[0 \, , + \infty )$. We find that
\[ \begin{array}{rl} & \displaystyle{\left | - \frac{1}{4 \alpha} \int_{\R} e^{- \alpha |y-z|} \left ( a_\omega^+ (z) + a_\omega^- (z) \right ) \, \text{d}z + \frac{1}{4 \alpha} \, e^{- \alpha |y|} \int_{\R} \left ( a_\omega^+ (z) + a_\omega^- (z) \right ) \, \text{d}z \right |} \\
\\ \leqslant & \displaystyle{\frac{C}{\alpha} \, e^{- \alpha |y|} \int_{\R} \left | e^{- \alpha (|y-z|-|y|)} -1 \right | \varepsilon_\omega Q_\omega^2 (z) \, \text{d}z \, \, \leqslant \, \, \frac{C \varepsilon_\omega}{\alpha} \, e^{- \alpha |y|} \int_{\R} \alpha ||y-z|-|y|| e^{\alpha ||y-z|-|y||} Q_\omega^2 (z) \, \text{d}z} \\
\\ \leqslant & \displaystyle{C \varepsilon_\omega e^{- \alpha |y|} \int_{\R} |z| e^{\alpha |z|} e^{-2|z|} \, \text{d}z \, \, \leqslant \, \, C \varepsilon_\omega e^{- \alpha |y|}.} \end{array} \]
At last, the final term is controlled as follows:
\[ \left | - \frac{1}{4 \alpha} \, e^{- \alpha |y|} \int_{\R} \left ( a_\omega^+ (z) + a_\omega^- (z) \right ) \, \text{d}z - e^{- \alpha |y|} \right | = \left | \frac{I_\omega}{4 \alpha} - 1 \right | e^{- \alpha |y|} \leqslant C \varrho_\omega e^{- \alpha |y|} \]
since we recall that $\alpha = \frac{I_\omega}{4} \left ( 1 + \mathcal{O} ( \varrho_\omega ) \right )$. Gathering all these estimates, we obtain:
\[ |X_1 (y) - e^{- \alpha |y|} | \leqslant C \varrho_\omega e^{- \alpha |y|}. \]
We have already proven previously that $|X_2(y)| \leqslant C \varepsilon_\omega e^{-|y|} \leqslant C \varepsilon_\omega e^{- \alpha |y|}$, therefore we get the desired estimate:
\[ |W_j (y) - e^{- \alpha |y|}| \leqslant C \varrho_\omega e^{- \alpha |y|}. \]
This bound enables us to prove the decay property of $W_j$ that we have not proven yet. Indeed, taking $\omega >0$ small enough, we get $|W_j(y)| \leqslant e^{- \alpha |y|} + C \varrho_\omega e^{- \alpha |y|} \leqslant C e^{- \alpha |y|}$. \\
\\ The estimates for $V_1$ and $V_2$ follow from the expressions $V_1 = (S^*)^2W_1$ and $V_2 = \lambda^{-1} L_+ V_1$. Let us give a little more details. 
\[ \begin{array}{rcl} |V_1 - (1-Q^2) e^{-\alpha|y|}| & \leqslant & \displaystyle{\left | \left ( 1 - Q_\omega^2 + \frac{g( \omega Q_\omega^2 )}{\omega} \right ) W_1 - (1-Q^2) e^{-\alpha |y|} \right | + \left | \frac{2 Q_\omega '}{Q_\omega} \right | \, |W_1'| + |W_1''|} \\
\\ & \leqslant & \displaystyle{|(1-Q_\omega^2) (W_1-e^{-\alpha |y|})| + \left | ((1-Q_\omega^2)-(1-Q^2)) e^{-\alpha |y|} \right | + \left | \frac{g( \omega Q_\omega^2 )}{\omega} \right | \, |W_1|} \\ & & \, \, \, \, \, \, \, \, \, \, \, \displaystyle{+ \, \left | \frac{2 Q_\omega '}{Q_\omega} \right | \, |W_1'| + |W_1''|} \\
\\ & \leqslant & \displaystyle{C ( \varepsilon_\omega + \varrho_\omega ) e^{- \alpha |y|}} \end{array} \]
after analysing each term. To control $V_2$ we first see that
\[ L_+ (1-Q_\omega^2) = 1 + r_\omega^1 \, \, \, \, \, \text{where} \, \, r_\omega^1 = (1-3Q^2) \frac{g(\omega Q_\omega^2)}{\omega} + 2(1-Q_\omega^2) Q_\omega^2 g'( \omega Q_\omega^2) -2 \frac{G ( \omega Q_\omega^2 )}{\omega^2 Q_\omega^2}. \]
In what follows, let us denote $g_\omega := \frac{g( \omega Q_\omega^2 )}{\omega}$, $G_\omega := \frac{G( \omega Q_\omega^2)}{\omega^2 Q_\omega^2}$, $dg_\omega := Q_\omega^2 g'(\omega Q_\omega^2)$ and $\xi_Q := \frac{Q_\omega '}{Q_\omega}$. After computations, we find that $L_+ V_1 = W_1 + r_\omega^2$ where
\[ \begin{array}{rcl} r_\omega^2 &=& r_\omega^1 W_1 + 4 Q_\omega Q_\omega ' W_1 ' - (1-Q_\omega^2) W_1'' - g_\omega '' W_1 - 2 g_\omega ' W_1 ' - g_\omega W_1 '' + (1-3Q_\omega^2 + g_\omega + 2dg_\omega ) g_\omega W_1 \\ \\ & & \, \, \, \, \, - 2 \xi_Q '' W_1 ' - 4 \xi_Q ' W_1 '' + 2 \xi_Q W_1 ''' + 2 (1 -3Q_\omega^2 + g_\omega + 2 dg_\omega ) \xi_Q W_1 ' - W_1 '''' + (1-3Q_\omega^2 + g_\omega + 2 dg_\omega ) W_1 ''. \end{array} \]
Using the estimates $|W_1 (y) - e^{- \alpha |y|} | \leqslant C \varrho_\omega e^{- \alpha |y|}$, $|r_\omega^2 (y)| \leqslant C \varepsilon_\omega e^{- \alpha |y|}$ and $\lambda^{-1} = 1 + \mathcal{O} ( \varrho_\omega )$, we finally obtain the desired estimate:
\[ |V_2 (y) - e^{- \alpha |y|} | \leqslant C \varrho_\omega e^{- \alpha |y|}. \]
The estimates for $| \langle W_1 \, , W_2 \rangle - \alpha^{-1}|$ and $| \langle V_1 \, , V_2 \rangle - \alpha^{-1}|$ follow easily by integration. Let us prove the second one for example. We use the estimates $|V_1 (y) - (1-Q^2) e^{- \alpha |y|}| \leqslant C \varrho_\omega e^{- \alpha |y|}$ and $|V_2 (y) - e^{- \alpha |y|} | \leqslant C \varrho_\omega e^{- \alpha |y|}$. From these estimates we see that $|V_1(y)| + |V_2 (y)| \leqslant C e^{- \alpha |y|}$. Then,
\[ \begin{array}{rl} & \displaystyle{\left | \langle V_1 \, , V_2 \rangle - \frac{1}{\alpha} \right | \, \, = \, \, \left | \int_{\R} \left ( V_1 V_2 - e^{-2 \alpha |y|} \right ) \, \text{d}y \right |} \\
\\ \leqslant & \displaystyle{\left | \int_{\R} V_1 (V_2 - e^{- \alpha |y|} ) \, \text{d}y \right | + \left | \int_{\R} e^{- \alpha |y|} (V_1 - (1-Q^2) e^{- \alpha |y|}) \, \text{d}y \right | + \left | \int_{\R} - e^{- \alpha |y|} Q^2(y) e^{- \alpha |y|} \, \text{d}y \right | } \\
\\ \leqslant & \displaystyle{C \varrho_\omega \int_{\R} e^{-2 \alpha |y|} \, \text{d}y + C \varrho_\omega \int_{\R} e^{-2 \alpha |y|} \, \text{d}y + C} \\
\\ \leqslant & \displaystyle{C \, \frac{\varrho_\omega}{\alpha} + C \, \, \leqslant \, \, C \, \frac{\varrho_\omega}{I_\omega} + C \, \, \leqslant \, \, C \, \frac{\varrho_\omega}{I_\omega},} \end{array} \]
the last inequality being true because $I_\omega \leqslant C \varepsilon_\omega \leqslant C \varrho_\omega$. \\
\\ The last estimate, about $W_j'$, is proven similarly to the one on $W_j$. Take $y>0$. We write that $|X_1'(y) + \alpha e^{- \alpha y}| \leqslant \mathbb{T}_1 + \mathbb{T}_2 + \mathbb{T}_3$ where
\[ \begin{array}{rl} & \displaystyle{\mathbb{T}_1 = \left | \frac{1}{4} \int_{\R} \text{sgn} (y-z) e^{- \alpha |y-z|} \left | a_\omega^+ (X_1 + X_2-1) + a_\omega^- (X_1 - X_2 - 1) \right ] (z) \, \text{d}z \right |,} \\
\\ & \displaystyle{\mathbb{T}_2 = \left | \frac{1}{4} \int_{\R} \text{sgn} (y-z) e^{- \alpha |y-z|} (a_\omega^+ + a_\omega^-)(z) - \frac{1}{4} e^{- \alpha y} \int_{\R} (a_\omega^+ + a_\omega^-)(z) \, \text{d}z \right |} \\
\\ \text{and} & \displaystyle{\mathbb{T}_3 = \left | \frac{1}{4} e^{- \alpha y} \int_{\R} (a_\omega^+ + a_\omega^-) (z) \, \text{d}z + \alpha e^{- \alpha y} \right |.} \end{array} \]
As previously, we see that $\mathbb{T}_1 \leqslant C \varepsilon_\omega^2 e^{- \alpha y} = C \varrho_\omega I_\omega e^{- \alpha y}$, $\mathbb{T}_2 \leqslant C \varepsilon_\omega I_\omega e^{- \alpha y} \leqslant C \varrho_\omega I_\omega e^{- \alpha y}$ and $\mathbb{T}_3 \leqslant C \varrho_\omega I_\omega e^{- \alpha y}$. Hence,
\[ |X_1'(y) + \alpha e^{- \alpha y}| \leqslant C \varrho_\omega I_\omega e^{- \alpha y}. \]
Recalling that $|X_2'(y)| \leqslant C \varepsilon_\omega e^{- \kappa y}$, we obtain the desired estimate on $|W_j' + \alpha e^{- \alpha y}|$. Since $W_j'$ is odd, a similar bound holds for $|W_j' - \alpha e^{- \alpha |y|}|$ for $y \leqslant 0$. \\
\\ \textit{(Derivatives with regards to $\omega$.)} The last estimates in Proposition 2 will require many more calculations. First, from Lemma 5 in \cite{Ri}, we know that $| \partial_\omega Q_\omega (y)| \leqslant \frac{C}{\omega} (1 +  |y|) e^{-|y|}$. Therefore, 
\begin{equation}
    | \partial_\omega a_\omega^{\pm}| \leqslant \frac{C \varepsilon_\omega}{\omega} (1 + |y|) e^{-2|y|}.
    \label{daomega}
\end{equation}
From the expression $M_{\alpha,\omega} = P_\omega^{\frac{1}{2}} N_\alpha P_\omega|^{\frac{1}{2}}$, we compute $M_{\alpha,\omega}^k = P_\omega^{\frac{1}{2}} (N_\alpha P_\omega)^{k-1} N_\alpha |P_\omega|^{\frac{1}{2}}$ for all $k \geqslant 1$. This leads to the following expression:
\[ r( \tilde{\alpha} \, , \omega ) = e_1 \cdot \int_{\R} P_\omega (1 + N_{\tilde{\alpha}} P_\omega)^{-1} e_1. \]
Thus, 
\[ \frac{\partial r}{\partial \omega} ( \tilde{\alpha} \, , \omega) = e_1 \cdot \left ( \int_{\R} (\partial_\omega P_\omega) (1 + N_\alpha P_\omega)^{-1} e_1 - \int_{\R} P_\omega N_\alpha (\partial_\omega P_\omega) (1 + N_\alpha P_\omega)^{-2} e_1 \right ). \]
Wee know from \eqref{daomega} that $| \partial_\omega P_\omega | \leqslant \frac{C \varepsilon_\omega}{\omega} (1 + |y|)e^{-2|y|}$. Using the estimates $|N_{\tilde{\alpha}} (y \, , z) | \leqslant C ( 1 + |y-z|)$ and $|P_\omega (z)| \leqslant C \varepsilon_\omega Q_\omega^2 (z) \leqslant C \varepsilon_\omega e^{-2|z|}$, we easily show by induction that
\[ \forall k \in \N , \, \, \, |(N_{\tilde{\alpha}} P_\omega)^k e_1 (y)| \leqslant C \varepsilon_\omega^k (1 + |y|). \]
Using Neumann expansion, we obtain that $|(1 + N_{\tilde{\alpha}} P_\omega )^{-1} e_1(y)| \leqslant C (1 + |y|)$ and then
\[ \left | \int_{\R} ( \partial_\omega P_\omega ) (y) \left ( (1 + N_{\tilde{\alpha}} P_\omega)^{-1} e_1 \right ) (y) \, \text{d}y \right | \leqslant \frac{C \varepsilon_\omega}{\omega}. \]
The second term in $\partial_\omega r$ is treated similarly: we see that $|(1 + N_{\tilde{\alpha}} P_\omega)^{-2} e_1(y)| \leqslant C (1 + |y|)$ and then
\[ \left | \int_{\R} P_\omega N_{\tilde{\alpha}} (\partial_\omega P_\omega) (1 + N_{\tilde{\alpha}} P_\omega)^{-2} e_1 \right | \leqslant \frac{C \varepsilon_\omega^2}{\omega}. \]
Combining these estimates, we obtain $| \partial_\omega r ( \tilde{\alpha} \, , \omega ) | \leqslant C \varepsilon_\omega / \omega$. Note that here $\tilde{\alpha}$ does not depend on $\omega$; this estimate is proven regardless of $\tilde{\alpha}$. Now we deduce the control of $\alpha ' ( \omega )$. Indeed, we recall that $\alpha ( \omega ) = - \frac{1}{2} r \left ( \alpha ( \omega ) \, , \omega \right )$. Thus,
\[ | \alpha ' ( \omega ) | = \left | \frac{- \partial_\omega r ( \alpha ( \omega ) \, , \omega )}{2 + \partial_\alpha r ( \alpha ( \omega ) \, , \omega )} \right | = \frac{|\partial_\omega r ( \alpha ( \omega ) \, , \omega )|}{|2 + \mathcal{O}( \varepsilon_\omega^2)|} \leqslant | \partial_\omega r( \alpha ( \omega ) \, , \omega )| \leqslant \frac{C \varepsilon_\omega}{\omega}, \]
which is the first result announced. Now, to control the difference $\alpha ' (\omega ) - \frac{1}{4} \partial_\omega I_\omega$, let us write that
\[ \alpha ' ( \omega ) - \frac{1}{4} \partial_\omega I_\omega  = \frac{-1/2}{1 + \frac{1}{2} \partial_\alpha r ( \alpha ( \omega ) \, , \omega )} \left ( \partial_\omega r ( \alpha ( \omega ) \, , \omega ) + \frac{1}{2} \partial_\omega I_\omega \right ) + \frac{\partial_\omega I_\omega}{4} \left ( \frac{1}{1 + \frac{1}{2} \partial_\alpha r ( \alpha ( \omega ) \, , \omega )} -1 \right ). \]
Since $\frac{\partial_\alpha r ( \alpha ( \omega ) \, , \omega )}{2} = \mathcal{O} ( \varepsilon_\omega )$, $\left | 1 + \frac{\partial_\alpha r ( \alpha ( \omega ) \, , \omega )}{2} \right |^{-1} \leqslant C$. Now, we recall that
\[ \partial_\omega r ( \alpha ( \omega ) \, , \omega ) = e_1 \cdot \int_{\R} (\partial_\omega P_\omega ) (1 + N_\alpha P_\omega)^{-1} e_1 - e_1 \cdot \int_{\R} P_\omega N_\alpha ( \partial_\omega P_\omega ) (1 + N_\alpha P_\omega)^{-2} e_1 \]
where we have already proven that
\[ \left | e_1 \cdot \int_{\R} P_\omega N_\alpha ( \partial_\omega P_\omega ) (1 + N_\alpha P_\omega)^{-2} e_1 \right | \leqslant \frac{C \varepsilon_\omega^2}{\omega}. \]
On the other hand,
\[ \left | e_1 \cdot ( \partial_\omega P_\omega )(1 + N_\alpha P_\omega)^{-1} e_1 + \frac{\partial_\omega I_\omega}{2} \right | = \left | e_1 \cdot \int_{\R} \partial_\omega P_\omega ((1 + N_\alpha P_\omega)^{-1} -1) e_1 \right | \leqslant \frac{C \varepsilon_\omega^2}{\omega} \]
thanks to the estimate $|((1 + N_\alpha P_\omega)^{-1} -1) e_1 (y)| \leqslant C \varepsilon_\omega (1+ |y|)$, established as previously thanks to Neumann expansion. Gathering these estimates, we obtain
\[ \left | \partial_\omega r ( \alpha ( \omega ) \, , \omega ) + \frac{1}{2} \partial_\omega I_\omega \right | \leqslant \frac{C \varepsilon_\omega^2}{\omega}. \]
We conclude by noticing that $| \partial_\omega I_\omega | \leqslant \frac{C \varepsilon_\omega}{\omega}$ and $\left | \frac{1}{1 + \frac{1}{2} \partial_\alpha r ( \alpha ( \omega ) \, , \omega )} -1 \right | \leqslant C \varepsilon_\omega$. This leads to the desired estimate:
\[ \left | \alpha' ( \omega ) - \frac{\partial_\omega I_\omega}{4} \right | \leqslant \frac{C \varepsilon_\omega^2}{\omega}. \]
Now, we control the terms $X_j$. To do so, we have to control first the terms $\widetilde{X}_j$. In what follows, $\alpha$ denotes $\alpha ( \omega )$. We recall that $\widetilde{X} = (N_0 - N_\alpha) Y_\omega - N_0 P_\omega ((1 + N_\alpha P_\omega)^{-1} -1) e_1$ where $Y_\omega = P_\omega (1 + N_\alpha P_\omega)^{-1} e_1$. Thus,
\[ \begin{array}{rcl} \partial_\omega \widetilde{X} &=& - \alpha '(\omega) (\partial_\alpha N_\alpha) Y_\omega + (N_0 - N_\alpha) \partial_\omega Y_\omega - N_0 ( \partial_\omega P_\omega ) \left ( (1 + N_\alpha P_\omega)^{-1} -1 \right ) e_1 \\
\\ & & \, \, \, \, \, \, \, \, + \, N_0 P_\omega \left ( \alpha ' ( \omega ) (\partial_\alpha N_\alpha ) P_\omega + N_\alpha \partial_\omega P_\omega \right ) (1 + N_\alpha P_\omega)^{-2} e_1. \end{array} \]
We recall that $| \partial_\alpha N_\alpha | \leqslant C (1 + y^2 + z^2)$. We also recall that $|(1 + N_\alpha P_\omega)^{-1} e_1(y)| \leqslant C (1 + |y|)$, which implies that $|Y_\omega (y)| \leqslant C \varepsilon_\omega (1 + |y|) Q_\omega^2 (y)$. Using all the previous bounds, we find successively that
\[ \begin{array}{rl} & \left | \alpha '(\omega) (\partial_\alpha N_\alpha) Y_\omega \right | \leqslant \frac{C \varepsilon_\omega^2}{\omega} (1 + y^2), \\
\\ & \left | N_0 P_\omega \alpha ' ( \omega ) (\partial_\alpha N_\alpha) P_\omega (1 + N_\alpha P_\omega)^{-2} e_1 \right | \leqslant \frac{C \varepsilon_\omega^3}{\omega} (1 +  |y|), \\
\\ & \left | N_0 P_\omega N_\alpha ( \partial_\omega P_\omega ) (1 + N_\alpha P_\omega)^{-2} e_1 \right | \leqslant \frac{C \varepsilon_\omega^2}{\omega} (1 + |y|) \\
\\ \text{and} & \left | N_0 ( \partial_\omega P_\omega ) \left ( (1 + N_\alpha P_\omega)^{-1} -1 \right )e_1 \right | \leqslant \frac{C \varepsilon_\omega^2}{\omega} (1 + |y|). \end{array} \]
To establish the last estimate, we have used the fact that $\left | \left ( (1 + N_\alpha P_\omega)^{-1} -1 \right ) e_1(z) \right | \leqslant C \varepsilon_\omega (1 + |z|)$, which is shown by Neumann expansion (as previously). The remaining term to be estimated is $\partial_\omega Y_\omega$. We have
\[ \partial_\omega Y_\omega = (\partial_\omega P_\omega ) (1 + N_\alpha P_\omega)^{-1} e_1 - P_\omega ( \alpha ' ( \omega ) ( \partial_\alpha N_\alpha ) P_\omega + N_\alpha (\partial_\omega P_\omega ) ) ( 1 + N_\alpha P_\omega )^{-2} e_1. \]
Controlling each one of these terms as we did previously, we obtain the following estimate:
\[ | \partial_\omega Y_\omega | \leqslant \frac{C \varepsilon_\omega}{\omega} (1 + y^2) Q_\omega^2 (y). \]
Recalling that $|(N_0 - N_\alpha)(y \, , z)| \leqslant C I_\omega (1 + y^2 + z^2)$, these estimates lead to
\[ \left | (N_0 - N_\alpha) \partial_\omega Y_\omega \right | \leqslant \frac{C \varepsilon_\omega I_\omega}{\omega} (1 + y^2). \]
Gathering all these estimates, we have proven that
\[ | \partial_\omega \widetilde{X} (y)| \leqslant \frac{C \varepsilon_\omega^2}{\omega} (1 + y^2). \]
In order to estimate $\partial_\omega X$, we recall that $X_1 = 1 +T_1 + \widetilde{X}_1$ and $X_2 = T_2 + \widetilde{X}_2$. Thanks to the explicit expressions of $T_1$ and $T_2$, and using \eqref{daomega}, we see that $| \partial_\omega T_1 | \leqslant \frac{C \varepsilon_\omega}{\omega} (1 + |y|)$ and $| \partial_\omega T_2 | \leqslant \frac{C \varepsilon_\omega}{\omega}$. This leads to
\[ | \partial_\omega X(y)| \leqslant \frac{C \varepsilon_\omega}{\omega} (1 + y^2). \]
Now, the proof resembles the one of the asymptotics of the eigenfunctions. We write that $-X_1 = e_1 \cdot H_\alpha^{-1} (P_\omega X) = \mathbf{T}_1 + \mathbf{T}_2 + \mathbf{T}_3$ where
\[ \begin{array}{rl} & \displaystyle{\mathbf{T}_1 = \frac{1}{2 \alpha} \int_{\R} e^{- \alpha |y-z|} (b_\omega^+ (X_1-1) + b_\omega^- X_2) (z) \, \text{d}z}, \\
\\ & \displaystyle{\mathbf{T}_2 = \frac{1}{2 \alpha} \int_{\R} (e^{- \alpha |y-z|} - e^{- \alpha |y|}) b_\omega^+ (z) \, \text{d}z} \\
\\ \text{and} & \displaystyle{\mathbf{T}_3 = \frac{e^{- \alpha |y|}}{2 \alpha} \int_{\R} b_\omega^+ (z) \, \text{d}z = - \frac{I_\omega}{4 \alpha} \, e^{- \alpha |y|}.} \end{array} \]
Using the estimate on $|4 \alpha '( \omega ) - \partial_\omega I_\omega |$, we establish that $\left | \partial_\omega \left ( \frac{I_\omega}{4 \alpha} \right ) \right | \leqslant \frac{C \varepsilon_\omega \varrho_\omega}{\omega \alpha}$. This leads to
\[ | \partial_\omega \mathbf{T}_3 | \leqslant \frac{C \varepsilon_\omega \varrho_\omega}{\omega \alpha} (1 + |y|) e^{- \alpha |y|}. \]
Now, we write $\partial_\omega \mathbf{T}_1 = \mathbf{T}_{1A} + \mathbf{T}_{1B} + \mathbf{T}_{1C}$, where
\[  \begin{array}{rl} & \displaystyle{\mathbf{T}_{1A} = - \frac{\alpha ' ( \omega )}{2 \alpha^2} \int_{\R} (1 + \alpha |y-z|) e^{- \alpha |y-z|} \left ( b_\omega^+ (X_1-1) + b_\omega^- X_2 \right ) (z) \, \text{d}z}, \\
\\ & \displaystyle{\mathbf{T}_{1B} = \int_{\R} \frac{e^{- \alpha |y-z|}}{2 \alpha} \left ( (\partial_\omega b_\omega^+) (X_1-1) + (\partial_\omega b_\omega^-)X_2 \right ) (z) \, \text{d}z} \\
\\ \text{and} & \displaystyle{\mathbf{T}_3 = \int_{\R} \frac{e^{- \alpha |y-z|}}{2 \alpha} \left ( b_\omega^+ \partial_\omega X_1 + b_\omega^- \partial_\omega X_2 \right ) (z) \, \text{d}z}. \end{array} \]
Using the previous known estimates, including $|X_1-1| \leqslant C \varepsilon_\omega (1+z^2)$, $|X_2| \leqslant C \varepsilon_\omega (1+z^2)$ and $| \partial_\omega X_j | \leqslant \frac{C \varepsilon_\omega}{\omega} (1+z^2)$, we find that
\[ | \mathbf{T}_{1A} | \leqslant \frac{C \varepsilon_\omega \varrho_\omega}{\omega \alpha} (1 + |y|) e^{- \alpha |y|} , \, \, \, \, | \mathbf{T}_{1B}| \leqslant \frac{C \varepsilon_\omega^2}{\omega \alpha} e^{- \alpha |y|} \, \, \, \, \text{and} \, \, \, \, | \mathbf{T}_{1C} | \leqslant \frac{C \varrho_\omega}{\omega} e^{- \alpha |y|}. \]
This leads to:
\[ | \partial_\omega \mathbf{T}_1 | \leqslant \frac{C \varepsilon_\omega \varrho_\omega}{\omega \alpha} (1 + |y|) e^{- \alpha |y|}. \]
Finally, we write $\partial_\omega \mathbf{T}_2 = \mathbf{T}_{2A} + \mathbf{T}_{2B}$, where
\[ \mathbf{T}_{2A} = \frac{\alpha ' ( \omega )}{2} \int_{\R} \partial_\alpha \left ( \frac{e^{- \alpha |y-z|} - e^{- \alpha |y|}}{\alpha} \right ) b_\omega^+ (z) \, \text{d}z \, \, \, \, \, \, \, \, \text{and} \, \, \, \, \, \, \, \, \mathbf{T}_{2B} = \int_{\R} \frac{e^{- \alpha |y-z|} - e^{- \alpha |y|}}{2 \alpha} \, \partial_\omega b_\omega^+ (z) \, \text{d}z. \]
Recalling that $| e^{- \alpha (|y-z| - |y|)} -1 | \leqslant \alpha |z| e^{\alpha |z|}$, we see that
\[ | \mathbf{T}_{2B} | \leqslant \frac{C \varepsilon_\omega}{\omega} \, e^{- \alpha |y|}. \]
For the term $\mathbf{T}_{2A}$, we see that
\[ \partial_\alpha  \left ( \frac{e^{- \alpha |y-z|} - e^{- \alpha |y|}}{\alpha} \right ) = \frac{1}{\alpha^2} \left [ e^{- \alpha |y|} - e^{- \alpha |y-z|} + \alpha \left ( |y| e^{- \alpha |y|} - |y-z| e^{- \alpha |y-z|} \right ) \right ] \]
where
\[ \begin{array}{rl} & |e^{- \alpha |y|} - e^{- \alpha |y-z|}| \leqslant \alpha (1 + |y|) e^{- \alpha |y|} |z| e^{\alpha |z|} \\
\\ \text{and} & \left | \, |y| e^{- \alpha |y|} - |y-z| e^{- \alpha |y-z|} \right | \leqslant |y| \, |e^{- \alpha |y|} - e^{- \alpha |y-z|} | + | \, |y| - |y-z| \, | e^{- \alpha |y-z|} \leqslant (1 + \alpha |y|) e^{- \alpha |y|} |z| e^{\alpha |z|}. \end{array} \]
Gathering these estimates, we evidently find that
\[ | \mathbf{T}_{2A}| \leqslant \frac{C \varrho_\omega}{\omega} (1 + |y|) e^{- \alpha |y|}. \]
Putting all the pieces together, we obtain:
\[ | \partial_\omega X_1 | \leqslant \frac{C \varepsilon_\omega \varrho_\omega}{\omega \alpha} (1 + |y|) e^{- \alpha |y|}. \]
Estimating $\partial_y^k \partial_\omega X_1$ requires the same proof, with minor adjustments. We differentiate the expression $X_1 = - e_1 \cdot H_\alpha^{-1} ( P_\omega X)$ with regards to $y$. For example, we shall write:
\[ \begin{array}{rcl} - \partial_y X_1 &=& \displaystyle{- \frac{1}{2} \int_{\R} \text{sgn} (y-z) e^{- \alpha |y-z|} (b_\omega^+ (X_1-1) + b_\omega^- X_2)(z) \, \text{d}z} \\
\\ & & \, \, \, \, \, \, \, \displaystyle{- \frac{1}{2} \int_{\R} (e^{- \alpha |y-z|} - e^{- \alpha |y|} ) \text{sgn} (y-z) b_\omega^+ (z) \, \text{d}z} \\
\\ & & \, \, \, \, \, \, \, \displaystyle{- \frac{e^{- \alpha |y|}}{2} \int_{\R} \text{sgn} (y-z) b_\omega^+ (z) \, \text{d}z.} \end{array} \]
Controlling as previously, we find that 
\[  | \partial_y \partial_\omega X_1 | \leqslant \frac{C \varepsilon_\omega}{\omega} (1 + |y|) e^{- \alpha |y|}. \]
More generally, a similar proof would show that, for any $k \geqslant 1$,
\[ | \partial_y^k \partial_\omega X_1 | \leqslant \frac{C \varepsilon_\omega^k}{\omega} (1 + |y|) e^{- \alpha |y|}. \]
As for $X_2$, the general idea is the same but the calculations are easier and we see that, for any $k \geqslant 0$,
\[ | \partial_y^k \partial_\omega X_2 | \leqslant \frac{C \varepsilon_\omega}{\omega} (1 + |y|) e^{- \kappa |y|}. \]
It follows that, for any $j \in \{ 1 \, , 2 \}$ and any $k \geqslant 1$,
\[ \begin{array}{rl} & | \partial_\omega W_j | \leqslant \frac{C \varepsilon_\omega \varrho_\omega}{\omega I_\omega} (1 + |y|) e^{- \alpha |y|} + \frac{C \varepsilon_\omega}{\omega} (1 + |y|) e^{- \kappa |y|} \leqslant \frac{C \varepsilon_\omega \varrho_\omega}{\omega I_\omega} (1 + |y|) e^{- \alpha |y|} \\
\\ \text{and} & | \partial_y^k \partial_\omega W_j | \leqslant \frac{C \varepsilon_\omega^k}{\omega} (1 + |y|) e^{- \alpha |y|} + \frac{C \varepsilon_\omega}{\omega} (1 + |y|) e^{- \kappa |y|} \leqslant \frac{C \varepsilon_\omega}{\omega} (1 + |y|) e^{- \alpha |y|}. \end{array} \]
One can notice that $\frac{\varepsilon_\omega}{\omega} \leqslant \frac{C \varepsilon_\omega \varrho_\omega}{\omega I_\omega}$. Now, differentiating the expressions $V_1 = (S^*)^2W_1$ and $V_2 = \lambda^{-1} L_+ V_1$ with regards to $\omega$ and using previous estimates (including $| \partial_\omega Q_\omega | \leqslant \frac{C}{\omega} (1 + |y|) e^{-|y|}$ and $| \alpha '( \omega ) | \leqslant \frac{C \varepsilon_\omega}{\omega}$), we evidently find that,
\[ | \partial_\omega V_1 | + | \partial_\omega V_2 | + | \partial_y \partial_\omega V_1 | + | \partial_y \partial_\omega V_2 | \leqslant \frac{C}{\omega} \left ( \frac{\varepsilon_\omega \varrho_\omega}{I_\omega} +1 \right ) (1 + |y|) e^{- \alpha |y|}, \]
which concludes the proof. Note that we are not able to compare $\frac{\varepsilon_\omega \varrho_\omega}{I_\omega}$ to $1$. \hfill \qedsymbol

\noindent \\ At this point the existence of the internal mode is established. Its uniqueness remains to be proven. To do so we procede to a new factorization, which will also be useful for the proof of the asymptotic stability. Set $\displaystyle{U = \partial_y - \frac{W_2 '}{W_2}}$. 

\begin{leftbar}
\noindent \textbf{Lemma 2.} For $\omega > 0$ small enough, the following factorization holds:
\[ U M_+ M_- = K U, \]
where $K = \partial_y^4 - 2 \partial_y^2 + K_2 \partial_y^2 + K_1 \partial_y + K_0 + 1$, with
\[ \begin{array}{rcl} K_2 &=& \displaystyle{1 - \lambda \frac{W_1}{W_2} + 3 \frac{W_2''}{W_2} - 4 \frac{(W_2')^2}{W_2^2} - a_\omega^+,} \\
\\ K_1 &=& \displaystyle{-3 \lambda \frac{W_1'}{W_2} + 3 \lambda \frac{W_1 W_2'}{W_2^2} + 3 \frac{W_2'''}{W_2} - 11 \frac{W_2'W_2''}{W_2^2} + 8 \frac{(W_2')^3}{W_2^3} - (a_\omega^+)',} \\
\\ K_0 &=& \displaystyle{ -2 \lambda \frac{W_1''}{W_2} + 5 \lambda \frac{W_1' W_2'}{W_2^2} + 2 \frac{(W_2')^2}{W_2^2} - 3 \lambda \frac{W_1 (W_2')^2}{W_2^3} + \lambda \frac{W_1 W_2''}{W_2^3} - \frac{W_2''}{W_2} + \frac{W_2''''}{W_2} -5 \frac{W_2' W_2'''}{W_2^2} } \\
\\ & & \displaystyle{ \, \, \, -3 \frac{(W_2'')^2}{W_2^2} + 15 \frac{W_2'' (W_2')^2}{W_2^3} -8 \frac{(W_2')^4}{W_2^4} - (a_\omega^+)' \frac{W_2'}{W_2} - a_\omega^+ \frac{W_2''}{W_2} + 2 a_\omega^+ \frac{(W_2')^2}{W_2^2} + \lambda^2 -1.} \end{array} \]
\end{leftbar}

\noindent \textit{Proof.} See the proof of Lemma 3 in \cite{Ma0}. The only difference is the start: in our case,
\[ \begin{array}{rcl} M_+ M_ h - \lambda^2 h &=& -2 \lambda W_1 ' k' - \lambda W_1 k'' + 2 W_2''' k' + 4 W_2''k'' + 2 W_2' k''' + W_2'' k'' \\ & & \, \, \, \, \, \, + 2 W_2' k''' + W_2 k'''' -2 W_2'k' - W_2k'' -2 a_\omega^+ W_2' k' - a_\omega^+ W_2 k'', \end{array} \]
where $h$ is any smooth function and $k = h/W_2$. Note that the potential $a_\omega^-$ does not appear, since the only time we actually use the operator $M_-$ is in the equality $M_- W_2 = W_1$. Starting from the relation above, the rest of the proof is entirely identical. \hfill \qedsymbol

\noindent \\ We need estimates on the functions $K_j$ involved in the operator $K$.

\begin{leftbar}
\noindent \textbf{Lemma 3.} For $\omega > 0$ small enough, for any $j \in \{ 0 \, , 1 \, , 2 \}$ and $k \in \{ 0 \, , j+1 \}$, on $\R$,
\[ |K_j^{(k)} | \leqslant C \varepsilon_\omega e^{- ( \kappa - \alpha ) |y|}. \]
\end{leftbar}

\noindent \textit{Proof.} We adapt the proof of Lemma 3 in \cite{Ma0}. First, taking $\omega$ small enough, since $|W_1 (y) - e^{- \alpha |y|}| \leqslant C \varrho_\omega e^{- \alpha |y|}$, we have $W_2 (y) \geqslant \frac{1}{2} e^{- \alpha |y|} > 0$. Exploiting Proposition 2, we see that
\[ \left | \frac{W_1}{W_2} -1 \right | + \left | \left ( \frac{W_1}{W_2} \right )^{(k)} \right | \leqslant C_k \varepsilon_\omega e^{-(\kappa - \alpha) |y|} \]
for any $k \geqslant 1$. Now, we rewrite the identity $M_- W_2 = \lambda W_1$ as $W_2 '' = \alpha^2 W_2 - w_0 W_2$ where $w_0 = \lambda \, \frac{W_1 - W_2}{W_2} - a_\omega^-$. Using the estimates on $W_1 - W_2$ and on $W_2$, as well as $| \lambda | \leqslant 1$ and $|a_\omega^-| \leqslant C \varepsilon_\omega e^{-2|y|} \leqslant C \varepsilon_\omega e^{- ( \kappa - \alpha ) |y|}$, we see that $|w_0^{(k)}| \leqslant C_k \varepsilon_\omega e^{-( \kappa - \alpha ) |y|}$ for any $k \in \{ 0 \, , ... \, , 3 \}$. This means that
\[ \left | \frac{W_2''}{W_2} - \alpha^2 \right | + \left | \left ( \frac{W_2 ''}{W_2} \right )^{(k)} \right | \leqslant C_k \varepsilon_\omega e^{-(\kappa - \alpha) |y|} \]
for any $k \in \{ 1 \, , ... \, , 3 \}$. Now, take $y \geqslant 0$. Multiplying the identity $W_2 '' = \alpha^2 W_2 - w_0 W_2$ by $W_2'$ and integrating on $[y \, , + \infty )$, we get
\[ (W_2')^2 (y) = \alpha^2 W_2^2(y) + 2 \int_y^{+ \infty} w_0 W_2' W_2. \]
Using the estimates on $w_0$, $W_2'$ and $W_2$, we find that
\[ \left | \frac{(W_2')^2}{W_2^2} - \alpha^2 \right | \leqslant C \varepsilon_\omega I_\omega e^{- ( \kappa - \alpha ) y} + C \varepsilon_\omega^2 e^{-2( \kappa - \alpha ) y}.  \]
For $y > \frac{1}{\kappa - \alpha} \ln \left ( \frac{C \varepsilon_\omega}{I_\omega} \right ) =: y_\omega$ with an appropriate constant $C$, we have both $C \varepsilon_\omega I_\omega e^{- ( \kappa - \alpha ) y} < \alpha^2$ and $\varepsilon_\omega^2 e^{-2( \kappa - \alpha ) y} \leqslant C \varepsilon_\omega I_\omega e^{- ( \kappa - \alpha ) y}$. Therefore, for $y > y_\omega$, 
\[ \left | \frac{(W_2')^2}{W_2^2} - \alpha^2 \right | \leqslant C \varepsilon_\omega I_\omega e^{- ( \kappa - \alpha ) y} < \alpha^2, \]
thus $W_2'(y) \neq 0$. Since $W_2 '(y) \sim - \alpha e^{- \alpha y}$ when $y \to + \infty$, we see that $W_2 ' < 0$ for $y > y_\omega$. For such $y$, $\left | \frac{W_2'}{W_2} - \alpha \right | = - \frac{W_2'}{W_2} + \alpha \geqslant \alpha > 0$ and then
\[ \left | \frac{W_2'}{W_2} + \alpha \right | = \frac{\left | \frac{(W_2')^2}{W_2^2} - \alpha^2 \right |}{\left | \frac{W_2'}{W_2} - \alpha \right |} \leqslant \frac{C \varepsilon_\omega I_\omega e^{- ( \kappa - \alpha ) y}}{\alpha} \leqslant C \varepsilon_\omega e^{-(\kappa - \alpha) |y|}. \]
Now, for $0 \leqslant y \leqslant y_\omega$, recall from Proposition 2 that $|W_2 '(y) + \alpha e^{- \alpha y}| \leqslant C \varrho_\omega I_\omega e^{- \alpha y} + C \varepsilon_\omega e^{- \kappa y}$ and $|W_2 (y) - e^{- \alpha y} | \leqslant C \varrho_\omega e^{-\alpha y}$. Also recalling that $W_2 \geqslant \frac{1}{2} e^{- \alpha y}$, this leads to
\[ \left | \frac{W_2'}{W_2} + \alpha \right | = \frac{1}{W_2} \left | W_2'(y) + \alpha e^{- \alpha y} - \alpha (e^{- \alpha y} - W_2 (y) ) \right | \leqslant C \varrho_\omega I_\omega + C \varepsilon_\omega e^{- (\kappa - \alpha )y}. \]
For $y \leqslant y_\omega$ we easily see that $\varrho_\omega I_\omega \leqslant C \varepsilon_\omega e^{- ( \kappa - \alpha ) y}$. Indeed, $\varepsilon_\omega e^{- ( \kappa - \alpha ) y} \geqslant \varepsilon_\omega C I_\omega / \varepsilon_\omega = C I_\omega \geqslant C I_\omega \varrho_\omega$. This proves that, for all $y \geqslant 0$,
\[ \left | \frac{W_2'}{W_2} + \alpha \right | \leqslant C \varepsilon_\omega e^{- ( \kappa - \alpha ) y}. \]
Then, using the relation $W_2^{(k+2)} = \alpha^2 W_2^{(k)} - (w_0 W_2)^{(k)}$, we deduce that
\[ \left | \frac{W_2^{(k)}}{W_2} - (- \alpha)^k \right | \leqslant C \varepsilon_\omega e^{-( \kappa - \alpha ) y} \]
for all $k \in \{ 1 \, , ... \, , 5 \}$ and all $y \geqslant 0$. For $y \leqslant 0$, the result must be adapted by taking into account the parity of $k$, since $W_2$ is an even function. By similar considerations we can show that
\[ \left | \frac{W_1^{(k)}}{W_2} - (- \alpha)^k \right | \leqslant C \varepsilon_\omega e^{-( \kappa - \alpha ) y} \]
for all $k \in \{ 1 \, , ... \, , 5 \}$ and all $y \geqslant 0$. Now, we can establish the estimates on $K_0$, $K_1$ and $K_2$. First, for $y \geqslant 0$,
\[ \begin{array}{rcl} |K_2| & = & |K_2 - (1- \lambda + 3 \alpha^2 - 4 \alpha^2) | \\
\\ & \leqslant & \displaystyle{|1-1| + \left | - \lambda \left ( \frac{W_1}{W_2} -1 \right ) \right | + 3 \left | \frac{W_2''}{W_2} - \alpha^2 \right | + 4 \left | \frac{(W_2')^2}{W_2^2} - \alpha^2 \right | + |a_\omega^+|} \\
\\ & \leqslant & C \varepsilon_\omega e^{-( \kappa - \alpha)y}. \end{array} \]
the proofs are identical for $K_1$ and $K_0$: we respectively use the identities $3 \lambda \alpha - 3 \lambda \alpha - 3 \alpha^3 + 11 \alpha^3 - 8 \alpha^3 = 0$ and $-2 \lambda \alpha^2 + 5 \lambda \alpha^2 + 2 \alpha^2 - 3 \lambda \alpha^2 + \lambda \alpha^2 - \alpha^2 + \alpha^4 -5 \alpha^4 - 3 \alpha^4 + 14 \alpha^4 -8 \alpha^4 + \lambda^2 -1=0$. The result for $y \leqslant 0$ holds by parity, and the generalisation for $k \geqslant 1$ does not present additional difficulties. \hfill \qedsymbol

\noindent \\ The virial computation $\int_{\R} (2yh'+h)Kh$ expands as follows.

\begin{leftbar}
\noindent \textbf{Lemma 4.} It holds formally
\[  \int_{\R} (2yh' + h)Kh = 4 \int_{\R} (h'')^2 + 4 \int_{\R} (h')^2 + \int_{\R} Y_1(h')^2 + \int_{\R} Y_0 h^2 \]
where the functions $Y_1 = -2K_2 -yK_2' +2yK_1$ and $Y_0 = \frac{1}{2} \left ( K_2'' - K_1' - 2yK_0' \right )$ satisfy, $|Y_1^{(k)}| \leqslant C_k \varepsilon_\omega e^{-|y|}$ for all $k \in \{ 0 \, , ... \, , 2 \}$ and $|Y_0| \leqslant C \varepsilon_\omega e^{-|y|}$ on $\R$. 
\end{leftbar}

\noindent \textit{Proof.} See the proof of Lemma 4 in \cite{Ma0}. \hfill \qedsymbol

\noindent \\ As in \cite{Si}, the coercivity of the virial functional relies on the sign of the integral $\int_{\R} Y_0$. The following lemma grants us this information.

\begin{leftbar}
\noindent \textbf{Lemma 5.} For $\omega > 0$ small, $\int_{\R} Y_0 = I_\omega \left ( 1 + \mathcal{O} (\varrho_\omega) \right )$.
\end{leftbar}

\noindent \textit{Proof.} We begin by writing $K_0 = -2 \left ( \frac{W_1''}{W_2} - \alpha^2 \right ) + \left ( \frac{W_2''''}{W_2} - \alpha^4 \right ) + \widetilde{K}_0$ where
\[ \begin{array}{rcl} \widetilde{K}_0 &=& \displaystyle{2 \alpha^2 \frac{W_1''}{W_2} + 5 \lambda \frac{W_1' W_2'}{W_2^2} + 2 \frac{(W_2')^2}{W_2^2} -3 \lambda  \frac{W_1 (W_2')^2}{W_2^3} + \lambda \left ( \frac{W_1}{W_2} -1 \right ) \frac{W_2''}{W_2}} \\
\\ & & \, \, \, \displaystyle{- \alpha^2 \frac{W_2''}{W_2} -5 \frac{W_2'}{W_2} \frac{W_2'''}{W_2} -3 \frac{(W_2'')^2}{W_2^2} + 15 \frac{W_2''}{W_2} \frac{(W_2')^2}{W_2^2} -8 \frac{(W_2')^4}{W_2^4}} \\
\\ & & \, \, \, \displaystyle{- (a_\omega^+)' \frac{W_2'}{W_2} - a_\omega^+ \frac{W_2''}{W_2} + 2 a_\omega^+ \frac{(W_2')^2}{W_2^2} + (\lambda^2 - 1) -2 \alpha^2 + \alpha^4.} \end{array} \]
We begin to control the terms in $\widetilde{K}_0$ as we did in the proof of Lemma 3. For example,
\[ \left | 2 \alpha^2 \frac{W_1''}{W_2} - 2 \alpha^4 \right | = 2 \alpha^2 \left | \frac{W_1''}{W_2} - \alpha^2 \right | \leqslant C \alpha^2 \varepsilon_\omega e^{-( \kappa - \alpha ) y} \]
and
\[ \left | 5 \lambda \frac{W_1' W_2'}{W_2^2} -5 \lambda \alpha^2 \right | \leqslant 5 \lambda \left (  \left | \frac{W_1'}{W_2} \left ( \frac{W_2'}{W_2} + \alpha \right ) \right | + \left | - \alpha \left ( \frac{W_1'}{W_2} + \alpha \right ) \right | \right ) \leqslant C \varepsilon_\omega^2 e^{-( \kappa - \alpha ) |y|}. \]
Proceeding similarly for the other terms, we obtain 
\[ | \widetilde{K}_0 | \leqslant C \varepsilon_\omega^2 e^{- ( \kappa - \alpha ) |y|}. \]
Now, we analyse the remaining terms. First,
\[ \int_{\R} \left ( \frac{W_1''}{W_2} - \alpha^2 \right ) \, \text{d}y = \int_{\R} \left ( \frac{W_1'}{W_2} \right )' \, \text{d}y +  \int_{\R} \left ( \frac{W_1' W_2'}{W_2^2} - \alpha^2 \right ). \]
Since $\left | \frac{W_1'}{W_2} + \alpha \right | \leqslant C \varepsilon_\omega e^{-( \kappa - \alpha )y} \longrightarrow 0$ as $y \to + \infty$, we have $\frac{W_1'}{W_2} \longrightarrow - \alpha$ as $y \to + \infty$. Hence, recalling that $W_1'/W_2$ is odd,
\[ \int_{\R} \left ( \frac{W_1'}{W_2} \right )' \, \text{d}y = \left [ \frac{W_1'}{W_2} \right ]_{- \infty}^{+ \infty} = -2 \alpha. \]
Besides, $\left | \frac{W_1' W_2'}{W_2^2} - \alpha^2 \right | \leqslant C \varepsilon_\omega^2 e^{- ( \kappa - \alpha ) |y|}$. Now, the last term can be written as follows:
\[ \int_{\R} \left ( \frac{W_2''''}{W_2} - \alpha^4 \right ) \, \text{d}y = 2 \int_0^{+ \infty} \left ( \left ( \frac{W_2'''}{W_2} \right ) ' + \frac{W_2'}{W_2} \left ( \frac{W_2'''}{W_2} + \alpha^3 \right ) - \alpha^3 \left ( \frac{W_2'}{W_2} + \alpha \right ) \right ). \]
We have $\int_0^{+ \infty} \left ( \frac{W_2'''}{W_2} \right ) ' = - \alpha^3$, $\left | \frac{W_2'}{W_2} \left ( \frac{W_2'''}{W_2} + \alpha^3 \right ) \right | \leqslant C \varepsilon_\omega^2 e^{-(\kappa - \alpha) y}$ and $\left | \frac{W_2'}{W_2} + \alpha \right | \leqslant C \varepsilon_\omega e^{-( \kappa - \alpha )y}$. Gathering all these estimates, we find that
\[ \int_{\R} Y_0 = \int_{\R} K_0 = 4 \alpha + \mathcal{O} ( \varepsilon_\omega^2 ) = I_\omega ( 1 + \mathcal{O} ( \varrho_\omega )) + \mathcal{O} ( I_\omega \varrho_\omega) = I_\omega \left ( 1 + \mathcal{O} ( \varrho_\omega ) \right ) \]
which is the desired result. \hfill \qedsymbol

\noindent \\ We can now prove a coercivity property on the operator $K$ which shows that the second transformed problem ($KZ = \mu Z$) has no solution.

\begin{leftbar}
\noindent \textbf{Lemma 6.} Assume hypotheses $(H_1)$ and $(H_2)$ hold. For $\omega > 0$ small, if $(\mu \, , Z) \in \R \times H^4 ( \R )$ satisfied $KZ = \mu Z$ then $Z=0$.
\end{leftbar}

\noindent \textit{Proof.} Let $( \mu \, , Z ) \in \R \times H^4 ( \R )$ be a solution of $KZ = \mu Z$. Since $\int_{\R} (2y Z' + Z)Z=0$, we deduce from Lemma 4 that
\[ 0 = 4 \int_{\R} (Z'')^2 + 4 \int_{\R} (Z')^2 + \int_{\R} Y_1 (Z')^2 + \int_{\R} Y_0 Z^2. \]
First, from Lemma 4 we know that $|Y_1| \leqslant C \varepsilon_\omega$ thus $\left | \int_{\R} Y_1 (Z')^2 \right | \leqslant C \varepsilon_\omega \int_{\R} (Z')^2$. Now we use Lemma 5 from \cite{Ma0} with $Y = \frac{Y_0}{C \varepsilon_\omega}$, $c=1$ and $h = Z$. It is correct since, for $\omega > 0$ small, $\int_{\R} \frac{Y_0}{C \varepsilon_\omega} \sim \frac{I_\omega}{C \varepsilon_\omega} > 0$. This lemma gives us
\[ 0 \leqslant \left ( \int_{\R} Y_0 \right ) \int_{\R} e^{-|y|} Z^2 \, \text{d}y \leqslant C \int_{\R} Y_0 Z^2 + \frac{C \varepsilon_\omega^2}{\int_{\R} Y_0} \int_{\R} (Z')^2 \leqslant C \int_{\R} Y_0 Z^2 + \frac{C \varepsilon_\omega^2}{I_\omega} \int_{\R} (Z')^2. \]
Hence, 
\[ - \int_{\R} Y_0 Z^2 \leqslant \frac{C \varepsilon_\omega^2}{I_\omega} \int_{\R} (Z')^2 = C \varrho_\omega \int_{\R} (Z')^2. \]
Putting these estimates together, we obtain
\[ 0 = 4 \int_{\R} (Z'')^2 + 4 \int_{\R} (Z')^2 + \int_{\R} Y_1 (Z')^2 + \int_{\R} Y_0 Z^2 \geqslant 4 \int_{\R} (Z'')^2 + (4 - C \varepsilon_\omega - C \varrho_\omega ) \int_{\R} (Z')^2. \]
Now, taking $\omega >0$ small enough so that $C \varepsilon_\omega + C \varrho_\omega < 1$, we have
\[ 0 \geqslant \int_{\R} (Z'')^2 + \int_{\R} (Z')^2 \]
which leads to $Z=0$. \hfill \qedsymbol

\noindent \\ Going back to the original problem, we establish the uniqueness of the internal mode.

\begin{leftbar}
\noindent \textbf{Lemma 7.} Assume hypotheses $(H_1)$ and $(H_2)$ hold. For $\omega > 0$ small, the only solutions $(\widetilde{\lambda} \, , \widetilde{V}_1 \, , \widetilde{V}_2 ) \in [0 \, , + \infty ) \times H^2 ( \R ) \times H^2 ( \R )$ of the eigenvalue problem \eqref{sysL} are
\begin{itemize}
    \item $(\mu \, , 0 \, , 0)$ for any $\mu \geqslant 0$,
    \item $(0 \, , a Q_\omega ' \, , b Q_\omega )$ for any $a,b \in \R$,
    \item $( \lambda \, , c V_1 \, , c V_2 )$ for any $c \in \R$, where $(\lambda \, , V_1 \, , V_2)$ is the internal mode constructed in Proposition 2.
\end{itemize}
\end{leftbar}

\noindent \textit{Proof.} See the proof of Lemma 8 in \cite{Ma0}. \hfill \qedsymbol

\noindent \\ Gathering Proposition 2 and Lemma 7, we obtain Theorem 1 (in its rescaled version). \hfill \qedsymbol

\section{Rescaled decomposition}
\noindent The two main points of this paper, in order to establish the asymptotic stability property, are the understanding of the internal mode (existence, uniqueness, properties, estimates) and the Fermi golden rule. The rest of this paper relies on an adaptation of the proof in \cite{Ma0} (which concerns the case $g(s)=s^2$). We point out the following three technical differences.
\begin{itemize}
	\item In \cite{Ma0}, $\alpha ( \omega ) \sim \frac{8}{9} \omega$ and $\varepsilon_\omega = C \omega$, which leads to more comfortable calculations. Thus, in most proofs, some occurrences of $\omega_0$ must be replaced by $\alpha ( \omega_0 )$ in our case. Other occurrences of $\omega_0$ must be replaced by $\varepsilon_{3 \omega_0 /2}$. Statements will be adapted in consequence. The definition of the weight function $\rho$ will also be adapted (see definition just before Lemma 10).
	\item The estimates on $\partial_y^k \partial_\omega V_j$ have simpler form in \cite{Ma0} and lead to easier calculations. In our case, recalling that $\frac{\omega_0}{2} \leqslant \omega \leqslant \frac{3 \omega_0}{2}$, we know that, for $j \in \{ 1 \, , 2 \}$,
\begin{equation}
    | \omega \partial_\omega V_j  | + | \omega \partial_y \partial_\omega V_j | \leqslant C \underline{\mathbf{V}} ( \omega_0 ) (1 + |y|) e^{- \alpha |y|}
    \label{dVomega}
\end{equation}
	where $\underline{\mathbf{V}} ( \omega_0 ) := \frac{\varepsilon_{3 \omega_0/2}}{\alpha ( \omega_0 )}$. 
	\item In Lemma 9 below we introduce $q_1$ and $q_2$ which are the nonlinear terms that appear when we linearize the evolution equation around the solitary wave. Denoting $\text{Taylor}_k (q_j)$ the $k^\text{th}$-order Taylor expansion of $q_j$, Lemma 10 in \cite{Ma0} shows that $|q_1 - \text{Taylor}_2 (q_1)| + |q_2 - \text{Taylor}_2 (q_2)| \leqslant C |u|^3$. In the general case, typically if $g(s) = s^\sigma$ with $\sigma \in (1 \, , 2)$, the bound $|q_2 - \text{Taylor}_2 (q_2)| \leqslant C |u|^3$ remains true but the bound $|q_1 - \text{Taylor}_2 (q_1)| \leqslant C |u|^3$ is not valid anymore. We can only obtain $|q_1 - \text{Taylor}_2 (q_1)| \leqslant C |u|^{7/3}$. The fact that the exponent is not $3$ does not matter; what matters is that the exponent is greater than $2$. Most proofs find themselves changed ($\epsilon^{1/3}$ instead of $\epsilon$ for instance) without further complication.
\end{itemize}
\noindent We introduce $\Lambda := \frac{1}{2} (1 + y \partial_y)$, $\Lambda^* = - \frac{y}{2} \partial_y$ and $\Lambda_\omega := \Lambda + \omega \partial_\omega$, such that $\omega \partial_\omega \phi_\omega (x) = \sqrt{\omega} \, \Lambda_\omega Q_\omega ( \sqrt{\omega} \, x)$. Set $\R_+^2 := \R \times (0 \, , + \infty )$. For $\varphi \in H^1 ( \R )$ and $\Pi = ( \gamma \, , \omega ) \in \R_+^2$, we define the function $\zeta [ \varphi \, , \Pi ] \, : \, \R \to \C$ by
\[ \zeta [ \varphi \, , \Pi ] (y) = \frac{e^{-i \gamma}}{\sqrt{\omega}} \varphi \left ( \frac{y}{\sqrt{\omega}} \right ). \]

\begin{leftbar}
\noindent \textbf{Lemma 8.} For any $\omega_0 > 0$ small and any $\epsilon >0$, there exists $\delta >0$ such that, for all even function $\varphi \in H^1 ( \R )$ with $|| \varphi - \phi_{\omega_0}||_{H^1 ( \R )} < \delta$, there exists a unique $\Pi = (\gamma \, , \omega ) \in \R_+^2$ such that $| \gamma | + | \omega - \omega_0 | < \epsilon$ and $u := \zeta [ \varphi \, , \Pi ] - Q_\omega$ satisfies
\[ ||u||_{H^1 ( \R)} < \epsilon \, \, \, \, \, \, \, \text{and} \, \, \, \, \, \, \, \langle u \, , i \Lambda_\omega Q_\omega \rangle = \langle u \, , Q_\omega \rangle = 0. \]
\end{leftbar}

\noindent \textit{Proof.} See the proof of Lemma 9 in \cite{Ma0}. We need to know that $\frac{\sqrt{\omega_0}}{2} \partial_\omega \left ( || \phi_\omega ||^2 \right )_{\omega = \omega_0}$ is positive for $\omega_0 >0$ small enough. This is proven in Lemma 5 in \cite{Ri}. \hfill \qedsymbol

\noindent \\ We now prove a technical lemma. 

\begin{leftbar}
\noindent \textbf{Lemma 9.} We set $f_\omega ( \psi ) := | \psi |^2 \psi + \frac{g ( \omega | \psi |^2 )}{\omega} \, \psi$. Let
\[ \begin{array}{rl} & q_1 = \text{Re} \left [ f_\omega (Q_\omega + u) - f_\omega (Q_\omega) - f_\omega '(Q_\omega) u \right ] \\ \\ \text{and} & q_2 = \text{Im} \left [ f_\omega (Q_\omega + u) - \frac{f_\omega (Q_\omega)}{Q_\omega} \, u \right ]. \end{array} \]
We have, for $u = u_1 + iu_2$ with $|u| < 1$,
\[ \begin{array}{rl} & \left | q_1 - \left [ Q_\omega (3 + 3 g'(\omega Q_\omega^2) + 2 \omega Q_\omega^2 g''(\omega Q_\omega^2))u_1^2 + Q_\omega (1 + g'( \omega Q_\omega^2)) u_2^2 \right ] \right | \leqslant C |u|^{7/3} \\
\\ \text{and} & \left | q_2 - 2 Q_\omega (1 + g'( \omega Q_\omega^2) ) u_1 u_2 \right | \leqslant C |u|^3. \end{array}  \]
\end{leftbar}

\noindent \textit{Proof.} Let us begin with $q_1$. First, consider the case $|u| \leqslant \frac{1}{100} Q_\omega^{3/2}$. We use Taylor's expansion and write that
\begin{equation}
    \begin{array}{rl} & q_1 - \left [ Q_\omega (3 + 3 g'(\omega Q_\omega^2) + 2 \omega Q_\omega^2 g''(\omega Q_\omega^2))u_1^2 + Q_\omega (1 + g'( \omega Q_\omega^2)) u_2^2 \right ] \\
    \\ =& |u|^2u_1 + 2u_1 |u|^2 Q_\omega^2 \omega g''( \omega Q_\omega^2) + u_1 |u|^2 g'( \omega Q_\omega^2) + 2 u_1^3 Q_\omega^2 \omega g''(\omega Q_\omega^2) \\
    \\ & \, \, \, \, \, + |u|^4 \omega Q_\omega \frac{g''(\omega Q_\omega^2)}{2} + 2 u_1^2 |u|^2 Q_\omega \omega g''(\omega Q_\omega^2) + |u|^4 u_1 \omega \frac{g''(\omega Q_\omega^2)}{2} + (Q_\omega + u_1) \frac{\text{IR}}{\omega} \end{array}
    \label{DLq1}
\end{equation}
where
\[ \text{IR} := \int_{\omega Q_\omega^2}^{\omega (Q_\omega^2 + 2 u_1 Q_\omega + |u|^2)} \frac{(\omega (Q_\omega + 2 u_1 Q_\omega + |u|^2) - s)^2}{2} \, g'''(s) \, \text{d}s. \]
For $\omega>0$ small enough and any $s \in [\omega Q_\omega^2 \, , \omega (Q_\omega^2 + 2 u_1 Q_\omega + |u|^2)]$, we have $|g'''(s)| \leqslant \varepsilon_\omega s^{-2}$ and $|(\omega (Q_\omega + 2 u_1 Q_\omega + |u|^2) - s)^2| \leqslant C \omega^2 (|u| Q_\omega + |u|^2)$. Using these bounds, we find that
\[ \left | \frac{\text{IR}}{\omega} \right | \leqslant C \varepsilon_\omega ( |u| Q_\omega + |u|^2 )^3 \frac{1}{Q_\omega^4} \left ( 1 + \frac{2 u_1}{Q_\omega} + \frac{|u|^2}{Q_\omega^2} \right )^{-1}. \]
The hypothesis $|u| \leqslant \frac{1}{100} Q_\omega^{3/2}$ implies that $1 + \frac{2u_1}{Q_\omega} + \frac{|u|^2}{Q_\omega^2} \geqslant 1 - \frac{Q_\omega^{1/2}}{50} \geqslant \frac{1}{2}$. It also implies that $\frac{1}{Q_\omega} \leqslant \frac{C}{|u|^{2/3}}$. This leads to
\[ \left | \frac{\text{IR}}{\omega} \right | \leqslant C \varepsilon_\omega \left ( \frac{|u|^3}{Q_\omega} + \frac{|u|^6}{Q_\omega^4} \right ) \leqslant C \varepsilon_\omega |u|^{7/3}. \]
Using the hypotheses on $g$, it is easy to check that the other terms in \eqref{DLq1} are also smaller (in module) than $C |u|^{7/3}$ (they are even controlled by $C |u|^3$). Consequently, in this first case,
\[  \left | q_1 - \left [ Q_\omega (3 + 3 g'(\omega Q_\omega^2) + 2 \omega Q_\omega^2 g''(\omega Q_\omega^2))u_1^2 + Q_\omega (1 + g'( \omega Q_\omega^2)) u_2^2 \right ]  \right | \leqslant C |u|^{7/3}. \]
Let us now consider the second case where $|u| \geqslant \frac{1}{100} Q_\omega^{3/2}$. This case is easier, as we simply estimate every term by triangular inequality. Using the hypotheses on $g$ and the bound $Q_\omega \leqslant C |u|^{2/3}$, we see that
\[ \left | Q_\omega (3 + 3 g'(\omega Q_\omega^2) + 2 \omega Q_\omega^2 g''(\omega Q_\omega^2))u_1^2 + Q_\omega (1 + g'( \omega Q_\omega^2)) u_2^2 \right | \leqslant C |u|^2 Q_\omega \leqslant C |u|^{8/3}. \]
As for the control of $q_1$, we write that, by definition of $f_\omega$, 
\[ q_1 = |u|^2 Q_\omega + |u|^2 u_1 + 2 u_1^2 Q_\omega - 2u_1 Q_\omega^2 g'( \omega Q_\omega^2 ) + \frac{Q_\omega + u_1}{\omega} \left ( g(\omega Q_\omega^2 + 2 \omega u_1 Q_\omega + \omega |u|^2 ) - g( \omega Q_\omega^2) \right ) \]
where $|g(\omega Q_\omega^2 + 2 \omega u_1 Q_\omega + \omega |u|^2 ) - g( \omega Q_\omega^2)| \leqslant C \omega |2 u_1 Q_\omega + |u|^2 |$ since $|g'(s)| \leqslant C$. We find that 
\[ |q_1| \leqslant C \left ( |u|^2 Q_\omega + |u|^3 + |u| Q_\omega^2 \right ) \leqslant C |u|^{7/3}. \]
Therefore, whatever case we are in, the following bound holds:
\[ \left | q_1 - \left [ Q_\omega (3 + 3 g'(\omega Q_\omega^2) + 2 \omega Q_\omega^2 g''(\omega Q_\omega^2))u_1^2 + Q_\omega (1 + g'( \omega Q_\omega^2)) u_2^2 \right ] \right | \leqslant C |u|^{7/3}. \]
Now, let us deal with $q_2$, which requires to distinguish two cases. First, consider the case $|u| \leqslant \frac{1}{100} Q_\omega$. We use Taylor's expansion and write that
\[ \begin{array}{rcl} q_2 - 2 Q_\omega (1 + g'(\omega Q_\omega^2)) u_1 u_2 &=& |u|^2 u_2 (1 + g'(\omega Q_\omega^2)) + 2 u_1^2 u_2 \omega Q_\omega^2 g''( \omega Q_\omega^2) 2 u_1 u_2 |u|^2 \omega Q_\omega g''( \omega Q_\omega^2) \\ \\ & & \, \, \, \, \, \, + \frac{1}{2} |u|^4 u_2 \omega g''(\omega Q_\omega^2) + u_2 \frac{\text{IR}}{\omega}, \end{array} \]
where IR is the same integral as above. Here the estimates are better, since we have $u$ in front of $\text{IR}$, and not $(Q_\omega + u)$. We find, reasoning as previously, that $|q_2 - 2 Q_\omega (1 + g'(\omega Q_\omega^2))u_1u_2| \leqslant C |u|^3$. The case $|u| \geqslant \frac{1}{100} Q_\omega$ does not present additional difficulty, and we find that the above estimate still holds. \hfill \qedsymbol

\noindent \\ Without additional hypotheses on $g$, $7/3$ is the best exponent that we can get. For the rest of the paper, we introduce the functions 
\[ \nu (y) = \text{sech} \left ( \frac{y}{10} \right ) \, \, \, \, \, \, \, \, \, \text{and} \, \, \, \, \, \, \, \, \, \rho (y) = \text{sech} \left ( \frac{\alpha ( \omega_0 )}{10} \, y \right ). \]
We give the following global decomposition result. Here, for a function $u$ depending on $s$, we denote $\dot{u} := \partial_s u$.

\begin{leftbar}
\noindent \textbf{Lemma 10.} For any $\omega_0 > 0$ small and any $\epsilon > 0$, there exists $\delta > 0$ such that, for all even function $\psi_0 \in H^1 ( \R )$ with $|| \psi_0 - \phi_{\omega_0}||_{H^1 ( \R )} < \delta$, there exists a unique $\mathscr{C}^1$ function $\Pi \, : \, [0 \, , + \infty ) \mapsto (\gamma \, , \omega ) \in \R_+^2$ such that, if $\psi$ is the solution of \eqref{NLS}, denoting
\[ u(s) := \zeta [ \psi ( \tau (s)) \, , \Pi (s) ] - Q_{\omega (s)} \, \, \, \, \, \, \text{where} \, \, \tau (s) := \int_0^s \frac{\text{d}s'}{\omega (s')}, \]
then the following properties hold, for all $s \in [0 \, , + \infty )$,
\begin{itemize}
    \item \textit{(Stability.)} $| \omega - \omega_0 | + ||u||_{H^1 ( \R )} \leqslant \epsilon$. 
    \item \textit{(Orthogonality relations.)} $\langle u \, , i \Lambda_\omega Q_\omega \rangle = \langle u \, , Q_\omega \rangle = 0$.
    \item \textit{(Equation.)} $u = u_1 + iu_2$ satisfies
    \begin{equation}
        \left \{ \begin{array}{ccl} \dot{u}_1 &=& L_- u_2 + \mu_2 + p_2 - q_2 \\ \dot{u}_2 & = & - L_+ u_1 - \mu_1 - p_1 + q_1 \end{array} \right.
        \label{sysu}
    \end{equation}
    where $m_\gamma := \dot{\gamma} -1$, $m_\omega := \dot{\omega} / \omega$, $\mu_1 = m_\gamma Q_\omega$, $\mu_2 = -m_\omega \Lambda_\omega Q_\omega$, $p_1 = m_\gamma u_1 + m_\omega \Lambda u_2$ and $p_2 = m_\gamma u_2 - m_\omega \Lambda u_1$. 
    \item \textit{(Control of the parameters.)} $|m_\gamma| + |m_\omega| \leqslant C || \nu u ||^2$.
\end{itemize}
\end{leftbar}

\noindent \textit{Proof.} See the proof of Lemma 11 in \cite{Ma0}. \hfill \qedsymbol

\noindent \\ In what follows, we will need the following remark. We recall that 
\[ | \alpha ' ( \omega ) | \leqslant \frac{C \varepsilon_\omega}{\omega} \leqslant \frac{C \varepsilon_{3 \omega_0 /2}}{\omega_0} \]
thanks to the definition of $\varepsilon_\omega$ and the bounds $\frac{\omega_0}{2} \leqslant \omega \leqslant \frac{3 \omega_0}{2}$. Thus, using Lemma 10 just above,
\[ | \alpha ( \omega ) - \alpha ( \omega_0 ) | \leqslant \frac{C \varepsilon_{3 \omega_0 /2}}{\omega_0} | \omega - \omega_0 | \leqslant \frac{C \varepsilon_{3 \omega_0 /2}}{\omega_0} \, \epsilon \leqslant \frac{1}{10} \alpha ( \omega_0 ) \]
if we take $\epsilon >0$ small enough (depending on $\omega_0$). Thus we can put ourselves in the case where $\alpha ( \omega ) \leqslant C \alpha ( \omega_0 )$ and $\alpha ( \omega )^{-1} \leqslant C \alpha ( \omega_0 )^{-1}$, and that is what we will do from now on. Recall from Proposition 2 that $\langle V_1 \, , V_2 \rangle \sim \alpha^{-1} > 0$. We introduce the notation
\[ h^\T := h - \frac{\langle h \, , V_1 \rangle}{\langle V_1 \, , V_2 \rangle} V_2 \, \, \, \, \, \, \, \, \text{and} \, \, \, \, \, \, \, \, h^\perp := h - \frac{\langle h \, , V_2 \rangle}{\langle V_1 \, , V_2 \rangle} V_1. \]
We decompose $u$ as the contribution of the internal mode and a term $v$ orthogonal to the internal mode.

\begin{leftbar}
\noindent \textbf{Lemma 11.} Under the assumptions of Lemma 10, possibly taking a smaller $\delta$, there exists a unique $\mathscr{C}^1$ function $b=b_1 + ib_2 \, : \, [ 0 \, , + \infty ) \to \C$ such that $v = v_1 + iv_2$, defined by
\[ u_1 = v_1 + b_1 V_1 \, \, \, \, \, \, \, \text{and} \, \, \, \, \, \, \, u_2 = v_2 + b_2 V_2, \]
satisfies, for all $s \in [0 \, , + \infty )$, the five following properties. 
\begin{itemize}
    \item \textit{(Stability.)} $||v||_{H^1} + |b| \leqslant \epsilon$. 
    \item \textit{(Orthogonality relations.)} $\langle v \, , i \Lambda_\omega Q_\omega \rangle = \langle v \, , Q_\omega \rangle = \langle v \, , i V_1 \rangle = \langle v \, , V_2 \rangle = 0$.
    \item \textit{(Control of the parameters.)}
    \begin{equation}
        |m_\gamma| + |m_\omega| \leqslant C \left ( ||\nu v||^2 + |b|^2 \right ).
        \label{m}
    \end{equation}
    \item \textit{(Equation of $v$.)} Setting $r_1 := - m_\omega b_2 \omega \partial_\omega V_2$ and $r_2 := m_\omega b_1 \omega \partial_\omega V_1$,
    \begin{equation}
        \left \{ \begin{array}{ccl} \dot{v}_1 &=& L_- v_2 + \mu_2 + p_2^\perp - q_2^\perp - r_2^\perp \\ \dot{v}_2 &=& -L_+v_1 - \mu_1 - p_1^{\T} + q_1^{\T} + r_1^{\T}. \end{array} \right.
        \label{sysv}
    \end{equation}
    \item \textit{(Equation of $b$.)} Setting $B_j := \frac{\langle p_j - q_j - r_j \, , V_j \rangle}{\langle V_1 \, , V_2 \rangle}$ for $j \in \{ 1 \, , 2 \}$, 
    \begin{equation}
        \left \{ \begin{array}{ccl} \dot{b}_1 &=& \lambda b_2 + B_2 \\ \dot{b_2} &=& - \lambda b_1 - B_1 \end{array} \right.
        \label{sysb}
    \end{equation}
    and
    \begin{equation}
        |B_1| + |B_2| \leqslant C \alpha ( \omega_0 ) ( |b|^2 + || \rho^4 v ||^2 ).
        \label{inegB}
    \end{equation}
\end{itemize}
\end{leftbar}

\noindent \textit{Proof.} Follow the proof of Lemma 12 in \cite{Ma0}. We define $b_1 = \frac{\langle u_1 \, , V_2 \rangle}{\langle V_1 \, , V_2 \rangle}$ and $b_2 = \frac{\langle u_2 \, , V_1 \rangle}{\langle V_1 \, , V_2 \rangle}$. We have
\[ \left | q_1 - \left [ \left ( Q_\omega (3 + 3 g'(\omega Q_\omega^2) + 2 \omega Q_\omega^2 g''( \omega Q_\omega^2) \right ) u_1^2 + Q_\omega (1 + g'(\omega Q_\omega^2)) u_2^2 \right ] \right | \leqslant C |v|^{7/3} + C |b|^{7/3} + C \nu^2 |v|^2 + C \nu^2 |b| \, |v|. \]
Setting $\widetilde{d}_1 ( \omega ) := \frac{1}{\langle V_1 \, , V_2 \rangle} \int_{\R} Q_\omega ( 3 + 3 g'(\omega Q_\omega^2) + 2 \omega Q_\omega^2 g''( \omega Q_\omega^2 )) V_1^3$ and $\widetilde{d}_2 ( \omega ) := \frac{1}{\langle V_1 \, , V_2 \rangle} \int_{\R} Q_\omega ( 1 + g'(\omega Q_\omega^2)) V_1 V_2^2$, this leads to
\[ \left | \frac{\int_{\R} q_1 V_1}{\langle V_1 \, , V_2 \rangle} - \widetilde{d}_1 ( \omega ) b_1^2 - \widetilde{d}_2 ( \omega ) b_2^2 \right | \leqslant C \alpha ( \omega_0 ) (  || \rho^4 v ||^2 + |b| \, || \nu v || ) + C |b|^{7/3}. \]
Since $|b|^{7/3} \leqslant \epsilon^{1/3} |b|^2 \leqslant C \alpha ( \omega_0 ) |b|^2$ for $\epsilon >0$ small enough (depending on $\omega_0$), we ultimately find \eqref{inegB}. \hfill \qedsymbol

\noindent \\ Estimating $|| \rho h^{ \perp / \top } ||$ is not different from estimating $|| \rho h ||$, as the following lemma shows.

\begin{leftbar}
\noindent \textbf{Lemma 12.} For all $k \in \{ 0 \, , ... \, , 2 \}$, $|(h^\perp)^{(k)} | + |(h^{\T})^{(k)}| \leqslant C |h^{(k)}| + C \sqrt{\alpha (\omega_0 )} || \rho^4 h || \rho^8$. In particular, $||\rho h^\perp|| + ||\rho h^\T || \leqslant C || \rho h ||$. 
\end{leftbar}

\noindent \textit{Proof.} See the proof of Lemma 13 in \cite{Ma0}. \hfill \qedsymbol

\noindent \\ We introduce now the functional $\mathscr{M}$ we will use at the end of this paper to conclude the proof of Theorem 1; the remainder of this paper is dedicated to the obtention of the convergence $\mathscr{M} (s) \, \underset{s \to + \infty}{\longrightarrow} \, 0$.

\begin{leftbar}
\noindent \textbf{Lemma 13.} We define $\mathscr{M} := |b|^4 + || \rho v ||^2$. For all $s \geqslant 0$,
\begin{equation}
    | \dot{\mathscr{M}} | \leqslant C \left ( |b|^4 + || \rho \partial_y v ||^2 + || \rho v ||^2 \right ).
    \label{eqM}
\end{equation}
\end{leftbar}

\noindent \textit{Proof.} Follow the proof of Lemma 16 in \cite{Ma0}. Recall that $| \omega \partial_\omega V_1 | + | \omega \partial_\omega V_1 | \leqslant C \alpha ( \omega_0 )^{-1} \underline{\mathbf{V}} ( \omega_0 )$. Using the definition of $r_1$ and $r_2$, as well as \eqref{m}, we find
\[ || \rho r_1 || + || \rho r_2 || \leqslant C |m_\omega | \, |b| \alpha ( \omega_0 )^{-1} \underline{\mathbf{V}} ( \omega_0 ) \leqslant C (|b|^2 + || \nu v ||^2 ) \epsilon \alpha ( \omega_0 )^{-1} \underline{\mathbf{V}} ( \omega_0 ) \leqslant C ( |b|^2 + || \nu v ||^2 ) \]
as long as we take $\epsilon >0$ small enough (depending on $\omega_0$). The rest of the proof is unchanged and gives the desired result. \hfill \qedsymbol

\section{Estimate at large scale}
\noindent We will use virial arguments, which require suitable functions that will be denoted as follows. The arguments and notation used here originate from \cite{Ma2}, \cite{Ma3}, \cite{Ma1}. We fix a smooth even function $\chi \, : \, \R \to \R$ such $\chi = 1$ on $[0 \, , 1]$, $\chi = 0$ on $[2 \, , + \infty )$ and $\chi ' \leqslant 0$ on $[ 0 \, , + \infty )$. Let $1 \ll B \ll A$ be large constants to be fixed later. We define
\[ \begin{array}{ll} \chi_A (y) := \chi \left ( \frac{y}{A} \right ) , & \eta_A (y) := \text{sech} \left ( \frac{2y}{A} \right ) \\
\\ \zeta_A (y) := \exp \left ( - \frac{|y|}{A} (1 - \chi (y)) \right ) , & \Phi_A (y) := \int_0^y \zeta_A^2. \end{array} \]
Note that $0 < \Phi_A' = \zeta_A^2 \leqslant 1$, $|\Phi_A| \leqslant |y|$ and $| \Phi_A | \leqslant CA$ on $\R$. We define the function $\Psi_{A,B} := \chi_A^2 \Phi_B$ and the virial operators as follows:
\[ \Theta_A := 2 \Phi_A \partial_y + \Phi_A ' \, \, \, \, \, \, \, \, \text{and} \, \, \, \, \, \, \, \, \Xi_{A,B} := 2 \Psi_{A,B} \partial_y + \Psi_{A,B} '. \]
The first virial estimate is given below.

\begin{leftbar}
\noindent \textbf{Proposition 3.} For all $s>0$,
\[ \int_0^s \left ( || \eta_A \partial_y v||^2 + \frac{1}{A^2} || \eta_A v ||^2 \right ) \leqslant C \epsilon + C \int_0^s \left ( || \rho^4 v ||^2 + |b|^4 \right ). \]
\end{leftbar}

\noindent \textit{Proof.} Follow the proof of Lemma 18 in \cite{Ma0}. Definitions, estimates and differences are condensed here. By \eqref{sysv}, we see that $\frac{\text{d}}{\text{d}s} \int_{\R} ( \Theta_A v_2)v_1 = \sum\limits_{j=1}^5 \textbf{i}_j$ where
\[ \begin{array}{rl} & \textbf{i}_1 = - \int_{\R} ( \Theta_A v_1) \partial_y^2 v_1 - \int_{\R} ( \Theta_A v_2) \partial_y^2 v_2 , \\
\\ & \textbf{i}_2 = \int_{\R} ( \Theta_A v_1 ) \mu_1 + \int_{\R} ( \Theta_A v_2 ) \mu_2 , \\
\\ & \textbf{i}_3 = \int_{\R} ( \Theta_A v_1 ) p_1^{\T} + \int_{\R} ( \Theta_A v_2 ) p_2^{\perp} , \\
\\ & \textbf{i}_4 = - \int_{\R} ( \Theta_A v_1 ) r_1^{\T} - \int_{\R} ( \Theta_A v_2 ) r_2^{\perp} \\
\\ \text{and} & \textbf{i}_5 = - \int_{\R} ( \Theta_A v_1 ) (f_\omega ' ( Q_\omega ) v_1 + q_1^{\T} ) - \int_{\R} ( \Theta_A v_2 ) \left ( \frac{f_\omega ( Q_\omega )}{Q_\omega} v_2 + q_2^{\perp} \right ). \end{array} \]
We see that 
\[ |\textbf{i}_1| \geqslant 2 || \partial_y \widetilde{v}||^2 - C || \nu v ||^2 \, \, \, \, \, \, \, \text{and} \, \, \, \, \, \, \, | \textbf{i}_2| \leqslant C ( || \nu v ||^2 +  |b|^4 ), \]
where $\widetilde{v} := \zeta_A v$. As for $\textbf{i}_3$, the proof is also identical, based on the fact that the inequality $(y^2+1)|V_j''| + (|y|+1)|V_j'| + |V_j| \leqslant C \rho^8$ remains true in our case since $\alpha \geqslant \frac{4}{5} \alpha ( \omega_0 )$. We find:
\[ \left | \textbf{i}_3 + m_\omega \int_{\R} ( \Theta_A v_2 ) v_1 \right | \leqslant \frac{1}{2} || \partial_y \widetilde{v} ||^2 + C || \rho^4 v ||^2 + C |b|^4. \]
Now, as for $\textbf{i}_4$, we recall \eqref{dVomega}. This leads to
\[ | \Theta_A \omega \partial_\omega V_j | \leqslant C \underline{\mathbf{V}} ( \omega_0 ) (1+y^2) e^{- \alpha |y|} \leqslant C \underline{\mathbf{V}} ( \omega_0 ) \rho^8 (y) \]
then
\[ \left | \int_{\R} ( \Theta_A v_1 ) r_1^{\T} \right | \leqslant C \frac{\underline{\mathbf{V}} ( \omega_0 )}{\alpha ( \omega_0 )} |b| \, |m_\omega | \, || \rho^4 v || \leqslant C \frac{\underline{\mathbf{V}} ( \omega_0 )}{\alpha ( \omega_0 )} \, \epsilon ( || \nu v ||^2 + |b|^2 ) || \rho^4 v ||^2 \leqslant C ( || \rho^4 v ||^2 + |b|^4 ) \]
choosing $\epsilon >0$ small enough (depending on $\omega_0$). The same proof holds for the term containing $r_2$ and we end up with
\[ | \textbf{i}_4 | \leqslant C \left ( || \rho^4 v ||^2 + |b|^4 \right ). \]
Finally, as for $\textbf{i}_5$, the proof is analogous: we consider
\[ \begin{array}{rl} & \widetilde{q}_1 = \text{Re} \left [ f_\omega (Q_\omega + v) - f_\omega (Q_\omega) \right ], \\
\\ & \widetilde{q}_2 = \text{Im} \left [ f_\omega ( Q_\omega + v ) - f_\omega (Q_\omega ) \right ], \\
\\ & \check{q}_1 = \text{Re} \left [ f_\omega (Q_\omega + u) - f_\omega (Q_\omega + v) - f_\omega '(Q_\omega) (u_1-v_1) \right ] \\
\\ \text{and} & \check{q}_2 = \text{Im} \left [ f_\omega (Q_\omega + u) - f_\omega (Q_\omega + v) - i \, \frac{f_\omega ( Q_\omega )}{Q_\omega} (u_2 - v_2) \right ]. \end{array} \]
We recall that $f_\omega ( \psi ) = | \psi |^2 \psi + \frac{g ( \omega | \psi |^2 )}{\omega}$ and we introduce $F_\omega ( \psi ) := \frac{| \psi |^4}{4} + \frac{G( \omega | \psi |^2 )}{2 \omega^2}$. Integrating by parts, we find that
\[ \int_{\R} ( \Theta_A v_1) \widetilde{q}_1 + \int_{\R} ( \Theta_A v_2 ) \widetilde{q}_2 = \textbf{i}_{5,1} + \textbf{i}_{5,2} + \textbf{i}_{5,3} \]
where
\[ \begin{array}{rl} & \textbf{i}_{5,1} = -2 \text{Re} \int_{\R} \Phi_A ' ( F_\omega (Q_\omega + v) - F_\omega (Q_\omega ) - f_\omega (Q_\omega ) v ), \\
\\ & \textbf{i}_{5,2} = -2 \text{Re} \int_{\R} \Phi_A Q_\omega ' (f_\omega (Q_\omega + v) - f_\omega (Q_\omega ) - f_\omega ' (Q_\omega ) v ) \\
\\ \text{and} & \textbf{i}_{5,3} = \text{Re} \int_{\R} \Phi_A ' \overline{v} (f_\omega (Q_\omega + v) - f_\omega (Q_\omega )). \end{array} \]
Estimating these integrals follow the same steps as in the proof of Proposition 3 in \cite{Ri} (the equivalent integrals are $I_1$, $I_2$ and $I_3$). It leads to
\[ \begin{array}{rl} & \displaystyle{| \textbf{i}_{5,1}| + | \textbf{i}_{5,3}| \leqslant C \int_{\R} \Phi_A ' ( |v|^4 + Q_\omega^2 |v|^2 ) \leqslant C || \zeta_A v^2 ||^2 + C || \nu v ||^2} \\
\\ \text{and} & \displaystyle{| \textbf{i}_{5,2}| \leqslant C \int_{\R} \Phi_A |Q_\omega '| ( Q_\omega |v|^2 + |v|^3 ) \leqslant C || \nu v ||^2.} \end{array} \]
Therefore, $|\textbf{i}_{5,1}| + | \textbf{i}_{5,2} | + | \textbf{i}_{5,3} | \leqslant C \left ( || \nu v ||^2 + || \zeta_A v^2 ||^2 \right ) \leqslant C || \nu v ||^2 + CA \epsilon || \partial_y \widetilde{v}||$ thus, choosing $\epsilon >0$ small enough (depending on $A$),
\[ \left | \int_{\R} ( \Theta_A v_1 ) \widetilde{q}_1 \right | + \left | \int_{\R} ( \Theta_A v_2 ) \widetilde{q}_2 \right | \leqslant \frac{1}{2} || \partial_y \widetilde{v} ||^2 + C || \rho^4 v ||^2. \]
Now, let us deal with the terms $\check{q}_1$ and $\check{q}_2$. We start with $\check{q}_1$ and we will need three different cases. First, suppose that $|u| \leqslant \frac{Q_\omega}{4}$ and $|v| \leqslant \frac{Q_\omega}{4}$. Using Taylor expansions, we see that
\[ \check{q}_1 = \text{Re} ( k_{1,\omega} (u) - k_{1,\omega} (v)) + |u|^2 u_1 - |v|^2 v_1 + Q_\omega (|u|^2 - |v|^2) + 2 Q_\omega (u_1^2 - v_1^2) + \frac{1}{\omega} \left [ (Q_\omega + v_1) (\text{IR}_u - \text{IR}_v) + (u_1 - v_1 ) \text{IR}_u \right ] \]
where 
\[ \begin{array}{rl} & \text{IR}_u = \int_{\omega Q_\omega^2}^{\omega | Q_\omega + u|^2} ( \omega | Q_\omega + u |^2 - t) g''(t) \, \text{d}t, \\ \\ & \text{IR}_v = \int_{\omega Q_\omega^2}^{\omega | Q_\omega + v|^2} ( \omega | Q_\omega + v |^2 - t) g''(t) \, \text{d}t \\ \\ \text{and} & k_{1,\omega} (u) = |Q_\omega + u|^2 ( Q_\omega + u) - Q_\omega^3 - 3 Q_\omega^2 u_1 - i Q_\omega^2 u_2. \end{array} \]
The notation $k_{1,\omega}$ is taken from the proof of Lemma 18 in \cite{Ma0}. It is shown that $|k_{1,\omega} (u) - k_{1,\omega} (v)| \leqslant C |u-v| (Q_\omega (|u| + |v|) + |u|^2 + |v|^2)$. Now, we decompose $\text{IR}_u - \text{IR}_v$ as follows and we use the bound $|g''(t)| \leqslant C/t$:
\[ \begin{array}{rcl} |\text{IR}_u - \text{IR}_v | & \leqslant & \displaystyle{\left | \int_{\omega | Q_\omega +v|^2}^{\omega | Q_\omega + u|^2} ( \omega | Q_\omega + v|^2 - t) g''(t) \, \text{d}t \right | + \left | \int_{\omega Q_\omega^2}^{\omega | Q_\omega + u |^2} (2 \omega Q_\omega (u_1-v_2) + \omega ( |u|^2 - |v|^2 ) ) g''(t) \, \text{d}t \right |} \\
\\ & \leqslant & \displaystyle{\omega \left | \, |Q_\omega + u|^2 - |Q_\omega + v|^2 \right | \, \left | \int_{\omega | Q_\omega + v|^2}^{\omega | Q_\omega + u |^2} \frac{C \, \text{d}t}{t} \right | + C \omega \left ( Q_\omega |u_1-v_1| + ||u|^2 - |v|^2| \right ) \left | \int_{\omega Q_\omega^2}^{\omega |Q_\omega + u|^2} \frac{C \, \text{d}t}{t} \right |.} \end{array} \]
We have $\left | \int_{\omega |Q_\omega + v|^2}^{\omega |Q_\omega + u|^2} \frac{C \, \text{d}t}{t} \right | \leqslant C \left | \ln \left | 1 + \frac{u}{Q_\omega} \right |^2 - \ln \left | 1 + \frac{v}{Q_\omega} \right |^2 \right |$ where $\left | 1 + \frac{u}{Q_\omega} \right |^2 = 1 + \frac{2u_1}{Q_\omega} + \frac{|u|^2}{Q_\omega^2}$, with $\frac{2u_1}{Q_\omega} + \frac{|u|^2}{Q_\omega^2} \geqslant - \frac{1}{2}$ thanks to the hypothesis. Since $\ln (1 + \cdot)$ is $C$-Lipschitz on $\left [ - \frac{1}{2} \, , + \infty \right )$, we have
\[ \left | \ln \left | 1 + \frac{u}{Q_\omega} \right |^2 - \ln \left | 1 + \frac{v}{Q_\omega} \right |^2 \right | \leqslant C \left ( \frac{|u_1 - v_1|}{Q_\omega} + \frac{||u|^2 - |v|^2|}{Q_\omega^2} \right ). \]
Moreover, we see that $||u|^2-|v|^2| \leqslant |u-v| (|u| + |v|) \leqslant Q_\omega |u-v|$, $||Q_\omega + u|^2 - |Q_\omega + v|^2| \leqslant Q_\omega |u-v| + |u-v| (|u| + |v|) \leqslant CQ_\omega |u-v|$ and
\[ \left | \int_{\omega Q_\omega^2}^{\omega | Q_\omega + u|^2} \frac{C \, \text{d}t}{t} \right | \leqslant C \left | \ln \left | 1 + \frac{u}{Q_\omega} \right |^2 - \ln 1 \right | \leqslant \frac{C |u|}{Q_\omega}. \]
Gathering all the previous estimates, we find
\[ |\text{IR}_u - \text{IR}_v | \leqslant C \omega |u-v| (|u| + |v|). \]
Moreover,
\[ | \text{IR}_u | \leqslant \left | \int_{\omega Q_\omega^2}^{\omega | Q_\omega + u|^2} ( \omega | Q_\omega + u |^2 - t ) g''(t) \, \text{d}t \right | \leqslant \omega ||Q_\omega + u|^2 - Q_\omega^2| \, \left | \int_{\omega Q_\omega^2}^{\omega |Q_\omega + u|^2} \frac{C \, \text{d}t}{t} \right | \leqslant C \omega \cdot Q_\omega |u| \cdot \frac{C |u|}{Q_\omega} \leqslant C \omega |u|^2. \]
Getting back to $\check{q}_1$, we evidently find that
\[ | \check{q}_1 | \leqslant C |u-v| \left ( Q_\omega ( |u| + |v|) + |u|^2 + |v|^2 \right ). \]
Now, let us consider the case $|u| \geqslant \frac{Q_\omega}{4}$. The situation is easier. We write
\[ \begin{array}{rcl} \check{q}_1 &=& \displaystyle{\text{Re} ( k_{1,\omega} (u) - k_{1,\omega} (v)) + \frac{Q_\omega + v_1}{\omega} \left ( g(\omega |Q_\omega+u|^2) - g( \omega | Q_\omega + v|^2) \right ) + \frac{u_1 - v_1}{\omega} g(\omega | Q_\omega + u|^2)} \\ \\ & & \, \, \, \, \, - (u_1 - v_1) \left ( 2 Q_\omega^2 g'( \omega Q_\omega^2 ) + \frac{g( \omega Q_\omega^2 )}{\omega} \right ). \end{array} \]
Since $|g'|$ is bounded, we have $| g( \omega | Q_\omega + u|^2 ) - g ( \omega | Q_\omega + v |^2 ) | \leqslant C \omega ||Q_\omega + u|^2 - |Q_\omega + v|^2| \leqslant C \omega |u-v| (Q_\omega + |u| + |v|)$. We also have, thanks to the hypothesis, $|Q_\omega + v_1| \leqslant C |v|$, $|g(\omega |Q_\omega + u|^2 )| \leqslant C \omega |Q_\omega + u|^2 \leqslant C \omega |u|^2$ and 
\[ \left | (u_1 - v_1) \left ( 2 Q_\omega^2 g'( \omega Q_\omega^2 ) + \frac{g( \omega Q_\omega^2 )}{\omega} \right ) \right | \leqslant |u-v| (CQ_\omega^2 + CQ_\omega^2) \leqslant C |u-v| \, |u|^2. \]
Therefore, gathering these estimates, we obtain:
\[ | \check{q}_1 | \leqslant  C |u-v| \left ( Q_\omega ( |u| + |v|) + |u|^2 + |v|^2 \right ). \]
The last case, namely $|v| \geqslant \frac{Q_\omega}{4}$, is treated analogously. The estimate above holds in any case. \\
\\ Now, as for $\check{q}_2$, we see that
\[ \check{q}_2 = \text{Im} (k_{1,\omega} (u) - k_{1,\omega} (v)) + \frac{Q_\omega + u_2}{\omega} g( \omega | Q_\omega + u|^2 ) - \frac{Q_\omega + v_2}{\omega} g ( \omega | Q_\omega + v|^2) - \frac{u_2-v_2}{\omega} g ( \omega Q_\omega^2) \]
where we already know that $|\text{Im} (k_{1,\omega} (u) - k_{1,\omega} (v))| \leqslant C |u-v| \left ( Q_\omega ( |u| + |v|) + |u|^2 + |v|^2 \right )$. For the remaining term, in the case $|u| \geqslant \frac{Q_\omega}{4}$ or $|v| \geqslant \frac{Q_\omega}{4}$, the proof is just as we did for $\check{q}_1$: $\left | \frac{u_2 - v_2}{\omega} g ( \omega Q_\omega^2 ) \right | \leqslant C Q_\omega^2 |u-v| \leqslant C |u|^2 |u-v|$ and $\left | \frac{Q_\omega + u_2}{\omega} g( \omega | Q_\omega + u|^2 ) - \frac{Q_\omega + v_2}{\omega} g( \omega |Q_\omega + v|^2 ) \right | \leqslant C |u-v| (|u|^2 + |v|^2)$. Therefore, in these cases,
\[ | \check{q}_2 | \leqslant  C |u-v| \left ( Q_\omega ( |u| + |v|) + |u|^2 + |v|^2 \right ). \]
Let us consider the remaining case: $|u| \leqslant \frac{Q_\omega}{4}$ and $|v| \leqslant \frac{Q_\omega}{4}$. Then, using Taylor expansions, we see that
\[ \begin{array}{rl} & \displaystyle{\frac{Q_\omega + u_2}{\omega} g( \omega | Q_\omega + u|^2 ) - \frac{Q_\omega + v_2}{\omega} g ( \omega | Q_\omega + v|^2) - \frac{u_2-v_2}{\omega} g ( \omega Q_\omega^2)} \\ \\ =& \displaystyle{g'( \omega Q_\omega^2 ) \, \text{Im} \left ( k_{1,\omega} (u) - k_{1,\omega} (v) \right ) + (Q_\omega + u_2) \frac{\text{IR}_u}{\omega} - (Q_\omega + v_2) \frac{\text{IR}_v}{\omega}} \end{array} \]
where $\text{IR}_u$ and $\text{IR}_v$ are the same integrals as before. We recall that $|g'|$ is bounded and we deal with the rest of the expression above as we did for $\check{q}_1$: the conclusion is the same and we have
\[ | \check{q}_2 | \leqslant  C |u-v| \left ( Q_\omega ( |u| + |v|) + |u|^2 + |v|^2 \right ) \]
in this last case too. Therefore, we have achieved to establish the estimate
\[ | \check{q}_1 | + | \check{q}_2 | \leqslant C |u-v| \left ( Q_\omega (|u| + |v|) + |u|^2 + |v|^2 \right ) \]
which leads to
\[ | \textbf{i}_5 | \leqslant || \partial_y \widetilde{v} ||^2 + C || \rho^4 v ||^2 + C |b|^4. \]
Setting $\textbf{I} := \omega \int_{\R} (\Theta_A v_2)v_1$, we have, on one hand,
\[ \frac{\text{d} \textbf{I}}{\text{d}s} \geqslant \omega_0 \left ( \frac{1}{2} || \partial_y \widetilde{v} ||^2 - C || \rho^4 v ||^2 - C |b|^4 \right ) \]
and, on the other hand, $| \textbf{I} (s) | \leqslant \omega_0 A ||v||_{H^1 ( \R)}^2 \leqslant \omega_0 A \epsilon^2$. The end of the proof is identical and leads to the desired result:
\[ \int_0^s \left ( || \eta_A \partial_y v ||^2 + \frac{1}{A^2} || \eta_A v ||^2 \right ) \leqslant C \epsilon + C \int_0^s ( || \rho^4 v ||^2 +  |b|^4 ). \]
\hfill \qedsymbol

\section{The Fermi golden rule}
\noindent For the Fermi golden rule, we need to construct a non trivial bounded solution $(g_1 \, , g_2)$ of
\begin{equation}
    \left \{ \begin{array}{ccl} L_+ g_1 &=& 2 \lambda g_2 \\ L_- g_2 &=& 2 \lambda g_1. \end{array} \right.
    \label{sysg}
\end{equation}
The reason for such a construction is the following. The nonlinear terms in the system \eqref{sysv} satisfied by $(v_1 \, , v_2)$ are, at main order, quadratic in $b$. Recall from the introduction that $b$ is to be thought as a function of time oscillating with pulsation $\lambda$. In that case, quantities quadratic in $b$ (such as $b_1b_2$, $b_1^2$, $b_2^2$) contain time-dependent functions oscillating with pulsation $2 \lambda$. This is, heuristically, why we consider the spectral problem \eqref{sysg} above. The pair $(g_1 \, , g_2)$ shall be used in the proof of Proposition 4 in order to control $\int_0^s |b|^4$. \\
\\ To construct the pair $(g_1 \, , g_2)$, notice that if $h_1$ satisfies $M_- M_+ h_1 = 4 \lambda^2 h_1$ then, setting $g_1 = (S^*)^2 h_1$ and $g_2 = \frac{1}{2 \lambda} L_+ g_1$, $(g_1 \, , g_2)$ satisfies \eqref{sysg} thanks to the relation $S^2 L_+ L_- = M_+ M_- S^2$. 

\begin{leftbar}
\noindent \textbf{Lemma 14.} Let $\tau := \sqrt{2 \lambda -1}$. For $\omega >0$ small enough, there exist smooth even functions $h_1$ and $h_2$ (depending on $\omega$) that satisfy
\begin{equation}
    \left \{ \begin{array}{ccl} M_+ h_1 &=& 2 \lambda h_2 \\ M_- h_2 &=& 2 \lambda h_1 \end{array} \right.
    \label{sysh}
\end{equation}
and, for all $k \in \N$, 
\[ | \partial_y^k (h_1 + \cos ( \tau y)) | + | \partial_y^k ( h_2 + \cos ( \tau y )) | \leqslant C \varepsilon_\omega \, \, \, \, \,  \,  \, \text{and} \, \, \, \, \, \, \, | \partial_y^k \partial_\omega h_1 | + | \partial_y^k \partial_\omega h_2 | \leqslant \frac{C \varepsilon_\omega}{\omega} (1 + |y|) \]
on $\R$. Setting $g_1 = (S^*)^2 h_1$ and $g_2 = \frac{1}{2 \lambda} L_+ g_1$, the pair $(g_1 \, , g_2)$ satisfies \eqref{sysg} and, for all $k \in \{ 0 \, , ... \,  ,2 \}$,
\[ \begin{array}{rl} & \left | \partial_y^k \left ( g_1 - \left ( \frac{2Q'}{Q} \sin ( \tau y ) + Q^2 \cos ( \tau y ) \right ) \right ) \right | + \left | \partial_y^k \left ( g_2 - \frac{2Q'}{Q} \sin ( \tau y ) \right ) \right | \leqslant C \varepsilon_\omega \\
\\ \text{and} & | \partial_\omega g_1 | + | \partial_\omega g_2 | + | \partial_y \partial_\omega g_1 | + | \partial_y \partial_\omega g_2 | \leqslant \frac{C}{\omega} ( 1 + |y| ) \end{array} \]
on $\R$. Moreover, the following orthogonality relations hold:
\begin{equation}
    \langle g_1 \, , Q_\omega \rangle = \langle g_2 \, , \Lambda_\omega Q_\omega \rangle = \langle g_1 \, , V_2 \rangle = \langle g_2 \, , V_1 \rangle = 0.
    \label{orthg}
\end{equation}
\end{leftbar}

\noindent \textit{Proof.} Follow the proof of Lemma 19 in \cite{Ma0}. Setting $\ell_1 := \frac{h_1 + h_2}{2}$ and $\ell_2 := \frac{h_1 - h_2}{2}$, we look for $(\ell_1 \, , \ell_2)$ satisfying
\[ \left \{ \begin{array}{l} - \ell_1 '' - ( 2 \lambda - 1) + b_\omega^+ \ell_1 + b_\omega^- \ell_2 = 0 \\ - \ell_2 '' + (2 \lambda + 1) \ell_2 + b_\omega^- \ell_1 + b_\omega^+ \ell_2 = 0. \end{array} \right. \]
Let $\check{\ell}_1 := \ell_1 + \cos ( \tau y )$ and $\check{\ell}_2 := \ell_2$. We look for $(\check{\ell}_1 \, , \check{\ell}_2)$ satisfying
\begin{equation}
    \left \{ \begin{array}{l} - \check{\ell}_1 '' - \tau^2 \check{\ell}_1 = - b_\omega^+ \check{\ell}_1 - b_\omega^- \check{\ell}_2 + b_\omega^+ \cos ( \tau y ) \\ - \check{\ell}_2 '' + (2 + \tau^2) \check{\ell}_2 = - b_\omega^- \check{\ell}_1 - b_\omega^+ \check{\ell}_2 + b_\omega^- \cos ( \tau y). \end{array} \right.
    \label{sysl}
\end{equation}
We define a bounded linear map $\check{\Upsilon} \, : \, ( \mathscr{C}_b ( \R ))^2 \to ( \mathscr{C}_b ( \R ))^2$, where $\mathscr{C}_b ( \R )$ is the space of bounded continuous functions on $\R$ equipped with the supremum norm $|| \cdot ||_\infty$, by setting
\[ \check{\Upsilon} \left ( \begin{array}{c} \check{\ell}_1 \\ \check{\ell}_2 \end{array} \right ) = \left ( \begin{array}{c} - \frac{1}{\tau} \int_0^y \sin (\tau (y-y')) (-b_\omega^+ \check{\ell}_1 - b_\omega^- \check{\ell}_2 ) (y') \, \text{d}y' \\ \\ \frac{1}{2 \sqrt{2 + \tau^2}} \int_{\R} e^{- \sqrt{2 + \tau^2} |y-y'|} (-b_\omega^- \check{\ell}_1 - b_\omega^+ \check{\ell}_2 ) (y') \, \text{d}y' \end{array} \right ). \]
We also define
\[ \begin{array}{rl} & \displaystyle{\check{f}_1 := - \frac{1}{\tau} \int_0^y \sin ( \tau (y-y')) b_\omega^+ (y') \cos ( \tau y') \, \text{d}y'} \\ \\ \text{and} & \displaystyle{\check{f}_2 := \frac{1}{2 \sqrt{2 + \tau^2}} \int_{\R} e^{- \sqrt{2 + \tau^2} |y-y'|} b_\omega^- (y') \cos ( \tau y' ) \, \text{d}y'.} \end{array} \]
That way, the integral formulation of the system \eqref{sysl} (for even functions satisfying $\check{\ell}_1 (0)=0$ by convention) is
\begin{equation}
    \left ( \begin{array}{c} \check{\ell}_1 \\ \check{\ell}_2 \end{array} \right ) = \check{\Upsilon} \left ( \begin{array}{c} \check{\ell}_1 \\ \check{\ell}_2 \end{array} \right ) + \left ( \begin{array}{c} \check{f}_1 \\ \check{f}_2 \end{array} \right ).
    \label{eqUpsilon}
\end{equation}
We easily see that $|| \check{f}_1 ||_\infty + || \check{f}_2 ||_\infty + ||| \check{\Upsilon} ||| \leqslant C \varepsilon_\omega$. Thus, the operator $\text{Id} - \check{\Upsilon}$ is invertible for $\omega >0$ small enough, and \eqref{eqUpsilon} becomes
\[ \left ( \begin{array}{c} \check{\ell}_1 \\ \check{\ell}_2 \end{array} \right ) = ( \text{Id} - \check{\Upsilon} )^{-1} \left ( \begin{array}{c} \check{f}_1 \\ \check{f}_2 \end{array} \right ) = \sum_{j=0}^{+ \infty} \left ( \begin{array}{c} \check{f}_1^j \\ \check{f}_2^j \end{array} \right ), \]
where $\left ( \begin{array}{c} \check{f}_1^j \\ \check{f}_2^j \end{array} \right ) := \check{\Upsilon}^j \left ( \begin{array}{c} \check{f}_1 \\ \check{f}_2 \end{array} \right )$. This proves the existence of a solution $( \check{\ell}_1 \, , \check{\ell}_2 )$ of \eqref{eqUpsilon}. Using the estimates $| \partial_\omega \tau | \leqslant \frac{C \alpha \varepsilon_\omega}{\omega}$ and \eqref{daomega}, we also find that $| \partial_\omega \check{f}_1 | + | \partial_\omega \check{f}_2 | \leqslant \frac{C \varepsilon
_\omega}{\omega} (1 +  |y|)$. From there, reasoning by induction, we show that, for all $j \in \N$,
\[ | \check{f}_1^j | + | \check{f}_2^j | \leqslant C \varepsilon_\omega^{j+1} \, \, \, \, \, \, \, \text{and} \, \, \, \, \, \, \, | \partial_\omega \check{f}_1^j | + | \partial_\omega \check{f}_2^j | \leqslant C \frac{\varepsilon_\omega^{j+1}}{\omega} ( 1+ |y|). \]
Differentiating with regards to $y$ only differentiates the term $\sin ( \tau (y-y'))$ or $e^{- \sqrt{2 + \tau^2} |y-y'|}$ in the integral. This does not change the estimates obtained previously, therefore, for all $j \in \N$ and all $k \in \N$,
\[ | \partial_y^k \check{f}_1^j | + | \partial_y^k \check{f}_2^j | \leqslant C \varepsilon_\omega^{j+1} \, \, \, \, \, \, \, \text{and} \, \, \, \, \, \, \, | \partial_y^k \partial_\omega \check{f}_1^j | + | \partial_y^k \partial_\omega \check{f}_2^j | \leqslant C \frac{\varepsilon_\omega^{j+1}}{\omega} ( 1+ |y|). \]
Getting back to the Neumann expansion, we get that, for all $k \in \N$,
\[ | \partial_y^k \check{\ell}_1 | + | \partial_y^k \check{\ell}_2 | \leqslant C \varepsilon_\omega \, \, \, \, \, \, \, \text{and} \, \, \, \, \, \, \, | \partial_y^k \partial_\omega \check{\ell}_1 | + | \partial_y^k \partial_\omega \check{\ell}_2 | \leqslant \frac{C \varepsilon_\omega}{\omega} (1 + |y|). \]
Until the end of this proof, $\tilde{\mathcal{O}}_p ( \varepsilon_\omega )$ denotes any function $\ell$ such that $| \partial_y^k \ell | \leqslant C_k \varepsilon_\omega$ for all $k \in \{ 0 \, , ... \, , p \}$. We define $h_1 = - \cos ( \tau y ) + ( \check{\ell}_1 + \check{\ell}_2 )$, $h_2 = - \cos ( \tau y) + ( \check{\ell}_1 + \check{\ell}_2 )$, $g_1 = (S^*)^2 h_1$ and $g_2 = \frac{1}{2 \lambda} L_+ g_1$. We check that $(h_1 \, , h_2)$ satisfies \eqref{sysh}, $(g_1 \, , g_2)$ satisfies \eqref{sysg}, and we have $h_1 = - \cos ( \tau y ) + \tilde{\mathcal{O}}_\infty ( \varepsilon_\omega )$ and $h_2 = - \cos ( \tau y ) + \tilde{\mathcal{O}}_\infty ( \varepsilon_\omega )$. Using the bounds $\left | \frac{Q_\omega '}{Q_\omega} - \frac{Q'}{Q} \right | = \tilde{\mathcal{O}}_5 ( \varepsilon_\omega )$, $\left | \frac{g ( \omega Q_\omega^2 )}{\omega} \right | = \tilde{\mathcal{O}}_4 ( \varepsilon_\omega )$, $| Q_\omega^2 - Q^2 | = \tilde{\mathcal{O}}_6 ( \varepsilon_\omega )$, $| \tau^2 -1 | \leqslant C \varepsilon_\omega^2$ and $| \tau -1 | \leqslant C \varepsilon_\omega^2$, we compute:
\[ \begin{array}{rcl} g_1 &=& \frac{Q_\omega ''}{Q_\omega} h_1 + \frac{2Q_\omega '}{Q_\omega} h_1' + h_1'' \\
\\ &=& ( 1-Q^2) h_1 + \frac{2Q'}{Q} h_1 ' + h_1'' + \tilde{\mathcal{O}}_4 ( \varepsilon_\omega ) \\
\\ &=& Q^2 \cos ( \tau y ) + \frac{2Q'}{Q} \sin ( \tau y) + \tilde{\mathcal{O}}_4 ( \varepsilon_\omega ). \end{array} \]
Then, using $| \lambda^{-1} -1 | \leqslant C \varepsilon_\omega^2$ and estimates like the ones above, we compute:
\[ g_2 = \frac{1}{2} \left ( -g_1 '' + g_1 - 3 Q^2 g_1 \right ) + \tilde{\mathcal{O}}_4 ( \varepsilon_\omega ) = \frac{2Q'}{Q} \sin ( \tau y ) + \tilde{\mathcal{O}}_2 ( \varepsilon_\omega ). \]
Differentiating with regards to $\omega$ the formulas $g_1 = (S^*)^2h_1$ and $g_2 = \frac{1}{2 \lambda} L_+ g_1$ and using estimates found in Proposition 2 (such as $| \partial_\omega Q_\omega | \leqslant \frac{C}{\omega}$ or $| \alpha ' ( \omega ) | \leqslant \frac{C \varepsilon_\omega}{\omega}$), we ultimately find that
\[ | \partial_y^k \partial_\omega g_1 | + | \partial_y^k \partial_\omega g_2 | \leqslant \frac{C}{\omega} (1 + |y|) \]
for any $k \in \{ 0 \, , 1 \}$ and all $y \in \R$. Finally, the orthogonality relations are proven using $L_+ \Lambda_\omega Q_\omega = - Q_\omega$, $L_- Q_\omega = 0$ and the equations of $(V_1 \, , V_2)$ and $(g_1 \, , g_2)$. \hfill \qedsymbol

\noindent \\ We now define
\[ \begin{array}{ll} G := V_1^2 Q_\omega ( 3 + 3 g'(\omega Q_\omega^2) + 2 \omega Q_\omega^2 g''( \omega Q_\omega^2 )) , & H := V_2^2 Q_\omega ( 1 + g'( \omega Q_\omega^2 )), \\ \\ G_1 = G-H, & G_2 = 2 V_1 V_2 Q_\omega (1 + g'( \omega Q_\omega^2 )), \\ \\ G_1^\T = G_1 - \frac{\langle G_1 \, , V_1 \rangle}{\langle V_1 \, , V_2 \rangle} V_2 , & G_2^\perp = G_2 - \frac{\langle G_2 \, , V_2 \rangle}{\langle V_1 \, , V_2 \rangle} V_1. \end{array} \]
The quantity $G$ above must not be confused with the function $G$. We keep the notation $G$ in order to fit the notation of \cite{Ma0}. As we will mostly have to deal with $G_1$ and $G_2$ (rather than $G$ itself), there should be no confusion. Lastly, we define 
\[ \Gamma ( \omega ) = \int_{\R} ( G_1^\T g_1 + G_2^\perp g_2). \]
The hypothesis $(H_3)$ presented in the introduction can be reformulated as follows:
\begin{equation} 
    \begin{array}{rcc} (H_3) & : & \text{there exists a positive quantity $\underline{\mathbf{\Gamma}} ( \omega_0 )$ depending only on $\omega_0$ such that,} \\ \\ & & | \omega - \omega_0 | \leqslant \frac{\omega_0}{2} \, \, \, \, \Longrightarrow \, \, \, \, \Gamma ( \omega ) \geqslant \underline{\mathbf{\Gamma}} ( \omega_0 ) > 0. \end{array}
    \label{H3}
\end{equation}
This hypothesis appears to be hard to verify, but \cite{Ma0} proves that it holds in the case $g(s) = s^2$. We shall investigate the case $g(s) = s^\sigma$ for $\sigma > 1$ a little further. For now, let us operate a simplification of $\Gamma ( \omega )$. 

\begin{leftbar}
\noindent \textbf{Lemma 15.} For $\omega > 0$ small enough, we have
\[ \Gamma ( \omega ) = \int_{\R} \left ( Q^2 \Delta_4 \cos y + \frac{2Q'}{Q} ( \Delta_4 + 2 \Delta_2 ) \sin y \right ) \, \text{d}y + \mathcal{O}( \varepsilon_\omega^2  ), \]
where 
\[ \begin{array}{rl} & \Delta_4 = 6Q (1-Q^2)R_1 + (1-Q^2)^2 \left ( 3 D_\omega + Q(3 g'(\omega Q^2) + 2 \omega Q^2 g''( \omega Q^2 )) \right ) -2Q R_2 - 3D_\omega - Q g'( \omega Q^2 ) \\ \\ & \, \, \, \, \, \, \, \, \, \, \, \, \, \, \, \, \, \, \, \, \, \, \, - \left ( \frac{g( \omega Q^2 )}{\omega} + 2 Q^2 g'( \omega Q^2) + 6Q D_\omega \right ) Q(1-Q^2) \\
\\ \text{and} & \Delta_2 = QR_1 + Q(1-Q^2) R_2 + (1-Q^2) (D_\omega + Q g'( \omega Q^2)). \end{array} \]
\end{leftbar}

\noindent \textit{Proof.} Follow the proof of Lemma 20 in \cite{Ma0}, in which a similar result is obtained. In what follows, $\tilde{\mathcal{O}}_p ( \varepsilon_\omega^2 )$ denotes any function $\ell$ such that $| \partial_y^k \ell | \leqslant C_k \varepsilon_\omega^2 (1 + y^2)$ for all $k \in \{ 0 \, , ... \, , p \}$. We first establish the following expansions, using Proposition 2: 
\[ \begin{array}{ccl} G &=& 3Q(1-Q^2)^2 + 6Q(1-Q^2)R_1 + 3D_\omega (1-Q^2)^2 + (3 g' (\omega Q^2 ) + 2 \omega Q^2 g''( \omega Q^2)) Q (1-Q^2)^2 + Q \, \tilde{\mathcal{O}}_2 ( \varepsilon_\omega^2 ) , \\
\\ H &=& Q + D_\omega + Q g'( \omega Q^2 ) + 2 Q R_2 + Q \, \tilde{\mathcal{O}}_2 ( \varepsilon_\omega^2 ) , \\
\\ G_1 &=& 3Q(1-Q^2)^2 - Q + \Delta_1 + Q \, \tilde{\mathcal{O}} ( \varepsilon_\omega^2 ) , \\
\\ G_2 &=& 2Q(1-Q^2) + 2 \Delta_2 + Q \, \tilde{\mathcal{O}} ( \varepsilon_\omega^2 ), \end{array} \]
where
\[ \Delta_1 := 6Q (1-Q^2)R_1 - 2 QR_2 - D_\omega - Q g'(\omega Q^2) + Q(1-Q^2)^2 \left ( 3 D_\omega + (3 g'(\omega Q^2) + 2 \omega Q^2 g''( \omega Q^2) \right ) \]
and $\Delta_2$ has the expression announced in the lemma. Then, recalling that $\lambda = 1 + \mathcal{O} ( \varepsilon_\omega^2 )$, we compute
\[ \begin{array}{rcl} \displaystyle{G_1 + \frac{1}{2 \lambda} L_+ G_2} &=& \displaystyle{G_1 + \frac{1}{2} \left [ -G_2 '' + G_2 - 3Q^2 G_2 - G_2 \left ( \frac{g ( \omega Q^2)}{\omega} + 2 Q^2 g'( \omega Q^2 ) + 6 QD_\omega \right ) + \tilde{\mathcal{O}}_0 ( \varepsilon_\omega^2 ) \right ]} \\
\\ &=& \displaystyle{2Q + \Delta_3 + Q \, \tilde{\mathcal{O}}_0 ( \varepsilon_\omega^2 ),} \end{array} \]
where
\[ \begin{array}{rcl} \Delta_3 &:=& \Delta_1 - \Delta_2 '' + \Delta_2 - 3Q^2 \Delta_2 - \left ( \frac{g ( \omega Q^2)}{\omega} + 2 Q^2 g'( \omega Q^2 ) + 6 QD_\omega \right ) Q (1-Q^2) \\ \\ &=& 2 D_\omega + (- \Delta_2 '' + \Delta_2 -3 Q^2 \Delta_2 ) + \Delta_4 \end{array} \]
where $\Delta_4$ turns out to be the quantity presented in the lemma. We find that
\[ \Gamma ( \omega ) = \int_{\R} g_1 ( \Delta_3 - 2 D_\omega ) + \mathcal{O}( \varepsilon_\omega^2 ) = \int_{\R} g_1 \Delta_4 + 2 \int_{\R} g_2 \Delta_2 + \mathcal{O} ( \varepsilon_\omega^2 ).  \]
We now use the expansions of $g_1$ and $g_2$ proven in Lemma 14. Noticing that $| \cos ( \tau y ) - \cos (y) | \leqslant |y| \, | \tau -1 | \leqslant C \varepsilon_\omega^2 |y|$ and that a similar estimate holds for $\sin$, we see that
\[ \int_{\R} Q^2 \cos ( \tau y ) \Delta_4 = \int_{\R} Q^2 \cos (y) \Delta_4 + \mathcal{O} ( \varepsilon_\omega^2 ) \]
for example. Combining these developments, we find the wanted formula:
\[ \Gamma ( \omega ) = \int_{\R} \left ( Q^2 \Delta_4 \cos y + \frac{2Q'}{Q} ( \Delta_4 + 2 \Delta_2 ) \sin y \right ) \, \text{d}y + \mathcal{O}( \varepsilon_\omega^2  ). \]
\hfill \qedsymbol

\noindent \\ It does not seem possible to go quite further from here for a general function $g$. But one can go further by considering the interesting and useful case $g(s) = s^\sigma$ where $\sigma >1$. This constitutes a generalisation of \cite{Ma0}, which deals with the case $\sigma = 2$.

\begin{leftbar}
\noindent \textbf{Lemma 16.} Here, we take $g(s) = as^\sigma$ with $\sigma >1$ and $a>0$. For $\omega > 0$ small enough, we have
\[ \Gamma ( \omega ) = a \, \Gamma_0 ( \sigma ) \, \omega^{\sigma-1} + \mathcal{O} \left ( \omega^{2(\sigma -1)} \right ) \]
where
\[ \Gamma_0 ( \sigma ) := \int_{\R} \left ( Q^2 \Delta_4^0 \cos y + \frac{2Q'}{Q} \Delta_5^0 \sin y \right ) \, \text{d}y  \]
with
\[ \begin{array}{rl} & \Delta_4^0 := (13Q^2 -16)Q^2 D^0 - 8 Q'Q(D^0)' + 2Q (3(1-Q^2)^2-1) T_1^0 \\ \\ & \, \, \, \, \, \, \, \, \, \, \, \, \, \, \, \, \, \, \, \,  + 6 Q(2-Q^2)^2 T_2^0 +4(2-3Q^2)Q' (T_1^0)' + 4 (4-3Q^2)Q' (T_2^0)' \\ \\ & \, \, \, \, \, \, \, \, \, \, \, \, \, \, \, \, \, \, \, \, + 2 \sigma^2 Q^{2 \sigma -1}- \left ( \frac{4( \sigma +2)}{\sigma +1} + (2 \sigma +1)^2 \right ) Q^{2 \sigma +1} + \left ( \frac{8}{\sigma +1} + ( \sigma +1) ( 2 \sigma +1) \right ) Q^{2 \sigma +3} , \\
\\ & \Delta_5^0 := (2-30Q^2 + 29Q^4 - 8 Q^6) D^0 -8 Q^3 Q' (D^0)' + 2Q(2-Q^2)(2-3Q^2) T_1^0 \\ \\ & \, \, \, \, \, \, \, \, \, \, \, \, \, \, \, \, \, \, \, \,  + 2Q(3Q^4-10Q^2+12)T_2^0 + 16(1-Q^2)Q' (T_1^0)' + 8(2-Q^2)Q' (T_2^0)' \\
\\ & \, \, \, \, \, \, \, \, \, \, \, \, \, \, \, \, \, \, \, \, +2 \sigma ( \sigma +1) Q^{2 \sigma -1} - \left ( \frac{16}{\sigma +1} + 2 \sigma + (2 \sigma +1)^2 \right ) Q^{2 \sigma +1} + \left ( \frac{4 ( 4 - \sigma )}{\sigma +1} + ( \sigma +1) (2 \sigma +1) \right ) Q^{2 \sigma +3} - \frac{4}{\sigma +1} Q^{2 \sigma +5} , \\
\\ & T_1^0(y) := - \frac{( \sigma -1)^2}{2 ( \sigma +1)} \int_{\R} |y-z| Q^{2 \sigma} (z) \, \text{d}z , \\
\\ & T_2^0 (y) := - \frac{\sqrt{2} \, \sigma ( \sigma -1)}{4 ( \sigma + 1)} \int_{\R} e^{- \sqrt{2} |y-z|} Q^{2 \sigma} (z) \, \text{d}z , \\
\\ & D^0 (y) := -Q'(y) \int_0^y AQ^{2 \sigma +1} + A(y) \frac{Q(y)^{2 \sigma +2}}{2 \sigma +2} \\
\\ \text{and} & A(y) := \frac{\sqrt{2}}{4 \, \text{cosh} (y)} \left ( 3 y  \, \text{tanh} (y) + \text{sinh}^2 (y) -2 \right ). \end{array} \]
\end{leftbar}

\noindent \textit{Proof.} We start with the following relation:
\[ L_+^0 D_\omega = a \omega^{\sigma-1} Q^{2 \sigma +1} + Z_\omega \]
where $Z_\omega := D_\omega^2 ( Q_\omega + 2Q) + a \omega^{\sigma -1} ( Q_\omega^{2 \sigma +1} - Q^{2 \sigma +1})$. We check that $| Z_\omega | \leqslant C_\sigma \omega^{2 ( \sigma -1)} e^{-3 |y|}$. Besides, we know how to invert the operator $L_+^0 = - \partial_y^2 + 1-3Q^2$, it is similar to the operator $I_+$ in \cite{Ma1}. We have
\[ (L_+^0)^{-1} [W] (y) = \left | \begin{array}{ll} -Q'(y) \int_0^y AW - A(y) \int_y^{+ \infty} Q' W & \text{if $y \geqslant 0$} \\ \\ Q'(y) \int_y^0 AW + A(y) \int_{- \infty}^y Q'W & \text{if $y<0$} \end{array} \right. \]
where $A$ denotes the even solution of $L_+^0 A=0$ such that $Q''A-Q'A'=1$ on $\R$. This solution is not bounded and verifies $| A^{(k)} | \leqslant C_k e^{|y|}$ on $\R$. Actually, we can compute $A$ explicitly:
\[ A(y) = \frac{\sqrt{2}}{4 \, \text{cosh} (y)} \left ( 3 y  \, \text{tanh} (y) + \text{sinh}^2 (y) -2 \right ). \]
This leads, for $y>0$, to
\[ D_\omega (y) = a \omega^{\sigma -1} D^0 (y) + \widetilde{D}_\omega (y), \]
where
\[ \begin{array}{rrcl} & D^0 (y) &=& - Q'(y) \int_0^y AQ^{2 \sigma +1} - A(y) \int_y^{+ \infty} Q' Q^{2 \sigma +1} \\
\\  & & = & -Q'(y) \int_0^y A Q^{2 \sigma +1} + A(y) \frac{Q(y)^{2 \sigma +2}}{2 \sigma + 2} \\
\\ \text{and} & \widetilde{D}_\omega (y) &=& -Q'(y) \int_0^y A Z_\omega - A(y) \int_y^{+ \infty} Q' Z_\omega. \end{array} \]
Using the bounds on $Z_\omega$, $Q'$ and $A$, we see that $| \widetilde{D}_\omega (y) | \leqslant C_\sigma \omega^{2 ( \sigma -1)} e^{-|y|}$. Similar estimates hold for $y<0$ (and anyway $D_\omega$ is even). \\
\\ The expressions of $R_1$ and $R_2$ involve $T_1$ and $T_2$, which involve $Q_\omega$. As we did before with many other expressions, we can replace these $Q_\omega$ by $Q$, at a cost of $\varepsilon_\omega^2$. We have:
\[ \begin{array}{ccl} T_1 (y) &=& \displaystyle{\frac{1}{2} \int_{\R} |y-z| \left ( 3 \frac{g ( \omega Q^2 )}{\omega} - 4  \frac{G ( \omega Q^2)}{\omega^2 Q^2} - Q^2 g'( \omega Q^2 ) \right ) \, \text{d}z + \tilde{\mathcal{O}}_1 ( \varepsilon_\omega^2 ),} \\
\\ T_2 (y) &=& \displaystyle{- \frac{\sqrt{2}}{4} \int_{\R} e^{- \sqrt{2} |y-z|} \left ( -2 \frac{g ( \omega Q^2 )}{\omega} +2 \frac{G ( \omega Q^2 )}{\omega^2 Q^2} + Q^2 g' ( \omega Q^2 ) \right ) \, \text{d}z + \tilde{\mathcal{O}}_1 ( \varepsilon_\omega^2 ),} \end{array} \]
and similar expansions for $T_1'$ and $T_2'$. In the lines above, $\tilde{\mathcal{O}}_p ( \varepsilon_\omega^2 )$ denotes any function $\ell$ such that $| \partial_y^k \ell | \leqslant C \varepsilon_\omega^2 (1 + |y|)$ for all $k \in \{ 0 \, , ... \, , p \}$. Let us finally notice that $\varepsilon_\omega = C_\sigma \omega^{\sigma -1}$. We eventually find that $R_1 = a \omega^{\sigma -1} R_1^0 + \tilde{\mathcal{O}}_0 ( \omega^{2 ( \sigma -1)} )$ and $R_2 = a \omega^{\sigma -1} R_2^0 + \tilde{\mathcal{O}}_0 ( \omega^{2 ( \sigma -1)} )$, where
\[ \begin{array}{rl} & R_1^0 = -2QD^0 + (1-Q^2)T_1^0 + (3-Q^2)T_2^0 + \frac{2Q'}{Q} (T_1^0)' + \frac{2Q'}{Q} (T_2^0)' - \frac{2 Q^{2 \sigma}}{\sigma +1} \\
\\ \text{and}& R_2^0 = -4(1-Q^2) QD^0 + 4Q' (D^0)' + T_1^0 -3 T_2^0 + \frac{2Q'}{Q} ((T_1^0)' - (T_2^0)') + \frac{2 ( \sigma -1)}{\sigma +1} Q^{2 \sigma} + \frac{2 Q^{2 \sigma +2}}{\sigma +1}, \end{array} \]
with
\[ \begin{array}{rl} & \displaystyle{T_1^0 (y) = - \frac{( \sigma -1)^2}{2 ( \sigma +1)} \int_{\R} |y-z| Q^{2 \sigma} (z) \, \text{d}z} \\
\\ \text{and} & \displaystyle{T_2^0 (y) = - \frac{\sqrt{2} \, \sigma ( \sigma -1)}{4 ( \sigma + 1)} \int_{\R} e^{- \sqrt{2} |y-z|} Q^{2 \sigma} (z) \, \text{d}z.} \end{array} \]
This leads to $\Delta_1 = a \omega^{\sigma -1} \Delta_1^0 + \tilde{\mathcal{O}}_0 ( \omega^{2 ( \sigma -1)} )$, $\Delta_2 = a \omega^{\sigma -1} \Delta_2^0 + \tilde{\mathcal{O}}_0 ( \omega^{2 ( \sigma -1)} )$ and $\Delta_4 = a \omega^{\sigma-1} \Delta_4^0 + \tilde{\mathcal{O}}_0 ( \omega^{2 ( \sigma -1)} )$, where:
\[ \begin{array}{rl} & \Delta_1^0 = 6Q(1-Q^2) R_1^0 -2QR_2^0 + (3(1-Q^2)^2 -1)D^0 + \sigma \left ( (2 \sigma +1) (1-Q^2)^2 -1 \right ) Q^{2 \sigma -1}, \\
\\ & \Delta_2^0 = QR_1^0 + Q(1-Q^2)R_2^0 + (1-Q^2)D^0 + \sigma (1 - Q^2) Q^{2 \sigma -1} \\
\\ \text{and} & \Delta_4^0 = \Delta_1^0 - ( 2 \sigma + 1) (1-Q^2) Q^{2 \sigma + 1} + 2 (3Q^2 (1-Q^2) +1) D^0. \end{array} \]
This leads to the desired expression:
\[ \Gamma ( \omega ) = a \, \Gamma_0 ( \sigma ) \, \omega^{\sigma-1} + \mathcal{O} \left ( \omega^{2(\sigma -1)} \right ) \]
where $\Gamma_0 ( \sigma ) = \int_{\R} \left ( Q^2 \Delta_4^0 \cos y + \frac{2Q'}{Q} \Delta_5^0 \sin y \right ) \, \text{d}y$ and $\Delta_5^0 := \Delta_4^0 + 2 \Delta_2^0$. Expanding $\Delta_1^0$, $\Delta_2^0$, $\Delta_4^0$ and $\Delta_5^0$, we find the expressions announced in the lemma. \hfill \qedsymbol

\noindent \\ This expression is explicit, and $(H_3)$ will hold if and only if $\Gamma_0 ( \sigma ) > 0$. The curves in Figure 1 show the function $\Gamma_0 ( \sigma )$. They have been obtained with \texttt{python}, and the error is $\simeq 10^{-4}$. For $\sigma = 2$, we find the value $\Gamma_0 (2) = \frac{32 \pi \sqrt{2}}{3 \text{cosh} ( \pi /2 )} \simeq 18.8870$ that has been computed in \cite{Ma0}. 

\begin{figure}
\centering
\includegraphics[width=10cm]{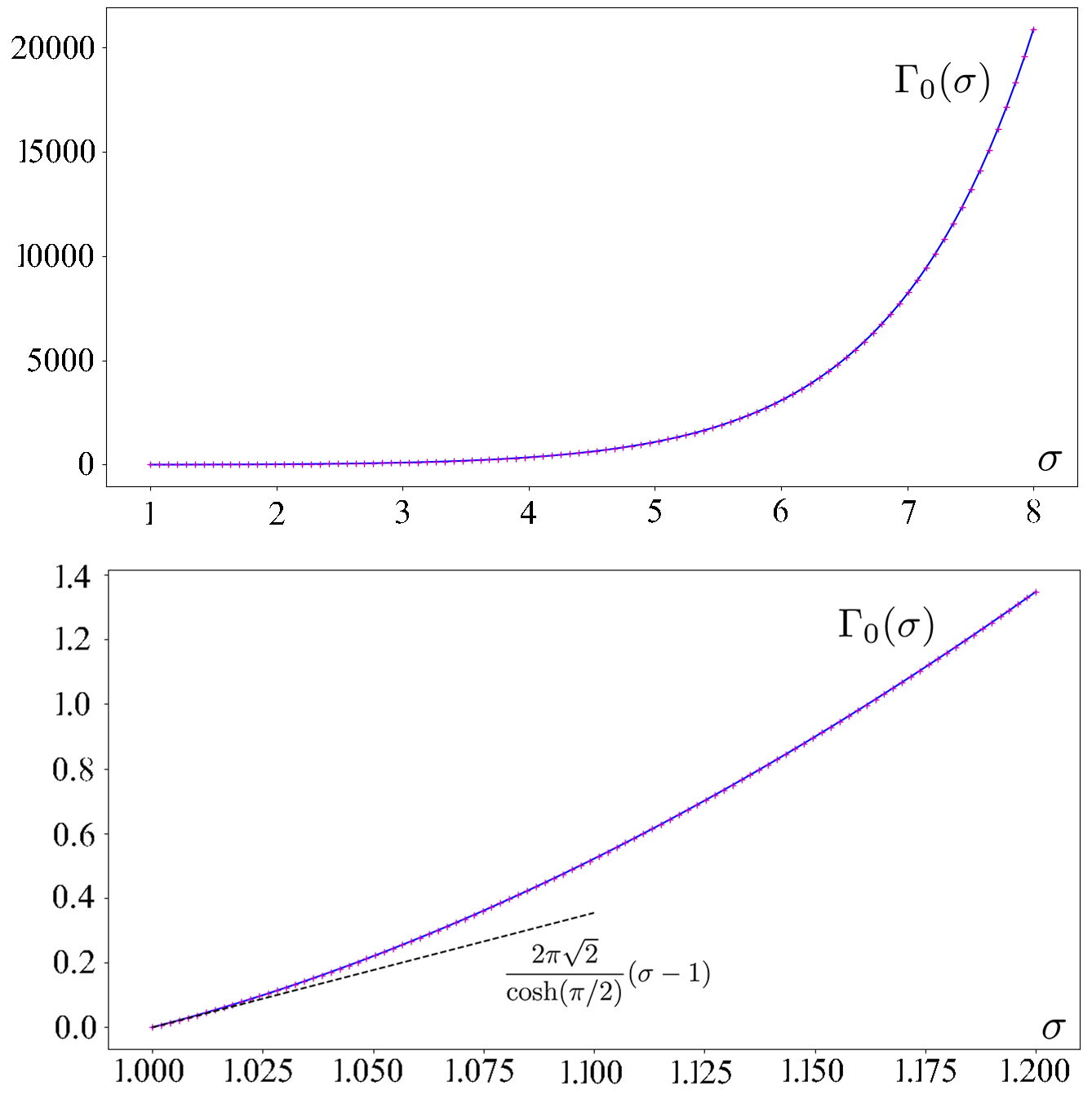}
\caption{Function $\Gamma_0 ( \sigma )$ for $\sigma \in [1 \, , 8]$ (first figure) and $\sigma \in [1 \, , 1.2]$ (second figure).}
\end{figure}

\begin{leftbar}
\noindent \textbf{Numerical check.} For all $\sigma > 1$, we have $\Gamma_0 ( \sigma ) >0$ and thus hypothesis $(H_3)$ holds.
\end{leftbar}

\noindent This cannot be stated as a proposition, as no proof of the positivity of $\Gamma_0 ( \sigma )$ for all $\sigma >1$ will be presented in this paper. However, it is not a conjecture: one can take whatever value of $\sigma >1$ and use the explicit expression from Lemma 16 to check numerically that indeed $\Gamma_0 ( \sigma ) > 0$, and thus hypothesis $(H_3)$ holds. The following lemma gives us the total understanding of $\Gamma_0 ( \sigma )$ for $\sigma \simeq 1^+$. 

\begin{leftbar}
\noindent \textbf{Lemma 17.} We have the following asymptotics:
\[ \Gamma_0 ( \sigma ) \, \underset{\sigma \to 1^+}{\sim} \, \frac{2 \pi \sqrt{2}}{\text{cosh} ( \pi /2)} ( \sigma -1). \]
In particular, $\Gamma_0 ( \sigma ) > 0$ for any $\sigma > 1$ close enough to $1$.
\end{leftbar}

\noindent \textit{Proof.} Until now, the function $\Gamma_0$ has been defined for $\sigma >1$ only, since this is the PDE frame we work in since the beginning of this paper. However, the expression of $\Gamma_0 ( \sigma )$ given in Lemma 16 still holds for $\sigma > 0$ and is continuous, and even $\mathscr{C}^1$, with regards to $\sigma$. We begin by checking that $\Gamma_0 (1)=0$. Let us take $\sigma=1$ for the moment. First, we can integrate explicitly and find that $D^0 = -\frac{Q}{2}$. This leads to $R_1^0 = R_2^0 = 0$, then to
\[ \begin{array}{ll} \Delta_1^0 = Q-3Q^3 + \frac{3Q^5}{2} , & \Delta_2^0 = \frac{Q}{2} - \frac{Q^3}{2} , \\ \\ \Delta_4^0 = 2Q-3Q^3 + \frac{3Q^5}{2} , & \Delta_4^0 + 2 \Delta_2^0 = 3Q - 4Q^3 + \frac{3Q^5}{2}. \end{array} \]
Setting $p_k := \int_{\R} Q^k \cos y \, \text{d}y$, we find that
\[ \Gamma_0 (1) = -6 p_1 + \frac{14}{3} p_3 - \frac{18}{5} p_5 + \frac{3}{2} p_7. \]
Using the relation $p_{k+2} = \frac{2(k^2+1)}{k(k+1)} p_k$, it follows that $\Gamma_0 (1)=0$. \\
\\ Now, let us differentiate $\Gamma_0$ with regards to $\sigma$. Let us start with $D^0$. We have
\[ \partial_\sigma D^0 = -2Q' \int_0^y AQ^{2 \sigma+1} \ln Q + A \left ( \frac{Q^{2 \sigma +2} \ln Q}{\sigma +1} - \frac{Q^{2 \sigma +2}}{2 ( \sigma + 1)^2} \right ). \]
Taking $\sigma=1$, we find that we can integrate explicitly the expression above and, after lengthy computations, we ultimately find that
\[ \left ( \partial_\sigma D^0 \right )_{\sigma =1} = \frac{Q}{4} - Q \ln Q + y Q'. \]
On the other hand, we see that $(\partial_\sigma T_1^0)_{\sigma=1} = 0$ and
\[ ( \partial_\sigma T_2^0 )_{\sigma=1} = - \frac{\sqrt{2}}{8} \int_{\R} e^{- \sqrt{2} |y-z|} Q^2(z) \, \text{d}z =: t_2. \]
Differentiating $\Delta_4^0$ and $\Delta_5^0$ with regards to $\sigma$, lengthy computations lead to
\[ \begin{array}{rcl} ( \partial_\sigma \Delta_4^0 )_{\sigma =1} &=& \frac{37}{4} Q^5 - 17 Q^3 + 4Q + 3 Q^5 \ln Q -6 Q^3 \ln Q + 4 Q \ln Q \\ & & \, \, + 21 y Q' Q^4 - 24 y Q' Q^2 + (6Q^5 - 24Q^3 + 24Q) t_2 + 4 (4-3Q^2) Q' t_2' \end{array} \]
and
\[ \begin{array}{rcl} ( \partial_\sigma \Delta_5^0 )_{\sigma=1} &=& \frac{29}{4} Q^5 - \frac{35}{2} Q^3 + \frac{13}{2} Q + 3Q^5 \ln Q - 8 Q^3 \ln Q + 6 Q \ln Q \\ & & \, \, + 21 y Q' Q^4 - 30 y Q' Q^2 + 2 y Q' + (6Q^5 - 20 Q^3 + 24Q) t_2 + (16-8Q^2) Q' t_2'. \end{array} \]
Setting $q_k := \int_{\R} Q^k \ln Q \cos y \, \text{d}y$, $r_k := \int_{\R} t_2 Q^k \cos y \, \text{d}y$, $s_k := \int_{\R} t_2 Q' Q^{k-1} \sin y \, \text{d}y$ and $m_k := \int_{\R} y Q^k \sin y \, \text{d}y$, integrating by parts we get
\[ \begin{array}{rcl} \int_{\R} Q^2 (\partial_\sigma \Delta_4^0)_{\sigma=1} \cos y \, \text{d}y &=& \frac{37}{4} p_7 - 17 p_5 + 4p_3 + 3q_7 - 6 q_5 + 4q_3 + 24 r_3 - 24 r_5 + 6r_7 + 16(2r_5-3r_3+s_3) \\ & & \, \, - 12 (3r_7-5r_5+s_5) + 21 \left ( - \frac{1}{7} p_7 + \frac{1}{7} m_7 \right ) - 24 \left ( - \frac{1}{5} p_5 + \frac{1}{5} m_5 \right ) \end{array} \]
and
\[ \begin{array}{rcl} \int_{\R} \frac{2Q'}{Q} (\partial_\sigma \Delta_5^0)_{\sigma =1} \sin y \, \text{d}y &=& -\frac{29}{10} p_5 + \frac{35}{3} p_3 - 13p_1 + 6 \left ( \frac{1}{25}p_5 - \frac{1}{5} q_5 \right ) - 16 \left ( \frac{1}{9} p_3 - \frac{1}{3} q_3 \right ) + 12 (p_1 - q_1) \\ & & \, \, + 4m_1 - 62 m_3 + 72m_5 - 21m_7 +12s_5 - 40s_3+48s_1 + 32(-s_1-r_1) \\ & & \, \, + 32(3s_3+r_3) + 8(-5s_5-r_5). \end{array} \]
An integration by parts gives the relation $q_{k+2} = \frac{2(k^2+1)}{k(k+1)} q_k + \frac{2(k^2 -2k-1)}{k^2(k+1)^2} p_k$. From the differential relation $t_2'' = 2t_2 + \frac{Q^2}{2}$ and integrating by parts, we also see that
\[ \begin{array}{rl} & r_{k+2} = \frac{2}{k(k+1)} \left ( (k^2-3)r_k -2ks_k - \frac{p_{k+2}}{2} \right ) \\ \\ \text{and} & s_{k+2} = \frac{2}{(k+1)(k+2)} \left [ (k^2-3)s_k + 2kr_k - (k+1)r_{k+2} + \frac{1}{2(k+2)} p_{k+2} \right ]. \end{array} \]
Integrating by parts one last time, we obtain the relation
\[ m_{k+2} = \frac{2}{k(k+1)} \left [ (k^2+1) m_k -2p_k \right ]. \]
Using these relations, we can express $\Gamma_0 '(1)$ only as a linear combination of $p_1$, $q_1$, $r_1$, $s_1$ and $m_1$. All the occurrences of $q_1$, $r_1$, $s_1$ and $m_1$ conveniently disappear, and we evidently find that
\[ \Gamma_0 '(1) = 2p_1 = \frac{2 \pi \sqrt{2}}{\text{cosh} ( \pi /2)}, \]
which gives us the asymptotic result we desired. \hfill \qedsymbol

\section{Estimate of the internal mode component}

\begin{leftbar}
\noindent \textbf{Proposition 4.} Assume hypothesis $(H_3)$ holds. For any $s > 0$,
\[ \int_0^s |b|^4 \leqslant C \epsilon + \frac{C}{A \alpha ( \omega_0 ) \underline{\mathbf{\Gamma}} ( \omega_0 )} \int_0^s || \rho^4 v ||^2. \]
\end{leftbar}

\noindent \textit{Proof.} Follow the proof of Lemma 21 in \cite{Ma0}. Definitions, estimates and differences are condensed here. We introduce $d_1 = b_1^2 - b_2^2$ and $d_2 = 2b_1b_2$, which verify the equations $\dot{d}_1 = 2 \lambda d_2 + D_2$ and $\dot{d}_2 = -2 \lambda d_1 + D_1$, where $D_2 := 2b_1B_2 + 2b_2B_1$ and $D_1 = 2b_2B_2 - 2b_1B_1$. Let 
\[ \Gamma_1 := \frac{1}{2} \int_{\R} (G^\perp + H^\perp )g_1 = \frac{1}{2} \int_{\R} (G+H)g_1 \, \, \, \, \, \, \, \text{and} \, \, \, \, \, \, \, \Gamma_2 := \frac{1}{4} \int_{\R} (G_1^\perp g_1 - G_2^\perp g_2) = \frac{1}{4} \int_{\R} (G_1 g_1 - G_2 g_2). \]
The equalities above come from the orthogonality relations $\langle V_1 \, , g_2 \rangle = \langle V_2 \, , g_1 \rangle = 0$. We also define
\[ \textbf{J} := d_1 \int_{\R} v_2 g_1 \chi_A - d_2 \int_{\R} v_1 g_2 \chi_A + \Gamma_1 \frac{d_2}{2 \lambda} |b|^2 + \Gamma_2 \frac{d_1 d_2}{2 \lambda}. \]
We compute $\dot{\textbf{J}} = \sum\limits_{j=1}^6 J_j$, where
\[ \begin{array}{rl} & \textbf{J}_1 = d_2 \int_{\R} q_2^\perp g_2 \chi_A + d_1 \int_{\R} q_1^\perp g_1 \chi_A - \Gamma_1 d_1 |b|^2 - \Gamma_2 (d_1^2 - d_2^2), \\
\\ & \textbf{J}_2 = d_2 \int_{\R} v_2 g_2 \chi_A '' + 2 d_2 \int_{\R} ( \partial_y v_2 ) g_2 \chi_A ' + d_1 \int_{\R} v_1 g_1 \chi_A '' + 2 d_1 \int_{\R} ( \partial_y v_1 ) g_1 \chi_A ', \\
\\ & \textbf{J}_3 = - d_2 \int_{\R} ( \mu_2 + p_2^\perp - r_2^\perp ) g_2 \chi_A - d_1 \int_{\R} ( \mu_1 + p_1^{\T} - r_1^{\T} ) g_1 \chi_A, \\
\\ & \textbf{J}_4 = D_2 \int_{\R} v_2 g_1 \chi_A - D_1 \int_{\R} v_1 g_2 \chi_A, \\
\\ & \textbf{J}_5 = \frac{\Gamma_1}{\lambda} \left ( b_1 (b_1 + b_2)^2 B_2 + b_2 ( b_1-b_2)^2 B_1 \right ) + \frac{\Gamma_2}{\lambda} \left ( b_1 |b|^2 B_1 + b_2 (3b_1^2 - b_2^2) B_2 \right ) \\
\\ \text{and} & \textbf{J}_6 = d_1 \int_{\R} v_2 \dot{g}_1 \chi_A - d_2 \int_{\R} v_1 \dot{g}_2 \chi_A + \dot{\Gamma}_1 \frac{d_2}{2 \lambda} |b|^2 + \dot{\Gamma}_2 \frac{d_1 d_2}{2 \lambda} - \dot{\lambda} \Gamma_1 \frac{d_2}{2 \lambda^2} |b|^2 - \dot{\lambda} \Gamma_2 \frac{d_1 d_2}{2 \lambda^2}. \end{array} \]
We have
\[ \begin{array}{rl} & | \textbf{J}_2 | \leqslant \frac{C |b|^2}{\sqrt{A}} \left ( || \eta_A \partial_y v||^2 + \frac{1}{A^2} || \eta_A v||^2 \right )^{1/2} , \\
\\ & | \textbf{J}_3 | \leqslant C \left ( e^{-A} + \epsilon A^{3/2} \right ) |b|^2 ( |b|^2 + || \nu v ||^2 ) , \\
\\ & | \textbf{J}_4 | \leqslant C \alpha ( \omega_0 ) |b| \left ( |b|^2 + || \rho^4 v ||^2 \right ) \sqrt{A} || \eta_A v || \\
\\ \text{and} & | \textbf{J}_5 | \leqslant C \alpha ( \omega_0 ) |b|^3 \left ( |b|^2 + || \rho^4 v ||^2 \right ). \end{array} \]
About $\textbf{J}_1$, we decompose it as: $\textbf{J}_1 = \textbf{J}_{1,1} + \textbf{J}_{1,2} + \textbf{J}_{1,3}$, with
\[ \begin{array}{rl} & \textbf{J}_{1,1} = d_1 \left ( b_1^2 \int_{\R} G^\T g_1 + b_2^2 \int_{\R} H^\T g_1 - \Gamma_1 |b|^2 \right ) + d_2 b_1 b_2 \int_{\R} G_2^\perp g_2 - \Gamma_2 ( d_1^2 - d_2^2 ) = \frac{\Gamma}{4} |b|^4, \\
\\ & \textbf{J}_{1,2} = d_2 \int_{\R} ( q_2^\perp \chi_A - b_1 b_2 G_2^\perp ) g_2 \\
\\ \text{and} & \textbf{J}_{1,3} = d_1 \int_{\R} ( q_1^\T \chi_A - b_1^2 G^\T - b_2^2 H^\T ) g_1. \end{array} \]
To estimate $J_{1,2}$ and $J_{1,3}$, first write that
\[ \begin{array}{rl} & q_1 = b_1^2 G + b_2^2 H + Q_\omega ( 3 + 3 g'(\omega Q_\omega^2) + 2 \omega Q_\omega^2 g''( \omega Q_\omega^2 )) ( 2b_1 V_1 v_1 + v_1^2) + Q_\omega (1 + \omega g'(\omega Q_\omega^2)) (2b_2 V_2 v_2 + v_2^2) + N_1 \\
\\ \text{and} & q_2 = b_1 b_2 G_2 + 2 Q_\omega (1 + g'(\omega Q_\omega^2)) (b_1 V_1 v_2 + b_2 V_2 v_1 + v_1 v_2) + N_2 \end{array} \]
where $|N_2| \leqslant C |u|^3$ but only $|N_1| \leqslant C |u|^{7/3}$. Hence, $|N_1| + |N_2| \leqslant C |v|^{7/3} + C |b|^{7/3} \rho^{16}$, wich leads to
\[ |J_{1,2}| + |J_{1,3}| \leqslant  C \left ( e^{- \frac{\alpha ( \omega_0) A}{2}} + \frac{\epsilon^{1/3}}{\alpha ( \omega_0 )} + \epsilon A^{3/2} \right ) |b|^4 + C (1 + \epsilon A^{3/2} ) |b|^2 || \rho^4 v ||^2 + C \epsilon^{1/3} |b|^2 || \eta_A v ||^2 + C |b|^3 || \nu v ||. \]
About $\textbf{J}_6$, we know from Lemma 14 that $| \partial_\omega g_1 | + | \partial_\omega g_2 | \leqslant C \omega_0^{-1} (1 +  |y|)$. Thus, $| \dot{g}_1 | + | \dot{g}_2 | \leqslant C \omega_0^{-1} | \dot{\omega} | (1 + |y|)$. We get
\[ \left | d_1 \int_{\R} v_2 \dot{g}_1 \chi_A - d_2 \int_{\R} v_1 \dot{g}_2 \chi_A \right | \leqslant C A^{3/2} |b|^2 \left ( || \nu v ||^2 + |b|^2 \right ) || \eta_A v ||. \]
Then we have to estimate $| \partial_\omega \Gamma_1 |$ and $| \partial_\omega \Gamma_2 |$. To do so, we have to estimate $\partial_\omega G_1$, $\partial_\omega G_2$, $\partial_\omega G$ and $\partial_\omega H$. Taking the last one for instance, we compute
\[ \partial_\omega H = 2 V_2 \partial_\omega V_2 Q_\omega (1 + g'( \omega Q_\omega^2 )) + V_2^2 \left ( \partial_\omega Q_\omega ( 1 + g'( \omega Q_\omega^2)) + Q_\omega ( Q_\omega^2 + 2 \omega Q_\omega \partial_\omega Q_\omega ) g'' (\omega Q_\omega^2 ) \right ). \]
We recall that $|V| \leqslant C e^{- \alpha |y|} \leqslant C$, $| \partial_\omega Q_\omega | \leqslant \frac{C}{\omega} (1 + |y|) e^{-|y|}$ and \eqref{dVomega}. This leads to
\[ | \partial_\omega H | \leqslant C \left ( \underline{\textbf{V}} (\omega_0 ) + \omega_0^{-1} \right ) (1 + |y|) e^{-|y|}. \]
We check similarly that
\[ | \partial_\omega G_1 | + | \partial_\omega G_2 | + | \partial_\omega G | \leqslant  C \left ( \underline{\textbf{V}} (\omega_0 ) + \omega_0^{-1} \right ) (1 + |y|) e^{-|y|}. \]
We also know that $|G| + |H| + |G_1 | + |G_2 | \leqslant C e^{-|y|}$. Combining these identities, we get
\[ | \dot{\Gamma}_1 | \leqslant C \left ( \underline{\textbf{V}} (\omega_0 ) + \omega_0^{-1} \right ) \dot{\omega}. \]
Besides, from $\lambda = 1 - \alpha^2$ we get $| \dot{\lambda} | \leqslant \frac{C \varepsilon_{3 \omega_0 /2} \alpha ( \omega_0 )}{\omega_0} | \dot{\omega} |$. Finally, we know from \eqref{m} that $| \dot{\omega} | = \omega | m_\omega | \leqslant C ( || \nu v ||^2 + |b|^2 )$. Gathering these estimates, we find that there exists a quantity $\underline{\textbf{E}} ( \omega_0 )$ depending only on $\omega_0$ such that
\[ \begin{array}{rl} & \left | \dot{\Gamma}_1 \frac{d_2}{\lambda} |b|^2 + \dot{\Gamma}_2 \frac{d_1 d_2}{2 \lambda} - \dot{\lambda} \Gamma_1 \frac{d_2}{2 \lambda^2} |b|^2 - \dot{\lambda} \frac{\Gamma_2}{2 \lambda^2} d_1 d_2 \right | \\
\\ \leqslant & C \underline{\textbf{E}} ( \omega_0 ) A^{3/2} |b|^2 ( || \nu v ||^2 + |b|^2 ) || \eta_A v || + C \underline{\textbf{E}} ( \omega_0 ) ( || \nu v ||^2 + |b|^2 ) |b|^4 \\
\\ \leqslant &  C \underline{\textbf{E}} ( \omega_0 ) A^{3/2} \epsilon ( || \nu v ||^2 + |b|^2 ) || \eta_A v || + C \underline{\textbf{E}} ( \omega_0 ) ( || \nu v ||^2 + |b|^2 ) \epsilon |b|^3. \end{array} \]
Now that all the terms constituting $\dot{\textbf{J}}$ are estimated, we can gather these bounds. Taking $A>0$ large enough (depending on $\omega_0$) and $\epsilon >0$ small enough (depending on $A$ and $\omega_0$), we find that
\[ \begin{array}{rcl} \left | \dot{\textbf{J}} - \frac{\Gamma}{4} |b|^4 \right | = | \dot{\textbf{J}} - J_{1,1} | & \leqslant & C \left ( e^{- \frac{\alpha ( \omega_0) A}{2}} + \frac{\epsilon^{1/3}}{\alpha ( \omega_0 )} + \epsilon A^{3/2} \right ) |b|^4 + C (1 + \epsilon A^{3/2} ) |b|^2 || \rho^4 v ||^2 + C \epsilon^{1/3} |b|^2 || \eta_A v ||^2 \\
\\ & & \, \, \, \, \, + C |b|^3 || \nu v || + \frac{C |b|^2}{\sqrt{A}} \left ( || \eta_A \partial_y v ||^2 + \frac{1}{A^2} || \eta_A v ||^2 \right )^{1/2} + C \sqrt{A} \, \alpha ( \omega_0 ) |b| \left ( |b|^2 + || \rho^4 v ||^2 \right ) || \eta_A v ||. \end{array} \]
We now use hypothesis $(H_3)$: since $\frac{\omega_0}{2} \leqslant \omega \leqslant \frac{3 \omega_0}{2}$, we have $\Gamma ( \omega ) \geqslant \underline{\mathbf{\Gamma}} ( \omega_0 ) > 0$. This leads to
\[ |b|^4 \leqslant \frac{C}{\underline{\mathbf{\Gamma}} ( \omega_0 )} \dot{\textbf{J}} + \frac{C}{A \alpha ( \omega_0 ) \underline{\mathbf{\Gamma}} ( \omega_0 )} \left ( || \eta_A \partial_y v ||^2 + \frac{1}{A^2} || \eta_A v ||^2 \right ). \]
To conclude the proof we integrate the inequality above on $[0 \, , s]$, we recall that $| \textbf{J} | \leqslant C \sqrt{A} \, \epsilon^3$ and we use the first virial result (Proposition 3). Thus we get
\[ \begin{array}{rcl} \displaystyle{\int_0^s |b|^4} & \leqslant & \displaystyle{\frac{C}{\underline{\mathbf{\Gamma}} ( \omega_0 )} \left ( |J(s)| + |J(0)| \right ) + \frac{C}{A \alpha ( \omega_0 ) \underline{\mathbf{\Gamma}} ( \omega_0 )} \int_0^s \left ( || \eta_A \partial_y v ||^2 + \frac{1}{A^2} || \eta_A v ||^2 \right )} \\
\\ & \leqslant & \displaystyle{\frac{C \sqrt{A} \, \epsilon^3}{\underline{\mathbf{\Gamma}} ( \omega_0 )} + \frac{C}{A \alpha ( \omega_0 ) \underline{\mathbf{\Gamma}} ( \omega_0 )} \left ( \epsilon + \int_0^s ( || \rho^4 v ||^2 + |b|^4 ) \right ),} \end{array} \]
from which we deduce that
\[ \left ( 1 - \frac{C}{A \alpha ( \omega_0 ) \underline{\mathbf{\Gamma}} ( \omega_0 )} \right ) \int_0^s |b|^4 \leqslant C \left ( \frac{\sqrt{A} \, \epsilon^2}{\underline{\mathbf{\Gamma}} ( \omega_0 )} + \frac{1}{A \alpha ( \omega_0 ) \underline{\mathbf{\Gamma}} ( \omega_0 )} \right ) \epsilon + \frac{C}{A \alpha ( \omega_0 ) \underline{\mathbf{\Gamma}} ( \omega_0 )} \int_0^s || \rho^4 v ||^2. \]
We now first choose $A>0$ large enough and then $\epsilon>0$ small enough such that $1 - \frac{C}{A \alpha ( \omega_0 ) \underline{\mathbf{\Gamma}} ( \omega_0 )} \geqslant \frac{1}{2}$ and $\frac{\sqrt{A} \, \epsilon^2}{\underline{\mathbf{\Gamma}} ( \omega_0 )} + \frac{1}{A \alpha ( \omega_0 ) \underline{\mathbf{\Gamma}} ( \omega_0 )} \leqslant C$. Consequently, we obtain the desired estimate:
\[ \int_0^s |b|^4 \leqslant C \epsilon + \frac{C}{A \alpha ( \omega_0 ) \underline{\mathbf{\Gamma}} ( \omega_0 )} \int_0^s || \rho^4 v ||^2. \]
\hfill \qedsymbol

\section{The transformed problem}
\noindent For $\theta > 0$ small to be fixed later, we set $X_\theta = (1 - \theta \partial_y^2)^{-1}$. We define $w_1 = X_\theta^2 M_- S^2 v_2$, $w_2 = -X_\theta^2 S^2 L_+ v_1$ and $w=w_1+iw_2$. This constitutes the first transformed problem. The operator $X_\theta$ aims at regularizing the solution after composition by a fourth-order differential operator. Setting $\xi_Q := Q_\omega ' / Q_\omega$, we compute $S^2 = \partial_y^2 -2 \xi_Q \partial_y + 1 + \frac{g( \omega Q_\omega^2 )}{\omega} - 2 \frac{G ( \omega Q_\omega^2 )}{\omega^2 Q_\omega^2}$,
\[ \begin{array}{rcl} M_- S^2 &=& - \partial_y^4 + 2 \partial_y^2 \cdot \xi_Q \cdot \partial_y + \partial_y \cdot \left ( -2 Q_\omega^2 g'(\omega Q_\omega^2) + 4 \frac{g ( \omega Q_\omega^2 )}{\omega} - 4 \, \frac{G ( \omega Q_\omega^2)}{\omega^2 Q_\omega^2} \right ) \cdot \partial_y \\
\\ & & \, \, + \left ( - 4 Q_\omega Q_\omega ' g' ( \omega Q_\omega^2 ) + 6 \, \xi_Q  \frac{g( \omega Q_\omega^2 )}{\omega} + 4 \omega Q_\omega ' Q_\omega^3 g''( \omega Q_\omega^2 ) - 4 \xi_Q \frac{G ( \omega Q_\omega^2 )}{\omega^2 Q_\omega^2} - 2 \xi_Q \right ) \cdot \partial_y \\
\\ & & \, \, + \, 1 + 2 \left ( - \frac{g ( \omega Q_\omega^2 )}{\omega} + Q_\omega^2 g'(\omega Q_\omega^2) - 2 \omega Q_\omega^4 g''( \omega Q_\omega^2) \right ) \\
\\ & & \, \, - \, 2 \frac{g'(\omega Q_\omega^2) G ( \omega Q_\omega^2 )}{\omega^2} - Q_\omega^4 g'(\omega Q_\omega^2) + 2 \omega Q_\omega^6 g''( \omega Q_\omega^2) + 4 Q_\omega^2 \frac{G(\omega Q_\omega^2) g''(\omega Q_\omega^2)}{\omega^2} \\
\\ & & \, \, + \, 2 Q_\omega^2 \frac{g ( \omega Q_\omega^2 )}{\omega} - 2 \frac{G ( \omega Q_\omega^2 )}{\omega^2} + \frac{g ( \omega Q_\omega^2 )^2}{\omega^2} \end{array} \]
and
\[ \begin{array}{rcl} S^2 L_+ &=& - \partial_y^4 + 2 \partial_y^2 \cdot \xi_Q \cdot \partial_y + \partial_y \cdot \left ( - Q_\omega^2 + 2 \frac{g ( \omega Q_\omega^2 )}{\omega} - 2 \, \frac{G ( \omega Q_\omega^2 )}{\omega^2 Q_\omega^2} - 2 Q_\omega^2 g'( \omega Q_\omega^2 ) \right ) \cdot \partial_y \\
\\ & & \, \, + \left ( -2 Q_\omega Q_\omega ' - 4 Q_\omega Q_\omega ' g' ( \omega Q_\omega^2 ) + 2 \xi_Q \frac{g( \omega Q_\omega^2 )}{\omega} - 4 \omega Q_\omega ' Q_\omega^3 g''( \omega Q_\omega^2 ) -2 \xi_Q \right ) \cdot \partial_y \\
\\ & & \, \, + \, 1 + \left ( -3 Q_\omega^2 - 20 \omega Q_\omega^4 g''(\omega Q_\omega^2) - 8 \omega^2 Q_\omega^6 g'''(\omega Q_\omega^2) - 2 Q_\omega^2 g'(\omega Q_\omega^2) - 2 \, \frac{G ( \omega Q_\omega^2)}{\omega^2 Q_\omega^2} \right ) \\
\\ & & \, \, + \, 3 Q_\omega^4 + 3 Q_\omega^2 \frac{g ( \omega Q_\omega^2)}{\omega} + 3 Q_\omega^4 g'(\omega Q_\omega^2) + 4 Q_\omega^2 \frac{g ( \omega Q_\omega^2 ) g'( \omega Q_\omega^2 )}{\omega} -2 \frac{g'( \omega Q_\omega^2 ) G( \omega Q_\omega^2 )}{\omega^2} \\
\\ & & \, \, + \, 12 \omega Q_\omega^6 g''(\omega Q_\omega^2) + 16 Q_\omega^2 \frac{G( \omega Q_\omega^2 ) g''(\omega Q_\omega^2)}{\omega} + 4 Q_\omega^4 g( \omega Q_\omega^2 ) g''( \omega Q_\omega^2 ) - 4  \omega^2 Q_\omega^8 g'''(\omega Q_\omega^2) \\
\\ & & \, \, + \, 8 Q_\omega^4 G ( \omega Q_\omega^2 ) g''' ( \omega Q_\omega^2 ) - \frac{g ( \omega Q_\omega^2 )^2}{\omega^2} + 2 \frac{g ( \omega Q_\omega^2 ) G ( \omega Q_\omega^2 )}{\omega^3 Q_\omega^2}. \end{array} \]
We introduce the operators $Q_-$ and $Q_+$, obtained respectively from $M_- S^2$ and $S^2 L_+$ by differentiation with respect to $\omega$ and then multiplication by $\omega$. Their exact expressions are given below.
\[ \begin{array}{rcl} Q_- &=& 2 \partial_y^2 \cdot \partial_\omega \xi_Q \cdot \partial_y + \partial_y \cdot \partial_\omega \left ( -2 Q_\omega^2 g'(\omega Q_\omega^2) + 4 \frac{g ( \omega Q_\omega^2 )}{\omega} - 4 \, \frac{G ( \omega Q_\omega^2)}{\omega^2 Q_\omega^2} \right ) \cdot \partial_y \\
\\ & & \, \, + \partial_\omega \left ( - 4 Q_\omega Q_\omega ' g' ( \omega Q_\omega^2 ) + 6 \, \xi_Q  \frac{g( \omega Q_\omega^2 )}{\omega} + 4 \omega Q_\omega ' Q_\omega^3 g''( \omega Q_\omega^2 ) - 4 \xi_Q \frac{G ( \omega Q_\omega^2 )}{\omega^2 Q_\omega^2} - 2 \xi_Q \right ) \cdot \partial_y \\
\\ & & \, \, + \partial_\omega \left [ 2 \left ( - \frac{g ( \omega Q_\omega^2 )}{\omega} + Q_\omega^2 g'(\omega Q_\omega^2) - 2 \omega Q_\omega^4 g''( \omega Q_\omega^2) \right ) \right. \\
\\ & & \, \, - \, 2 \frac{g'(\omega Q_\omega^2) G ( \omega Q_\omega^2 )}{\omega^2} - Q_\omega^4 g'(\omega Q_\omega^2) + 2 \omega Q_\omega^6 g''( \omega Q_\omega^2) + 4 Q_\omega^2 \frac{G(\omega Q_\omega^2) g''(\omega Q_\omega^2)}{\omega^2} \\
\\ & & \left. \, \, + \, 2 Q_\omega^2 \frac{g ( \omega Q_\omega^2 )}{\omega} - 2 \frac{G ( \omega Q_\omega^2 )}{\omega^2} + \frac{g ( \omega Q_\omega^2 )^2}{\omega^2} \right ] \end{array} \]
and
\[ \begin{array}{rcl} Q_+ &=& 2 \partial_y^2 \cdot \partial_\omega \xi_Q \cdot \partial_y + \partial_y \cdot \partial_\omega \left ( - Q_\omega^2 + 2 \frac{g ( \omega Q_\omega^2 )}{\omega} - 2 \, \frac{G ( \omega Q_\omega^2 )}{\omega^2 Q_\omega^2} - 2 Q_\omega^2 g'( \omega Q_\omega^2 ) \right ) \cdot \partial_y \\
\\ & & \, \, + \partial_\omega \left ( -2 Q_\omega Q_\omega ' - 4 Q_\omega Q_\omega ' g' ( \omega Q_\omega^2 ) + 2 \xi_Q \frac{g( \omega Q_\omega^2 )}{\omega} - 4 \omega Q_\omega ' Q_\omega^3 g''( \omega Q_\omega^2 ) -2 \xi_Q \right ) \cdot \partial_y \\
\\ & & \, \, + \partial_\omega \left [ \left ( -3 Q_\omega^2 - 20 \omega Q_\omega^4 g''(\omega Q_\omega^2) - 8 \omega^2 Q_\omega^6 g'''(\omega Q_\omega^2) - 2 Q_\omega^2 g'(\omega Q_\omega^2) - 2 \, \frac{G ( \omega Q_\omega^2)}{\omega^2 Q_\omega^2} \right ) \right. \\
\\ & & \, \, + \, 3 Q_\omega^4 + 3 Q_\omega^2 \frac{g ( \omega Q_\omega^2)}{\omega} + 3 Q_\omega^4 g'(\omega Q_\omega^2) + 4 Q_\omega^2 \frac{g ( \omega Q_\omega^2 ) g'( \omega Q_\omega^2 )}{\omega} -2 \frac{g'( \omega Q_\omega^2 ) G( \omega Q_\omega^2 )}{\omega^2} \\
\\ & & \, \, + \, 12 \omega Q_\omega^6 g''(\omega Q_\omega^2) + 16 Q_\omega^2 \frac{G( \omega Q_\omega^2 ) g''(\omega Q_\omega^2)}{\omega} + 4 Q_\omega^4 g( \omega Q_\omega^2 ) g''( \omega Q_\omega^2 ) - 4  \omega^2 Q_\omega^8 g'''(\omega Q_\omega^2) \\
\\ & & \left. \, \, + \, 8 Q_\omega^4 G ( \omega Q_\omega^2 ) g''' ( \omega Q_\omega^2 ) - \frac{g ( \omega Q_\omega^2 )^2}{\omega^2} + 2 \frac{g ( \omega Q_\omega^2 ) G ( \omega Q_\omega^2 )}{\omega^3 Q_\omega^2} \right ]. \end{array} \]
The fully developed versions of the operators $Q_-$ and $Q_+$ (not rescaled) can be found at the beginning of section 3.3 in \cite{Ri}. \\
\\ From \eqref{sysv} and the identity $S^2 L_+ L_- = M_+ M_- S^2$, the function $w$ satisfies the system
\begin{equation}
    \left \{ \begin{array}{ccl} \dot{w}_1 &=& M_- w_2 + \left [ X_\theta^2 \, , \frac{a_\omega^-}{\omega} \right ] S^2 L_+ v_1 + X_\theta^2 n_2 \\
    \\ \dot{w}_2 &=& -M_+w_1 - \left [ X_\theta^2 \, , \frac{a_\omega^+}{\omega} \right ] M_- S^2 v_2 - X_\theta^2 n_1 \end{array} \right.
    \label{sysw}
\end{equation}
where $[X_\theta^2 \, , a ] = X_\theta^2 a - a X_\theta^2$, $n_1 = -S^2 L_+ p_2^\perp + S^2 L_+ q_2^\perp + S^2 L_+ r_2^\perp + \dot{\omega} Q_+ v_1$ and $n_2 = -M_- S^2 p_1^\T + M_- S^2 q_1^\T + M_- S^2 r_1^\T + \dot{\omega} Q_- v_2$. \\
\\ Now we set the second transformed problem, whose goal is to suppress the internal mode: for $\vartheta > \theta$ small to be chosen (later we will eventually choose $\vartheta = \theta^{1/4}$), we define $z_1 = X_\vartheta U w_2$, $z_2 = - X_\vartheta U M_+ w_1$ and $z=z_1 + iz_2$. We denote $\xi_W := W_2' / W_2$, which implies that $U = \partial_y - \xi_W$ and
\[ UM_+ = - \partial_y^3 + \partial_y \cdot \xi_W \cdot \partial_y + \partial_y - \xi_W ' \partial_y + a_\omega^+ \partial_y - \xi_W - \xi_W a_\omega^+ + (a_\omega^+)'. \]
From \eqref{sysw} and the identity $UM_+ M_- = KU$ (see Lemma 2), the function $z$ satisfies the system
\begin{equation}
    \left \{ \begin{array}{ccl} \dot{z}_1 &=& z_2 - X_\vartheta U \left [ X_\theta^2 \, , \frac{a_\omega^+}{\omega} \right ] M_- S^2 v_2 - X_\vartheta U X_\theta^2 n_1 + \dot{\omega} X_\vartheta P_+ w_2 \\
    \\ \dot{z}_2 &=& -Kz_1 - [X_\vartheta \, , K] Uw_2 - X_\vartheta UM_+ \left [ X_\theta^2 \, , \frac{a_\omega^-}{\omega} \right ] S^2 L_+ v_1 - X_\vartheta U M_+ X_\theta^2 n_2 - \dot{\omega} X_\vartheta P_- w_1 \end{array} \right.
    \label{sysz}
\end{equation}
where $P_+ = - \partial_\omega \xi_W$ and
\[ P_- = \partial_y \cdot \partial_\omega \xi_W \cdot \partial_y - ( \partial_\omega \xi_W ' ) \partial_y + ( \partial_\omega a_\omega^+ ) \partial_y + \partial_\omega \left ( - \xi_W - a_\omega^+ \xi_W + (a_\omega^+)' \right ). \]
Before going further, we will need the following technical lemma in order to estimate $\xi_W$ and its derivatives.

\begin{leftbar}
\noindent \textbf{Lemma 18.} We have the following bounds on $\xi_W$:
\begin{itemize}
    \item for any $k \in \N$, $| \partial_y^k \xi_W | \leqslant C \varepsilon_\omega$ on $\R$;
    \item for any $k \in \{ 0 \, , ... \, , 3 \}$, there exists a quantity $\underline{\mathbf{\zeta}}_{k} ( \omega_0 )$ depending only on $k$ and $\omega_0$ such that $| \partial_y^k \partial_\omega \xi_W | \leqslant C \underline{\mathbf{\zeta}}_k ( \omega_0 )$ on $\R$. 
\end{itemize}
\end{leftbar}

\noindent \textit{Proof.} The first point is obtained easily thanks to the estimates $| W_2^{(k)} (y)| \leqslant C \varepsilon_\omega e^{- \alpha |y|}$ and $W_2 (y) \geqslant \frac{1}{2} e^{- \alpha |y|}$. For the second point, we take $y>0$ and recall the following identity established in the proof of Lemma 3:
\[ \xi_W (y) = - \sqrt{\alpha^2 + \frac{2}{W_2^2 (y)} \int_y^{+ \infty} w_0 W_2' W_2}, \]
where we recall that $w_0 = \lambda \frac{W_1 - W_2}{W_2} - a_\omega^-$. Thus,
\[ \partial_\omega \xi_W = -\frac{1}{2 \xi_W} \left ( 2 \alpha \alpha '( \omega) - \frac{4 \partial_\omega W_2}{W_2^3} \int_y^{+ \infty} w_0 W_2' W_2 + \frac{2}{W_2^2} \int_y^{+ \infty} ( \partial_\omega w_0 W_2' W_2 + w_0 \partial_\omega W_2 ' W_2 + w_0 W_2 ' \partial_\omega W_2 ) \right ). \]
We recall the following estimates from Proposition 2 and the proof of Lemma 3: $| \alpha ' (\omega) | \leqslant \frac{C \varepsilon_\omega}{\omega}$, $| \partial_\omega W_j | \leqslant \frac{C \varepsilon_\omega \varrho_\omega}{\omega \alpha} (1 + |y|) e^{- \alpha |y|}$, $| \partial_\omega W_2' | \leqslant \frac{C \varepsilon_\omega}{\omega} (1+|y|) e^{- \alpha |y|}$, $|W_2| \leqslant Ce^{- \alpha |y|}$, $|W_2'| \leqslant C \varepsilon_\omega e^{-\alpha |y|}$, $|w_0| \leqslant C \varepsilon_\omega e^{-|y|}$, $| \partial_\omega \lambda | \leqslant \frac{C \varepsilon_\omega \alpha}{\omega}$, $| \partial_\omega a_\omega^- | \leqslant \frac{C \varepsilon_\omega}{\omega} (1 + |y|) e^{-2|y|}$, $| W_1 - W_2 | \leqslant C \varepsilon_\omega e^{- \kappa |y|}$ and $| \partial_\omega (W_1 - W_2) | \leqslant \frac{C \varepsilon_\omega}{\omega} (1 + |y|) e^{- \kappa |y|}$. For this last one, one has to check the proof of Proposition 2 and recall that $W_1 - W_2 = 2X_2$. Gathering all these estimates, we find that
\[ | \partial_\omega \xi_W | \leqslant \frac{C \varrho_\omega^2}{\omega | \xi_W |}. \]
Now, we also find, thanks to the same estimates, that $\left | \frac{2}{W_2^2} \int_y^{+ \infty} w_0 W_2' W_2 \right | \leqslant C \varepsilon_\omega^2 e^{-|y|}$. Thus, for $y \geqslant y_\omega^1 := \ln \left ( \frac{2C \varepsilon_\omega^2}{\alpha^2} \right )$, we have $\left | \frac{2}{W_2^2} \int_y^{+ \infty} w_0 W_2' W_2 \right | \leqslant \frac{\alpha^2}{2}$ and thus $| \xi_W | \geqslant \alpha^2 / 2$. For such $y$, we have
\[ | \partial_\omega \xi_W | \leqslant \frac{C \varrho_\omega^2}{\omega \alpha^2} \leqslant \frac{C \varepsilon_\omega^4}{\omega \alpha^4} \leqslant \frac{C \varepsilon_{3 \omega / 2}^4}{\omega_0 \alpha ( \omega_0 )^4}. \]
Now, take $0 < y < y_\omega^1$. Recalling that $| \partial_\omega W_2 | \leqslant \frac{C \varepsilon_\omega \varrho_\omega}{\omega \alpha} (1 + |y|) e^{- \alpha |y|}$, $| \partial_\omega W_2' | \leqslant \frac{C \varepsilon_\omega}{\omega} (1+|y|) e^{- \alpha |y|}$, $|W_2| \leqslant Ce^{- \alpha |y|}$, $|W_2'| \leqslant C \varepsilon_\omega e^{-\alpha |y|}$, $|w_0| \leqslant C \varepsilon_\omega e^{-|y|}$ and $W_2 \geqslant \frac{1}{2} e^{- \alpha |y|}$, an elementary calculation of $\partial_\omega \xi_W$ shows that
\[ | \partial_\omega \xi_W | \leqslant \frac{C \varrho_\omega}{\omega} (1 + |y|) \leqslant \frac{C \varrho_\omega}{\omega} (1 + y_\omega^1) \leqslant \frac{C \varepsilon_{3 \omega_0/2}^2}{\omega_0 \alpha ( \omega_0 )} \left ( 1 + \ln \left ( \frac{2C \varepsilon_{3 \omega_0 / 2}^2}{\alpha ( \omega_0 )^2} \right ) \right ). \]
Similar considerations hold for $y<0$. Setting $\underline{\mathbf{\zeta}}_0 ( \omega_0 ) := \max \left ( \frac{C \varepsilon_{3 \omega / 2}^4}{\omega_0 \alpha ( \omega_0 )^4} \, , \frac{C \varepsilon_{3 \omega_0/2}^2}{\omega_0 \alpha ( \omega_0 )} \left ( 1 + \ln \left ( \frac{2C \varepsilon_{3 \omega_0 / 2}^2}{\alpha ( \omega_0 )^2} \right ) \right ) \right )$, we get the desired result for $k=0$. The result for larger values of $k$ is obtained similarly. It does not matter, for later proofs, that the quantities $\underline{\mathbf{\zeta}}_k ( \omega_0 )$ do not vanish as $\omega_0 \to 0$. \hfill \qedsymbol

\noindent \\ Now we follow \cite{Ma0} (see Lemmas 22 to 27) to give useful technical lemmas about the operators $X_\theta$. The proofs are globally unchanged. 

\begin{leftbar}
\noindent \textbf{Lemma 19.} For $\theta >0$ small enough and all $h \in L^2 ( \R )$,
\[ \begin{array}{lll} || X_\theta h || \leqslant C ||h|| , & || \partial_y X_\theta^{1/2} h || \leqslant C \theta^{-1/2} ||h|| , & || \rho X_\theta h || \leqslant C || X_\theta ( \rho h ) || , \\
\\ || \eta_A^{-1} X_\theta ( \eta_A h) || \leqslant C || X_\theta h || , & || \eta_A X_\theta h || \leqslant C || X_\theta ( \eta_A h )||, & || \eta_A X_\theta \partial_y h || \leqslant C \theta^{-1/2} || \eta_A h || , \\
\\ || \eta_A X_\theta \partial_y^2 h || \leqslant C \theta^{-1} || \eta_A h || , & ||\rho^{-1} X_\theta ( \rho h ) || \leqslant C || X_\theta h || , & || \rho^{-1} X_\theta \partial_y ( \rho h ) || \leqslant C \theta^{-1/2} ||h|| , \\
\\ || \rho^{-1} X_\theta \partial_y^2 ( \rho h ) || \leqslant C \theta^{-1} ||h||. \end{array} \]
\end{leftbar}

\begin{leftbar}
\noindent \textbf{Lemma 20.} For $\theta >0$ small enough and all $h \in H^1 ( \R )$,
\[ \begin{array}{l} ||\eta_A X_\theta^2 M_- S^2 h || + || \eta_A X_\theta^2 S^2 L_+ h || \leqslant C \theta^{-2} || \eta_A h ||, \\
\\ ||\eta_A X_\theta^2 M_- S^2 h || + || \eta_A X_\theta^2 S^2 L_+ h || \leqslant C \theta^{-3/2} || \eta_A \partial_y h || + C || \eta_A h || \\
\\ || \eta_A \partial_y X_\theta^2 M_- S^2 h || + || \eta_A \partial_y X_\theta^2 S^2 L_+ h || \leqslant C \theta^{-2} || \eta_A \partial_y h || + C || \eta_A h || , \\
\\ || \eta_A \partial_y^2 X_\theta Uh || + || \eta_A \partial_y X_\theta Uh || + || \eta_A X_\theta U h || \leqslant C \theta^{-1} || \eta_A \partial_y h || + C ||\eta_A h || , \\
\\ ||\eta_A X_\theta M_+ h || \leqslant C \theta^{-1} || \eta_A h || , \\
\\ || \eta_A X_\theta U M_+ h || \leqslant C \theta^{-1} || \eta_A \partial_y h || + || \eta_A h ||. \end{array} \]
\end{leftbar}

\noindent \textit{Proof.} See \cite{Ma0}: the first three points are analogous to Lemma 23; the last three points are analogous to Lemma 24. The proof is identical and requires the bound $| \xi_W | \leqslant C$, which is proven in our Lemma 18 above. \hfill \qedsymbol

\noindent \\ Applying the estimates above to the definitions of $v$ and $w$, we find the following result.

\begin{leftbar}
\noindent \textbf{Lemma 21.} For $0 < \theta < \vartheta^2$ small enough, and for all $s \geqslant 0$,
\[ \begin{array}{l} || \eta_A \partial_y w || + || \eta_A w || \leqslant C \theta^{-2} || \eta_A \partial_y v || + C|| \eta_A v ||, \\
\\ || \eta_A \partial_y^2 z_1|| + || \eta_A \partial_y z_1|| + || \eta_A z_1 || \leqslant C \vartheta^{-1} || \eta_A \partial_y w_2 || + C || \eta_A w_2 || , \\
\\ || \eta_A z_2 || \leqslant C \vartheta^{-1} || \eta_A \partial_y w_1 || + || \eta_A w_1 ||. \end{array} \]
\end{leftbar}

\noindent In \cite{Ri} (see Lemma 11) one can find the proof of the following lemma (it is the same result, here rescaled).

\begin{leftbar}
\noindent \textbf{Lemma 22.} For $\theta >0$ small engouh and any $h \in H^1 ( \R )$,
\[ ||\eta_A X_\theta^2 Q_- h || + || \eta_A X_\theta^2 Q_+ h || \leqslant C \theta^{-1} || \eta_A \partial_y h || + C || \eta_A h ||. \]
\end{leftbar}

\noindent The last technical lemma is the following.

\begin{leftbar}
\noindent \textbf{Lemma 23.} There exists a quantity $\underline{\mathbf{P}} ( \omega_0 )$ depending only on $\omega_0$ such that, for $\theta >0$ small enough and any $h \in H^1 ( \R )$,
\[ \begin{array}{rl} & ||\eta_A X_\theta P_- h || \leqslant C \underline{\mathbf{P}} ( \omega_0 ) \left ( \theta^{-1/2} || \eta_A \partial_y h || + || \eta_A h || \right ) \\
\\ \text{and}& ||\eta_A P_+ h ||\leqslant C \underline{\mathbf{P}} ( \omega_0 ) ||\eta_A h ||. \end{array} \]
\end{leftbar}

\noindent \textit{Proof.} See the proof of Lemma 27 in \cite{Ma0}: the difference comes from the fact that, here, we do not have $| \partial_y^k \partial_\omega \xi_W | \leqslant C$ but simply $| \partial_y^k \partial_\omega \xi_W | \leqslant C \underline{\mathbf{\zeta}}_k ( \omega_0 )$. This implies the presence of the factor $\underline{\mathbf{P}} ( \omega_0 )$ in the estimates above. We will ultimately find that this factor, depending only on $\omega_0$, does not hinder the proofs to come. \hfill \qedsymbol

\begin{leftbar}
\noindent \textbf{Lemma 24.} Let $\widetilde{z} := \chi_A \zeta_B z$. For all $s \geqslant 0$,
\[ || \rho \partial_y^2 z_1 || + || \rho \partial_y z_1 || + || \rho z_1 || \leqslant C \left ( || \partial_y^2 \widetilde{z}_1 || + || \partial_y \widetilde{z}_1 || + || \rho^{1/2} \widetilde{z}_1 || + A^{-2} \theta^{-5/2} ( ||\eta_A \partial_y v || + || \eta_A v || ) \right ). \]
\end{leftbar}

\noindent \textit{Proof.} See the proof of Lemma 28 in \cite{Ma0}. \hfill \qedsymbol

\section{Coercivity of the transformed problem}
\noindent The goal of this section is to control $w$ thanks to $z$, in other words to go back from $z$ to $w$. 

\begin{leftbar}
\noindent \textbf{Lemma 25.} For all $s \geqslant 0$,
\[ \begin{array}{rl} & || \rho^2 \partial_y w_2 || + || \rho^2 w_2 || \leqslant C \left ( \vartheta ||\rho \partial_y^2 z_1 || + \vartheta || \rho \partial_y z_1 || + \alpha ( \omega_0 )^{-1} || \rho z_1 || \right ) \\
\\ \text{and} & || \rho^2 \partial_y w_1 || + || \rho^2 w_1 || \leqslant C \alpha ( \omega_0 )^{-3/2} || \rho z_2 ||. \end{array} \]
\end{leftbar}

\noindent \textit{Proof.} Follow the proof of Lemma 29 in \cite{Ma0}. We first begin by checking that $| \langle w_1 \, , W_2 \rangle | \leqslant C \theta \varepsilon_{3 \omega_0 / 2} || \rho^2 w_1 ||$ and $| \langle w_2 \, , W_1 \rangle | \leqslant C \theta \varepsilon_{3 \omega_0 /2} || \rho^2 w_2 ||$. We show that
\begin{equation}
    w_2 = a W_2 - \vartheta \partial_y z_1 - \vartheta \frac{W_2'}{W_2} z_1 + W_2 \int_0^y \frac{m_2 z_1}{W_2}
    \label{lem21}
\end{equation}
with $| \langle z_1 \, , W_1' \rangle | \leqslant C \sqrt{\alpha ( \omega_0 )} || \rho z_1 ||$, $\left | \left \langle z_1 \, , \frac{W_2' W_1}{W_2} \right \rangle \right | \leqslant C \sqrt{\alpha ( \omega_0 )} ||\rho z_1||$, 
\[ \left | \int_0^y \frac{m_2 z_1}{W_2} \right | \leqslant \frac{C || \rho z_1 ||}{\sqrt{\alpha ( \omega_0 )}} \rho^{-1} e^{\alpha |y|} \, \, \, \, \, \, \, \text{and} \, \, \, \, \, \, \, \left | \left \langle W_2 \int_0^y \frac{z_1 m_2}{W_2} \, , W_1 \right \rangle \right | \leqslant \frac{C ||\rho z_1||}{\alpha ( \omega_0 )^{3/2}}. \]
This leads to the estimate $|a| \leqslant C \left ( \theta \alpha ( \omega_0 ) \varepsilon_{3 \omega_0 / 2} || \rho^2 w_2|| + \alpha ( \omega_0 )^{-1/2} || \rho z_1 || \right )$. Then we multiply \eqref{lem21} by $\rho^2$ and we find
\[ (1 - C \theta \sqrt{\alpha ( \omega_0 )} \varepsilon_{3 \omega_0 /2} ) || \rho^2 w_2 || \leqslant \frac{C}{\alpha ( \omega_0 )} ||\rho z_1|| + C \vartheta || \rho^2 \partial_y z_1 ||, \]
which gives the result by taking $\theta >0$ small enough (depending on $\omega_0$). We differentiate \eqref{lem21} with regards to $y$ in order to get the similar estimate for $\partial_y w_2$. The proof for $w_1$ and $\partial_y w_1$ is similar but requires the introduction of $H_1$ and $H_2$, solutions to $M_+ H=0$ that satisfy $H_1'H_2-H_1H_2'=1$, $|H_1^{(k)}(y)| \leqslant Ce^{-y}$ and $|H_2^{(k)}(y)| \leqslant Ce^y$ on $\R$. The existence of these two functions is established in \cite{Ri} (see Lemma 3). The rest of the proof does not present any complication. \hfill \qedsymbol

\begin{leftbar}
\noindent \textbf{Lemma 26.} For all $s \geqslant 0$,
\[ \begin{array}{rl} & || \rho^4 v_1 || \leqslant C || \rho^2 w_2 || \leqslant C \left ( \vartheta ||\rho \partial_y^2 z_1|| + \vartheta || \partial_y z_1 || + \alpha ( \omega_0 )^{-1} || \rho z_1 || \right ) \\
\\ \text{and} & || \rho^4 v_2 || \leqslant C || \rho^2 w_1 || \leqslant C \alpha ( \omega_0 )^{-3/2} ||\rho z_2 ||. \end{array} \]
\end{leftbar}

\noindent \textit{Proof.} The analogous result in \cite{Ma0} (Lemma 30) is established by adapting the proof of Proposition 19 in \cite{Ma1}. We follow the same idea, adapting instead the proof of Proposition 5 in \cite{Ri}. It does not present any complication. \hfill \qedsymbol

\section{Estimate on the transformed problem}
\noindent We here give the last virial argument that we will use, the one concerning the transformed problem \eqref{sysz}. It relies on the repulsive nature of the potential of the operator $K$. 

\begin{leftbar}
\noindent \textbf{Proposition 5.} Assume hypotheses $(H_1)$, $(H_2)$ and $(H_3)$ hold. For all $s \geqslant 0$,
\[ \int_0^s \left ( || \rho \partial_y^2 z_1 ||^2 + || \rho \partial_y z_1 ||^2 + || \rho z_1 ||^2 + || \rho z_2 ||^2 \right ) \leqslant C \sqrt{\epsilon} + \frac{C}{\sqrt{A}} \int_0^s ||\rho^4 v||^2. \]
\end{leftbar}

\noindent \textit{Proof.} Follow the proof of Lemma 31 in \cite{Ma0}. Definitions, estimates and differences are condensed here. The parameters will be chosen in the following order, in order to complete the proof: first $\omega_0 >0$ small enough, then $B>0$ large enough (depending on $\omega_0$), then $\theta>0$ small enough (depending on $\omega_0$ and $B$), then $\vartheta = \theta^{1/4} >0$, then $A>0$ large enough (depending on all the previous parameters), and finally $\epsilon >0$ small enough (depending on all the previous parameters). In short:
\[ \omega_0 \longrightarrow B  \longrightarrow \theta \longrightarrow \vartheta \longrightarrow A \longrightarrow \epsilon. \]
Now, let
\[ \textbf{K} := - \int_{\R} ( \Xi_{A,B} z_1 ) z_2 \, \, \, \, \, \, \, \, \, \text{and} \, \, \, \, \, \, \, \, \, \textbf{L} := \int_{\R} \rho^2 z_1 z_2. \]
We have $| \textbf{K} | + | \textbf{L} | \leqslant C \epsilon$, by taking $\epsilon >0$ small enough (depending on $B$, $\theta$ and $\vartheta$). Then we compute: $\dot{\textbf{K}} = \sum\limits_{j=1}^5 \textbf{K}_j$, $\dot{\textbf{L}} = \sum\limits_{j=1}^5 \textbf{L}_j$ and $\textbf{K}_1 = \textbf{P} + \sum\limits_{j=1}^9 \textbf{R}_j$, where
\[ \begin{array}{ccl} \textbf{K}_1 &=& \int_{\R} ( \Xi_{A,B} z_1) Kz_1, \\
\\ \textbf{K}_2 &=& \int_{\R} ( \Xi_{A,B} z_1) [ X_\vartheta \, , K] U w_2 , \\
\\ \textbf{K}_3 &=& -\int_{\R} ( \Xi_{A,B} z_2 ) X_\vartheta U \left [ X_\theta^2 \, , \frac{a_\omega^+}{\omega} \right ] M_- S^2 v_2 + \int_{\R} ( \Xi_{A,B} z_1 ) X_\vartheta U M_+ \left [ X_\theta^2 \, , \frac{a_\omega^-}{\omega} \right ] S^2 L_+ v_1 , \\
\\ \textbf{K}_4 &=& -\int_{\R} ( \Xi_{A,B} z_2) X_\vartheta U X_\theta^2 n_1 + \int_{\R} ( \Xi_{A,B} z_1 ) X_\vartheta UM_+ X_\theta^2 n_2, \\
\\ \textbf{K}_5 &=& \dot{\omega} \int_{\R} ( \Xi_{A,B} z_2 ) X_\vartheta P_+ w_2 + \dot{\omega} \int_{\R} ( \Xi_{A,B} z_1 ) X_\vartheta P_- w_1 , \\
\\ \textbf{L}_1 &=& \int_{\R} \rho^2 (z_2^2 - z_1 Kz_1 ), \\
\\ \textbf{L}_2 &=& - \int_{\R} \rho^2 z_1 [ X_\vartheta \, , K] Uw_2 , \\
\\ \textbf{L}_3 &=& - \int_{\R} \rho^2 z_2 X_\vartheta U \left [ X_\theta^2 \, , \frac{a_\omega^+}{\omega} \right ] M_- S^2 v_2 - \int_{\R} \rho^2 z_1 X_\vartheta UM_+ \left [ X_\theta^2 \, , \frac{a_\omega^-}{\omega} \right ] S^2 L_+ v_1 , \\
\\ \textbf{L}_4 &=& - \int_{\R} \rho^2 z_2 X_\vartheta U X_\theta^2 n_1 - \int_{\R} \rho^2 z_1 X_\vartheta U M_+ X_\theta^2 n_2 , \\
\\ \textbf{L}_5 &=& \dot{\omega} \int_{\R} \rho^2 z_2 X_\vartheta P_+ w_2 - \dot{\omega} \int_{\R} \rho^2 z_1 X_\vartheta P_- w_1 , \\
\\ \textbf{P} &=& \int_{\R} \left ( 4 ( \partial_y^2 \widetilde{z}_1 )^2 + (4 + \xi_B ) ( \partial_y \widetilde{z}_1)^2 + Y_0 \widetilde{z}_1^2 \right ) , \\
\\ \textbf{R}_1 &=& 4 \int_{\R} ( \chi_A ')^2 \Phi_B ( \partial_y^2 z_1)^2 -4 \int_{\R} ((\chi_A \zeta_B)'''' - \chi_A \zeta_B '''') \chi_A \zeta_B z_1^2 \\
\\ & & \, \, \, \, \, + 8 \int_{\R} \left ( 2 ((\chi_A \zeta_B)'' - \chi_A \zeta_B '') \chi_A \zeta_B - (((\chi_A \zeta_B)')^2 - \chi_A^2 (\zeta_B')^2) \right ) ( \partial_y z_1)^2 , \\
\\ \textbf{R}_2 &=& -3 \int_{\R} \left ( 3 ( \chi_A^2)' ( \zeta_B^2)' + 3 ( \chi_A^2)'' \zeta_B^2 + ( \chi_A^2 ) ''' \Phi_B \right ) ( \partial_y z_1 )^2 , \\
\\ \textbf{R}_3 &=& \frac{1}{2} \int_{\R} ( \Psi_{A,B}''''' - \chi_A^2 ( \zeta_B^2)'''' ) z_1^2 , \\
\\ \textbf{R}_4 &=& 4 \int_{\R} ( \chi_A ')^2 \Phi_B ( \partial_y z_1)^2 -2 \int_{\R} ( \chi_A^2 )' \Phi_B K_2 ( \partial_y z_1)^2 , \\
\\ \textbf{R}_5 &=& - \int_{\R} ( \Psi_{A,B}''' - \chi_A^2 ( \zeta_B^2 ) '' ) z_1^2 + \frac{1}{2} \int_{\R} ( \Psi_{A,B} ''' - \chi_A^2 ( \zeta_B^2) '' ) K_2 z_1^2 \\
\\ & & \, \, \, \, \, + \int_{\R} \left ( 2 ( \chi_A')^2 \zeta_B^2 + ( \chi_A^2 )'' \Phi_B \right ) K_2 z_1^2 + \frac{1}{2} \int_{\R} (\chi_A^2)' \Phi_B K_2'' z_1^2, \end{array} \]
\[  \begin{array}{ccl} \textbf{R}_6 &=& - \frac{1}{2} \int_{\R} \left ( 2 ( \chi_A ')^2 \zeta_B^2 + ( \chi_A^2 )'' \Phi_B \right ) K_1 z_1^2 - \frac{1}{2} \int_{\R} (\chi_A^2)' \Phi_B K_1^2 z_1^2 , \\
\\ \textbf{R}_7 &=& 4 \int_{\R} \chi_a \zeta_B ( \chi_A '' \beta_B + 2 \chi_A ' \zeta_B ') z_1^2 + \int_{\R} \chi_A \zeta_B ( \chi_A '' \zeta_B \xi_B + 2 \chi_A ' \zeta_B ' \xi_B + \chi_A ' \zeta_B \xi_B ' ) z_1^2 , \\
\\ \textbf{R}_8 &=& \int_{\R} \chi_A^2 ( y \zeta_B^2 - \Phi_B ) K_0 ' z_1^2 + \frac{1}{2} \int_{\R} \chi_A^2 \left ( (\zeta_B^2)'' K_2 + ( \zeta_B')^2 (2 K_2' - K_1) \right ) z_1^2 , \\
\\ \textbf{R}_9 &=& \int_{\R} \chi_A^2 \left ( 2 \zeta_B '' \zeta_B - 2 ( \zeta_B')^2 - 3 \zeta_B '''' \zeta_B + 4 \zeta_B ''' \zeta_B + 3 ( \zeta_B '')^2 + \zeta_B \zeta_B '' \xi_B + \zeta_B \zeta_B ' \xi_B ' \right ) z_1^2, \end{array} \]
with $\xi_B = 10 \frac{\zeta_B ''}{\zeta_B} - 14 \frac{( \zeta_B ')^2}{\zeta_B^2} -2 K_2 - \frac{\Phi_B}{\zeta_B^2} K_2 ' + 2 \frac{\Phi_B}{\zeta_B^2} K_1$. \\
\\ We begin by applying Lemma 5 from \cite{Ma0} with $c=1$ and $Y = Y_0 / C \varepsilon_{3 \omega_0 / 2}$, then with $c = \alpha ( \omega_0 ) /10$ and $Y = e^{-|y|}$. Recalling the crucial fact that $\int_{\R} Y_0 \geqslant \frac{I_\omega}{2} \geqslant C \alpha ( \omega_0 ) > 0$ (from Lemma 5), we obtain
\[ \alpha ( \omega_0 )^2 \int_{\R} \rho h^2 \leqslant \alpha ( \omega_0 )^2 \int_{\R} e^{- \frac{\alpha ( \omega_0 )}{10} |y|} h^2 \leqslant C \alpha ( \omega_0 ) \int_{\R} e^{-|y|} h^2 + C \int_{\R} (h')^2 \leqslant C \int_{\R} Y_0 h^2 + \int_{\R} (h')^2 \]
for any $h \in H^1 ( \R )$. This leads to
\[ \textbf{P} \geqslant C \left ( ||\partial_y^2 \widetilde{z}_1 ||^2 + || \partial_y \widetilde{z}_1 ||^2 + \alpha ( \omega_0 )^2 || \sqrt{\rho} \, \widetilde{z}_1 ||^2 \right ) \]
and then, using Lemma 24, 
\[ \alpha ( \omega_0 )^2 \left ( || \rho \partial_y^2 z_1 ||^2 + || \rho \partial_y z_1 ||^2 + || \rho z_1 ||^2 \right ) \leqslant C \textbf{P} + C A^{-4} \theta^{-5} \left ( ||\eta_A \partial_y v ||^2 + || \eta_A v ||^2 \right ). \]
We estimate the terms $\textbf{R}_j$ for $1 \leqslant j \leqslant 9$ as follows:
\[ \sum_{k=1}^7 | \textbf{R}_j | \leqslant \frac{CB}{A \theta^5} \left ( || \eta_A \partial_y v ||^2 + \frac{B^2}{A^2} || \eta_A v ||^2 \right ) , \, \, \, \, \, \text{and} \, \, \, \, \, | \textbf{R}_8| + | \textbf{R}_9 | \leqslant \frac{C}{B} || \nu z_1 ||^2 . \]
We ultimately find that
\[ C_1 \alpha ( \omega_0 )^2 \left ( || \rho \partial_y^2 z_1 ||^2 + || \rho \partial_y z_1 ||^2 + || \rho z_1 ||^2 \right ) \leqslant \textbf{K}_1 + \frac{CB}{A \theta^5} \left ( || \eta_A \partial_y v ||^2 + \frac{B^2}{A^2} || \eta_A v ||^2 \right ). \]
We estimate $\textbf{L}_1$ as follows: $\textbf{L}_1 \geqslant ||\rho z_2||^2 - C \left ( || \rho \partial_y^2 z_1||^2 + || \rho \partial_y z_1||^2 + || \rho z_1 ||^2 \right )$. Setting $\mathbf{Z} := || \rho \partial_y^2 z_1 ||^2 + || \rho \partial_y z_1 ||^2 + || \rho z_1 ||^2 + || \rho z_2 ||^2$, we get that
\[ \alpha ( \omega_0 )^2 \textbf{Z} \leqslant C \left [ \textbf{K}_1 + \alpha ( \omega_0 )^2 \textbf{L}_1 + \frac{B}{A \theta^5} \left ( || \eta_A \partial_y v ||^2 + \frac{B^2}{A^2} || \eta_A v ||^2 \right ) \right ]. \]
Now, let us control the other $\textbf{K}_j$ and $\textbf{L}_j$. We estimate $\textbf{K}_2$ and $\textbf{L}_2$ as follows: taking $\vartheta = \theta^{1/4}$, we have $| \textbf{K}_2 | \leqslant CB \theta^{1/8} \textbf{Z}$ and $| \textbf{L}_2 | \leqslant C \theta^{1/8} \textbf{Z}$. \\
\\ As for $\textbf{K}_3$ and $\textbf{L}_3$, write that $| \xi_W | \leqslant C \varepsilon_{3 \omega_0 / 2}$ and $\left | \frac{a_\omega^+}{\omega} \right | \leqslant \frac{C \varepsilon_{3 \omega_0 /2}}{\omega_0} e^{-2|y|}$. We find that
\[ \begin{array}{rl} & | \textbf{K}_3 | \leqslant CB \theta^{1/4} \left ( 1 + \frac{\varepsilon_{3 \omega_0 /2}}{\omega_0} \right ) \alpha ( \omega_0 )^{-3/2} \textbf{Z} \leqslant C \theta^{1/8} \textbf{Z} \\
\\ \text{and} & | \textbf{L}_3 | \leqslant C\theta^{1/4} \left ( 1 + \frac{\varepsilon_{3 \omega_0 /2}}{\omega_0} \right ) \alpha ( \omega_0 )^{-3/2} \textbf{Z} \leqslant C \theta^{1/8} \textbf{Z}, \end{array} \]
taking $\theta > 0$ small enough (depending on $B$ and $\omega_0$). It is for these estimates that we need the entire hypothesis $(H_1)$: here $g$ must be differentiated $5$ times, and the assumption that $s^4g^{(5)} (s)$ is bounded is enough. This leads to
\begin{equation}
    \alpha ( \omega_0 )^2 \textbf{Z} \leqslant C \left [ \textbf{K}_1 + \textbf{K}_2 + \textbf{K}_3 + \alpha ( \omega_0 )^2 ( \textbf{L}_1 + \textbf{L}_2 + \textbf{L}_3 ) + \frac{B}{A \theta^5} \left ( ||\eta_A \partial_y v ||^2 + \frac{B^2}{A^2} || \eta_A v ||^2 \right ) \right ].
    \label{prop6A}
\end{equation}
Now, as for $\textbf{K}_4$ and $\textbf{L}_4$, we write $q_1 = q_{1,1} + q_{1,2}$ and $q_2 = q_{2,1} + q_{2,2}$ where $q_{1,1} = b_1^2 G + b_2^2 H$, $q_{2,1} = b_1 b_2 G_2$,
\[ \begin{array}{rl} & q_{1,2} = Q_\omega ( 3 + 3 g'(\omega Q_\omega^2 ) + 2 \omega Q_\omega^2 g''( \omega Q_\omega^2 )) (2 b_1 V_1 v_1 + v_1^2) + Q_\omega ( 1 + g'( \omega Q_\omega^2 )) (2b_2 V_2 v_2 + v_2^2 ) + N_1 \\
\\ \text{and} & q_{2,2} = 2 Q_\omega ( 1 + g'( \omega Q_\omega^2 )) (2b_1 V_1 v_2 + 2b_2 V_2 v_1 + v_1v_2) + N_2, \end{array} \]
with $|N_1| + |N_2| \leqslant C |u|^{7/3} \leqslant C |b|^{7/3} \rho^{16} + C |v|^{7/3}$. We define $n_{1,1} = S^2 L_+ q_{2,1}^\perp$, $n_{2,1} = M_- S^2 q_{1,1}^\T$, $n_{1,2} = -S^2 L_+ p_2^\perp + S^2 L_+ q_{2,2}^\perp + S^2 L_+ r_2^\perp + \dot{\omega}Q_+ v_1$ and $n_{2,2} = - M_- S^2 p_1^\T + M_- S^2 q_{1,2}^\T + M_- S^2 r_1^\T + \dot{\omega} Q_- v_2$. \\
\\ We successively prove that $|n_{1,1}^{(k)}| + |n_{2,1}^{(k)}| \leqslant C ( \nu + \sqrt{\alpha ( \omega_0 )} \, \rho^8 )$,
\[ \left | \int_{\R} ( \Xi_{A,B} z_2 ) X_\vartheta U X_\theta^2 n_{1,1} \right | \leqslant C B |b|^2 || \rho z_2 || \, \, \, \, \, \text{and} \, \, \, \, \, \left | \int_{\R} ( \Xi_{A,B} z_1 ) X_\vartheta U M_+ X_\theta^2 n_{2,1} \right | \leqslant C B |b|^2 || \rho z_1 ||. \]
Taking $A$ large enough (depending on $\omega_0$), we have $|| \eta_A p_2^\perp || \leqslant C A \alpha ( \omega_0 )^{-1/2} \epsilon ( || \nu v ||^2 + |b|^2 ) \leqslant C A^2 \epsilon ( || \nu v ||^2 + |b|^2 )$ and $|| \eta_A q_{2,2}^\perp || \leqslant C \alpha ( \omega_0 )^{-1/2} \epsilon^{1/3} ( || \eta_A v || + |b|^2 ) \leqslant CA \epsilon^{1/3} ( || \eta_A v || + |b|^2 )$. Moreover, using \eqref{dVomega}, we have
\[ || \eta_A r_2^\perp || \leqslant C \underline{\mathbf{V}} ( \omega_0 ) \alpha ( \omega_0 )^{-1/2} \epsilon ( ||\nu v ||^2 + |b|^2 ) \leqslant C A \epsilon ( ||\nu v ||^2 + |b|^2 ), \]
taking $A$ large enough (depending on $\omega_0$). The rest of the proof is unchanged and we eventually find that
\[ \begin{array}{rl} & | \textbf{K}_4 | \leqslant CB|b|^2 || \rho z || + C A^2 B \theta^{-9/4} \epsilon^{1/3} \left ( || \eta_A \partial_y^2 z_1 || + ||\eta_A \partial_y z_1 || + || \eta_A z_1 || + || \eta_A z_2 || \right ) \left ( || \eta_A v || + |b|^2 \right ) \\
\\ \text{and} & | \textbf{L}_4 | \leqslant C|b|^2 || \rho z || + C A^2 \theta^{-9/4} \epsilon^{1/3} \left ( || \eta_A \partial_y^2 z_1 || + ||\eta_A \partial_y z_1 || + || \eta_A z_1 || + || \eta_A z_2 || \right ) \left ( || \eta_A v || + |b|^2 \right ). \end{array} \]
As for $\textbf{K}_5$ and $\textbf{L}_5$, the proof relies on Lemma 22. We find that
\[ \begin{array}{rl} & | \textbf{K}_5 | \leqslant CB \theta^{-9/4} \underline{\mathbf{P}} ( \omega_0 ) \epsilon ( || \nu v ||^2 + |b|^2 ) \left ( || \eta_A \partial_y^2 z_1 || + || \eta_A \partial_y z_1 || + || \eta_A z_1 || + || \eta_A z_2 || \right ) \\
\\ \text{and} & | \textbf{L}_5 | \leqslant C \theta^{-9/4} \underline{\mathbf{P}} ( \omega_0 ) \epsilon ( || \nu v ||^2 + |b|^2 ) \left ( || \eta_A \partial_y^2 z_1 || + || \eta_A \partial_y z_1 || + || \eta_A z_1 || + || \eta_A z_2 || \right ). \end{array} \]
Gathering these last estimates and using Lemma 20, we have
\begin{equation}
     | \textbf{K}_4 | + | \textbf{K}_5 | + | \textbf{L}_4 | + | \textbf{L}_5 | \leqslant CB |b|^2 \textbf{Z}^{1/2} + CA^2 B \theta^{-9/4} \epsilon^{1/3} ( 1 + \underline{\mathbf{P}} ( \omega_0)) ( || \eta_A v || + |b|^2 ) ( || \eta_A \partial_y v || + ||\eta_A v || ).
    \label{prop6B}
\end{equation}
Combining \eqref{prop6A} and \eqref{prop6B}, and taking $\epsilon>0$ small enough (depending on $\theta$ and $A$), we obtain:
\[ \alpha ( \omega_0 )^2 \textbf{Z} \leqslant C \dot{\textbf{K}} + C \alpha ( \omega_0 )^2 \dot{\textbf{L}} + \frac{CB}{A \theta^5} \left ( || \eta_A \partial_y v ||^2 + \frac{B^2}{A^2} || \eta_A v ||^2 \right ) + \frac{CB^2}{\alpha ( \omega_0 )^2} |b|^4. \]
We integrate this inequality on $[0 \, , s]$ and recall that $| \textbf{K} | + | \textbf{L} | \leqslant C \epsilon$. Using Proposition 3 and Proposition 4, it leads to:
\[ \begin{array}{rcl} \displaystyle{\int_0^s \textbf{Z}} & \leqslant & \displaystyle{\frac{C}{\alpha ( \omega_0 )^2} \left ( \epsilon + \frac{B^3}{A \theta^5} \int_0^s \left ( || \eta_A \partial_y v ||^2 + \frac{1}{A^2} || \eta_A v ||^2 \right ) + \frac{B^2}{\alpha ( \omega_0 )^2} \int_0^s |b|^4 \right )} \\
\\ & \leqslant & \displaystyle{C \left ( \frac{1}{\alpha ( \omega_0 )^2} + \frac{B^3}{A \theta^5 \alpha ( \omega_0 )^2} + \frac{B^2}{\alpha ( \omega_0 )^4} \right ) \epsilon + \frac{C}{A} \left ( \frac{B^3}{\alpha ( \omega_0 )^2 \theta^5} + \frac{B^2}{\alpha ( \omega_0 )^5 \underline{\mathbf{\Gamma}} ( \omega_0 )} + \frac{B^3}{A \alpha ( \omega_0 )^3 \theta^5 \underline{\mathbf{\Gamma}} ( \omega_0 )} \right ) \int_0^s || \rho^4 v ||^2.} \end{array} \]
Taking $A>0$ large enough (depending on $\omega_0$, $B$ and $\theta$) in order to control the second term, then $\epsilon >0$ small enough (depending on $\omega_0$, $B$, $\theta$ and $A$) in order to control the first term, we obtain
\[ \int_0^s \textbf{Z} \leqslant C \sqrt{\epsilon} + \frac{C}{\sqrt{A}} \int_0^s || \rho^4 v ||^2, \]
which is the announced result. \hfill \qedsymbol

\section{Final estimates and convergence of $\omega$}
\noindent We finish the proof of Theorem 1 as in \cite{Ma0} (section 11). Recall that
\[ \begin{array}{l} \displaystyle{u(s \, , y) = \frac{e^{-i \gamma (s)}}{\sqrt{\omega (s)}} \psi \left ( \tau(s) \, , \frac{y}{\sqrt{\omega (s)}} \right ) - Q_{\omega (s)},} \\
\\ \begin{array}{ll} u_1 = v_1 + b_1 (t) V_1 , & u_2 = v_2 + b_2(t) V_2 , \\
\\ w_1 = X_\theta^2 M_- S^2 v_2 , & w_2 = - X_\theta^2 S^2 L_+ v_1 , \\
\\ z_1 = X_\vartheta U w_2 , & z_2 = - X_\vartheta U M_+ w_1 . \end{array} \end{array} \]
Proposition 3 gives a large-scale estimate on $v$. Proposition 4 gives a control on $|b|^4$. Proposition 5 gives a very good control on $z$, obtained thanks to the coercivity of $K$, which is the reason why we transformed the system twice (factorizing $(L_+ \, , L_-) \longrightarrow (M_+ \, , M) \longrightarrow K$). Lemma 26 enables to control $v$ once we control $z$. It is time to gather these results to conclude. \\
\\ We combine Lemma 26 and Proposition 5 to get that, for all $s \geqslant 0$,
\[ \alpha ( \omega_0 )^3 \int_0^s ||\rho^4 v ||^2 \leqslant C \int_0^s \left ( || \rho \partial_y^2 z_1||^2 + ||\rho \partial_y z_1||^2 + || \rho z_1||^2 + ||\rho z_2||^2 \right ) \leqslant C \sqrt{\epsilon} + \frac{C}{\sqrt{A}} \int_0^s ||\rho^4 v ||^2. \]
Taking $A$ large enough (depending on $\omega_0$) and then $\epsilon >0$ small enough (depending on $\omega_0$), we obtain
\[ \int_0^s || \rho^4 v ||^2 \leqslant C \alpha ( \omega_0 )^{-3} \sqrt{\epsilon} \leqslant \epsilon^{1/4} \leqslant 1. \]
Passing to the limit $s \to + \infty$ in Propositions 3 and 4 and taking $A>0$ large enough (depending on $\omega_0$), we have successively $\int_0^{+ \infty} |b|^4 \leqslant C$ and
\[ \begin{array}{rcl} \displaystyle{\int_0^{+ \infty} \left ( |b|^4 + || \rho \partial_y v ||^2 + || \rho v ||^2 \right )} & \leqslant & \displaystyle{C \left ( \epsilon + \frac{1}{A \alpha ( \omega_0 ) \underline{\mathbf{\Gamma}} ( \omega_0 )} \int_0^{+ \infty} ||\rho^4 v ||^2 \right ) + CA^2 \left ( \epsilon + \int_0^{+ \infty} || \rho^4 v ||^2 + \int_0^{+ \infty} |b|^4 \right )} \\
\\ & \leqslant & C + CA^2 \, \, \, \leqslant \, \, \, CA^2. \end{array} \]
In particular, there exists a sequence $s_n \to + \infty$ such that
\[ |b(s_n)|^4 + || \rho \partial_y v(s_n) ||^2 + || \rho v(s_n)||^2 \, \, \underset{n \to + \infty}{\longrightarrow} \, \, 0. \]
Recall that, setting $\mathscr{M} = |b|^4 + || \rho v ||^2$, Lemma 13 states that $| \dot{\mathscr{M}} | \leqslant C \left ( |b|^4 + || \rho \partial_y v ||^2 + || \rho v ||^2 \right )$. For $s>0$ and $n$ such that $s_n>s$, we integrate on $(s \, , s_n)$ to find that
\[ \mathscr{M} (s) \leqslant \mathscr{M} (s_n) + \int_s^{s_n} | \dot{\mathscr{M}} | \leqslant \mathscr{M} (s_n) + C \int_s^{s_n} \left ( |b|^4 + || \rho \partial_y v ||^2 + || \rho v ||^2 \right ). \]
Passing to the limit $n \to + \infty$, we find $\displaystyle{\mathscr{M} (s) \leqslant C \int_s^{+ \infty} \left ( |b|^4 + || \rho \partial_y v ||^2 + || \rho v ||^2 \right )}$. Therefore, 
\[ \mathscr{M} (s) \, \, \underset{s \to + \infty}{\longrightarrow} \, \, 0. \]
All that is left to prove is that $\omega (s)$ converges as $s \to + \infty$. We follow ideas from \cite{Cu3}, as the method in \cite{Ma0} does not seem to apply to the framework of the present paper. First consider
\[ \mathbf{M} (s) := \sqrt{\omega (s)} \int_{\R} \frac{|v(s)|^2}{2}. \]
From \eqref{sysv} we have
\[ \frac{\text{d} \mathbf{M}}{\text{d}s} = \frac{\dot{\omega}}{4 \sqrt{\omega}} \int_{\R} |v|^2 + \sqrt{\omega} \int_{\R} \left ( v_1 L_- v_2 + v_1 \mu_2 + v_1 p_2^\perp - v_1 q_2^\perp - v_1 r_2^\perp - v_2 L_+ v_1 - v_2 \mu_1 + v_2 p_1^\top + v_2 q_1^\top + v_2 r_1^\top \right ). \]
First, integrating by parts, $\int_{\R} (v_1 L_- v_2 - v_2 L_+ v_1) = \int_{\R} v_1 v_2 (2Q_\omega^2 + 2Q_\omega^2 g'(\omega Q_\omega^2))$. Thus $\left | \int_{\R} (v_1 L_- v_2 - v_2 L_+ v_1) \right | \leqslant C || \rho v ||_{L^2}^2$. Note that the constants here (and in the following computations) can depend on $\omega_0$: it does not matter, $\omega_0$ is considered as fixed. Now, recalling the definitions of $\mu_1$ and $\mu_2$, we have
\[ \left | \int_{\R} v_1 \mu_2 \right | \leqslant |m_\omega| \int_{\R} |v \Lambda_\omega Q_\omega | \leqslant C ( || \nu v ||^2 + |b|^2 ) \int_{\R} | \rho v | \cdot \rho \leqslant C (  || \rho v ||^2 + |b|^4 ). \]
A similar estimate holds for $\int_{\R} v_2 \mu_1$. Then, recalling the definitions of $p_1$ and $p_2$ we have
\[ \int_{\R} (v_1 p_2^\perp - v_2 p_1^\top ) = - m_\omega \int_{\R} (u_1 \Lambda u_1 + u_2 \Lambda u_2 ) + \sum\limits_{j=1}^2 (-1)^j \left [ b_j \int_{\R} V_j \left ( p_{3-j} - \frac{\langle p_{3-j} \, , V_j \rangle}{\langle V_1 \, , V_2 \rangle} V_{3-j} \right ) + \frac{\langle p_{3-j} \, , V_j \rangle}{\langle V_1 \, , V_2 \rangle} \int_{\R} u_j V_{3-j} \right ]. \]
First, $m_\omega \int_{\R} \left ( u_1 \Lambda u_1 + u_2 \Lambda u_2 \right ) = \frac{\dot{\omega}}{4 \omega} \int_{\R} |u|^2$. Then, $\langle p_2 \, , V_1 \rangle = m_\gamma b_2 \langle V_1 \, , V_2 \rangle + \frac{m_\omega}{2} \langle v_1 \, , y \partial_y V_1 \rangle + \frac{m_\omega}{2} b_1 \langle V_1 \, , y \partial_y V_1 \rangle$ hence, thanks to the estimates on $V$ ($|V| \leqslant C \rho^8$ for instance), $| \langle p_2 \, , V_1 \rangle | \leqslant C \left ( || \nu v ||^2 + |b|^2 \right ) \left ( |b| + || \rho v || \right )$. Writing $u_1 = v_1 + b_1 V_1$, it follows that
\[ \left | \frac{\langle p_2 \, , V_1 \rangle}{\langle V_1 \, , V_2 \rangle} \int_{\R} u_1 V_2 \right | \leqslant C \left ( || \nu v ||^2 + |b|^2 \right ) \left ( |b| + || \rho v || \right ) \left ( C || \rho v || + C |b|^2 \right ) \leqslant C \left ( || \rho v ||^2 + |b|^4 \right ). \]
Similarly, $\left | b_1 \int_{\R} V_1 \left ( p_2 - \frac{\langle p_2 \, , V_1 \rangle}{\langle V_1 \, , V_2 \rangle} V_2 \right ) \right | \leqslant C |b| (|m_\gamma| + |m_\omega| + || \rho v ||^2 + |b|^2) ( || \rho v || + |b| ) \leqslant C ( || \rho v ||^2 + |b|^4 )$. Same estimates hold for $j=2$: we find that
\[ \left | \int_{\R} (v_1 p_2^\perp - v_2 p_1^\top ) + \frac{\dot{\omega}}{4 \omega} \int_{\R} |u|^2 \right | \leqslant C \left ( || \rho v ||^2 + |b|^4 \right ). \]
Writing $u_j = v_j + b_j V_j$, we have
\[ \begin{array}{rcl} \displaystyle{\left | \frac{\dot{\omega}}{4 \omega} \int_{\R} |u|^2 - \frac{\dot{\omega}}{4 \omega} \int_{\R} |v|^2 \right |} &=& \displaystyle{\frac{|m_\omega|}{4} \left | \int_{\R} \left ( 2v_1b_1V_1 + 2v_2b_2V_2 + b_1^2 V_1^2 + b_2^2V_2^2 \right ) \right |} \\ \\ & \leqslant & C ( || \rho v ||^2 + |b|^2 ) ( |b| \, || \rho v || + |b|^2 ) \, \, \leqslant \, \, C ( || \rho v ||^2 + |b|^4 ) \end{array} \]
thus $\left | \int_{\R} (v_1 p_2^\perp - v_2 p_1^\top ) + \frac{\dot{\omega}}{4 \omega} \int_{\R} |v|^2 \right | \leqslant C ( || \rho v ||^2 + |b|^4 )$. Now, recalling the definitions of $r_1$ and $r_2$ and the estimates \eqref{dVomega}, we have $| \langle r_1 \, , V_1 \rangle | \leqslant C|b| \, |m_\omega|$ and 
\[ \left | \int_{\R} v_2 r_1^\top \right | \leqslant C |b| \, |m_\omega| \, ||\rho v || \leqslant C ( || \rho v ||^2 + |b|^4 ). \]
A similar estimate holds for $\int_{\R} v_1 r_2^\perp$. At this point we have proven that
\[ \left | \frac{\text{d} \mathbf{M}}{\text{d}s} \right | \leqslant C ( || \rho v ||^2 + |b|^4 ) + \sqrt{\omega} \left | \int_{\R} (v_2 q_1^\top - v_1 q_2^\perp ) \right |. \]
The last term requires a little more care. Set $\tilde{g}_\omega (s) = s + \frac{g( \omega s )}{\omega}$ so that $f_\omega ( \psi ) = \tilde{g}_\omega ( | \psi |^2 ) \psi$. We see that
\[ \begin{array}{rccl} & q_1 &=& (Q_\omega + b_1V_1) \tilde{g}_\omega ((Q_\omega + v_1 + b_1V_1)^2 + (v_2+b_2V_2)^2) - (v_1+b_1V_1)( \tilde{g}_\omega (Q_\omega^2) + 2 Q_\omega^2 \tilde{g}_\omega ' ( Q_\omega^2 )) \\ & & & \, \, \, \, - Q_\omega \tilde{g}_\omega (Q_\omega^2) + v_1 \tilde{g}_\omega ((Q_\omega + v_1+b_1V_1)^2 + (v_2+b_2V_2)^2) \\
\\ \text{and} & q_2 &=& b_2V_2 \tilde{g}_\omega ((Q_\omega + v_1+b_1V_1)^2 + (v_2+b_2V_2)^2) - \tilde{g}_\omega (Q_\omega^2)(v_2+b_2V_2) \\ & & & \, \, \, \, + v_2 \tilde{g}_\omega ((Q_\omega + v_1+b_1V_1)^2 + (v_2 + b_2 V_2)^2). \end{array} \]
This leads to (note the cancellation of the last terms in $v_2q_1$ and in $v_1q_2$):
\[ v_2q_1 - v_1q_2 = ( Q_\omega v_2 + v_2 b_1 V_1 - v_1b_2V_2) \left [ \tilde{g}_\omega ((Q_\omega + v_1+b_1V_1)^2 + (v_2 + b_2V_2)^2) - \tilde{g}_\omega (Q_\omega^2) \right ]  - 2 v_2 (v_1 + b_1V_1) Q_\omega^2 \tilde{g}_\omega ' (Q_\omega^2). \]
$\star$ First case: assume $\frac{1}{100} Q_\omega \leqslant |v| + |b_1V_1| + |b_2V_2|$. In this case we simply write 
\[ | \tilde{g}_\omega ((Q_\omega + v_1+b_1V_1)^2 + (v_2+b_2V_2)^2) - \tilde{g}_\omega (Q_\omega^2)| \leqslant C (|v|^2 + |b_1V_1|^2 + |b_2V_2|^2 ), \]
which leads to
\[ | v_2q_1 - v_1q_2 | \leqslant C (Q_\omega |v| + |v| (|b_1V_1|+|b_2V_2|)) (|v|^2 + |b_1V_1|^2 + |b_2V_2|^2) + C|v| (|v| + |b_1V_1|) Q_\omega^2 \leqslant C ( | \rho v |^2 + |b|^4 \rho^8 ), \]
using estimates on $V$ and $Q_\omega$ and the fact that $Q_\omega^2 \leqslant C Q_\omega ( |v| + |b_1V_1| + |b_2V_2| ) \leqslant C ( |v| + |b|) \rho^8$ in this situation. \\
\\ $\star$ Second case: assume $|v| + |b_1V_1| + |b_2V_2| \leqslant \frac{1}{100} Q_\omega$. Here we write Taylor's expansion
\[ \begin{array}{rl} & \tilde{g}_\omega ((Q_\omega + v_1+b_1V_1)^2 + (v_2+b_2V_2)^2) \\ =& \tilde{g}_\omega (Q_\omega^2) +  \tilde{g}_\omega ' (Q_\omega^2) \left ( v_1^2 + b_1^2V_1^2  + 2 Q_\omega v_1 + 2v_1b_1V_1  + 2 Q_\omega b_1V_1 + v_2^2 + b_2^2 V_2^2 + 2v_2b_2V_2 \right ) + \mathbf{IR} \end{array} \]
where $\mathbf{IR} = \int_{Q_\omega^2}^{(Q_\omega + v_1 + b_1V_1)^2  (v_2+b_2V_2)^2} \left ( (Q_\omega + v_1 + b_1V_1 )^2 + (v_2 + b_2V_2)^2 - z \right ) \tilde{g}_\omega '' (z) \, \text{d}z$. It follows from a convenient cancellation that
\[ \begin{array}{rcl} v_2q_1 - v_1q_2 &=& Q_\omega v_2 (v_1^2 + b_1^2V_1^2 + 2v_1b_1V_1 + v_2^2 + b_2^2V_2^2 + 2v_2b_2V_2) \tilde{g}_\omega ' (Q_\omega^2) \\ & & \, \, \, \, \, + (v_2b_1V_1 - v_1b_2V_2)(v_1^2 + b_1^2 V_1^1 + 2 Q_\omega v_1 + 2v_1b_1V_1 + 2 Q_\omega b_1 V_1 + v_2^2 + b_2^2V_2^2 + 2v_2b_2V_2) \tilde{g}_\omega ' (Q_\omega^2) \\ & & \, \, \, \, \, + (Q_\omega v_2 + v_2b_1V_1 - v_1b_2V_2) \mathbf{IR}. \end{array} \]
Since $|v| + |b_1V_1| + |b_2V_2| \leqslant \frac{1}{100} Q_\omega$, we have $| \tilde{g}_\omega ''(z) | \leqslant \frac{C}{z} \leqslant \frac{C}{Q_\omega^2}$ for any $z \in [Q_\omega^2 \, , (Q_\omega + v_1 + b_1V_1)^2 + (v_2 + b_2V_2)^2 ]$. This leads to
\[ | \mathbf{IR} | \leqslant \frac{C}{Q_\omega^2} \left |(Q_\omega + v_1 + b_1V_1)^2+(v_2+b_2V_2)^2-Q_\omega^2 \right |^2 \leqslant \frac{C}{Q_\omega^2} \left [ Q_\omega (|v| + |b_1V_1| + |b_2V_2|) \right ]^2 \leqslant C ( |v| + |b_1V_1| + |b_2V_2| )^2. \]
Using again our hypothesis $|v| + |b_1V_1| + |b_2V_2| \leqslant Q_\omega$, we evidently find
\[ \begin{array}{rcl} |v_2q_1 - v_1q_2| & \leqslant & C Q_\omega |v| (|v|^2 + |b_1V_1|^2 + |b_2V_2|^2 ) + C(|v|^2 + |b_1V_1|^2 + |b_2V_2|^2 ) Q_\omega (|v| + |b_1V_1| + |b_2V_2|) \\ & & \, \, \, \, \, + C Q_\omega |v| ( |v|^2 + |b_1V_1|^2 + |b_2V_2|^2 ) \\ \\ & \leqslant & C (|\rho v|^2 + |b|^4 \rho^8 ). \end{array} \]
Gathering both cases, we see that the estimate $|v_2 q_1 - v_1 q_2| \leqslant C (| \rho v |^2 + |b|^4 \rho^8 )$ always holds, hence 
\[ \int_{\R} |v_2 q_1 - v_1 q_2 | \leqslant C ( || \rho v ||^2 + |b|^4 ). \]
The term in $\dot{\mathbf{M}}$ is actually 
\[ \int_{\R} (v_2 q_1^\top - v_1 q_2^\perp ) = \int_{\R} (v_2 q_1 - v_1 q_2 )  + \int_{\R} \left ( v_1 \frac{\langle q_2 \, , V_2 \rangle}{\langle V_1 \, , V_2 \rangle} V_1 - v_2 \frac{\langle q_1 \, , V_1 \rangle}{\langle V_1 \, , V_2 \rangle} V_1 \right ). \]
Lemma 9 assures that $|q_1| + |q_2| \leqslant C (Q_\omega |u|^2 + |u|^{7/3} ) \leqslant C |u|^2 \leqslant C (|v|^2 + |b|^2)$. Since $|V| \leqslant \rho^8$, this leads to $| \langle q_2 \, , V_2 \rangle | \leqslant C ( || \rho v ||^2 + |b|^2 )$ then $\left | \int_{\R} v_1 \frac{\langle q_2 \, , V_2 \rangle}{\langle V_1 \, , V_2 \rangle} V_1 \right | \leqslant C ( || \rho v ||^2 + |b|^4 )$. The same estimate holds for $\int_{\R} v_2 \frac{\langle q_1 \, , V_1 \rangle}{\langle V_1 \, , V_2 \rangle} V_2$. Finally, 
\[ \left | \int_{\R} (v_2 q_1^\top - v_1 q_2^\perp ) \right | \leqslant C ( || \rho v ||^2 + |b|^4 ) \, \, \, \, \, \, \, \, \, \, \, \, \text{hence} \, \, \, \, \, \, \, \, \, \, \, \, \left | \frac{\text{d} \mathbf{M}}{\text{d}s} \right | \leqslant C \left ( || \rho v ||^2 + |b|^4 \right ). \]
Since $\int_0^{+ \infty} ( || \rho v  ||^2 + |b|^4 ) \, \text{d}s < + \infty$, the quantity $\mathbf{M} (s)$ converges as $s \to + \infty$. Now, recalling Lemma 10, write that
\[ e^{-i \gamma (s(t))} \psi (t \, , x) = \phi_{\omega (s(t))} + \sqrt{\omega (t(s))} \, v ( \sqrt{\omega (s(t))} \, x \, , s(t)) + \sum\limits_{j=1}^2 \sqrt{\omega (t(s))} \, b_j(s(t)) V_j ( \sqrt{\omega (s(t))} \, x \, , s(t) ). \]
Take the square of the $L^2$ norm on both sides: on the left side, $||e^{-i \gamma} \psi (t)||^2 = ||\psi_0||^2$ thanks to the conservation of the mass. On the right side, we expand and analyse each term. In what follows, take $j \in \{ 1 \, , 2 \}$.
\begin{itemize}
	\item First, $|| \sqrt{\omega} \, v(\sqrt{\omega} \, x \, , s(t)) ||^2 = 2 \mathbf{M} (s(t))$. This quantity converges as $t \to + \infty$.
	\item Then, $|| \sqrt{\omega} \, b_j(s(t)) V_j( \sqrt{\omega} \, x \, , s(t)) ||^2 = |b_j(s(t))|^2 \sqrt{\omega} \int_{\R} |V_j(y \, , s(t))|^2 \, \text{d}y \leqslant C |b(s(t))|^2 \, \underset{t \to + \infty}{\longrightarrow} \, 0$. 
	\item Similarly, we prove that $\int_{\R} \left | \phi_\omega (x) \cdot \sqrt{\omega} \, b_j(s(t)) V_j( \sqrt{\omega} \, x \, , s(t)) \right | \text{d}x  \, \underset{t \to + \infty}{\longrightarrow} \, 0$ thanks to the convergence $|b(s(t))| \, \underset{t \to + \infty}{\longrightarrow} \, 0$ and the bound on $\phi_{\omega}$.
	\item Similarly, we prove that $\int_{\R} \left | \sqrt{\omega} \, b_1(s(t)) V_1( \sqrt{\omega} \, x \, , s(t)) \cdot \sqrt{\omega} \, b_2(s(t)) V_2( \sqrt{\omega} \, x \, , s(t)) \right | \text{d}x  \, \underset{t \to + \infty}{\longrightarrow} \, 0$ thanks to the convergence $|b(s(t))| \, \underset{t \to + \infty}{\longrightarrow} \, 0$.
	\item Similarly, we prove that $\int_{\R} \left | \sqrt{\omega} \, v(\sqrt{\omega} \, x \, , s(t)) \cdot \sqrt{\omega} \, b_j(s(t)) V_j( \sqrt{\omega} \, x \, , s(t)) \right | \text{d}x  \, \underset{t \to + \infty}{\longrightarrow} \, 0$ thanks to the convergence $|b(s(t))| \, \underset{t \to + \infty}{\longrightarrow} \, 0$ and the estimate $||v||_{L^\infty} \leqslant C ||v||_{H^1} \leqslant C \epsilon \leqslant C$. 
	\item Moreover, $\int_{\R} \left | \phi_\omega (x) v ( \sqrt{\omega} \, x \, , s(t)) \right | \text{d}x = \int_{\R} Q_\omega (y) |v(y \, , s(t))| \, \text{d}y \leqslant \int_{\R} \rho^8 |v(s(t))| \leqslant C || \rho v(s(t)) || \, \underset{t \to + \infty}{\longrightarrow} \, 0$ thanks to Cauchy-Schwarz inequality and the convergence $|| \rho v(s(t))|| \, \underset{t \to + \infty}{\longrightarrow} \, 0$. 
\end{itemize}
The only remaining term on the expansion of the right-hand side is $||\phi_{\omega(s(t))}||^2$. What precedes shows that $||\phi_{\omega(s(t))}||^2$ converges as $t \to + \infty$. Recall the well-know property $\partial_\omega ||\phi_\omega||^2 > 0$, thus $\omega \mapsto ||\phi_\omega||^2$ is a strictly monotonic function. As a consequence, from the convergence of $||\phi_{\omega(s(t))}||^2$ we deduce the convergence of $\omega (s(t))$ as $t \to + \infty$. We denote by $\omega_+$ the limit of $\omega (s(t))$ as $t \to + \infty$. \\
\\ The estimate $|m_\gamma| \leqslant C || \nu u ||^2 \, \underset{s \to + \infty}{\longrightarrow} \, 0$ shows that $\dot{\gamma} (s) \, \underset{s \to + \infty}{\longrightarrow} \, 1$. Hence, $\tilde{\gamma} (t) = \gamma (s(t))$ satisfies
\[ \frac{\text{d} \tilde{\gamma}}{\text{d}t} = \frac{\text{d} \gamma}{\text{d}s} \, \frac{\text{d}s}{\text{d}t} = \dot{\gamma} (s(t)) \omega (s(t)) \, \underset{t \to + \infty}{\longrightarrow} \, \omega_+ . \]
With our notation here, the "$\gamma$" presented in the statement of Theorem 2 is actually our $\tilde{\gamma}$. Theorem 2 is now fully proven. \hfill \qedsymbol

\end{document}